\documentclass{article}

\usepackage{microtype}
\usepackage{graphicx,array}
\usepackage{subfigure,tabularx}
\usepackage{booktabs} 
\usepackage{amssymb,amsmath,amsfonts,mathrsfs}
\usepackage{etoolbox,verbatim,color}
\usepackage{url}
\newcolumntype{C}[1]{>{\centering}m{#1}}
\usepackage{hyperref}

\newtheorem{assumption}{Assumption}
\newtheorem{definition}{Definition}
\newtheorem{proposition}{Proposition}
\newtheorem{remark}{Remark}
\newtheorem{theorem}{Theorem}
\newtheorem{lemma}{Lemma}
\def\QEDmark{\ensuremath{\square}}
\def\proof{\paragraph{Proof:}}
\def\endproof{\hfill\QEDmark}
\usepackage{multicol}
\usepackage{multirow}
\catcode`~=11 \def\UrlSpecials{\do\~{\kern -.15em\lower .7ex\hbox{~}\kern .04em}} \catcode`~=13 

\newcommand{\iprod}[2]{\left\langle {#1}, {#2} \right\rangle}
\newcommand{\norm}[1]{\left\lVert {#1} \right\rVert}

\newcommand{\lrpar}[1]{\left(#1\right)}

\newcommand{\lrbrace}[1]{\left\{#1\right\}}

\newcommand{\calA}{\mathcal{A}}

\newcommand{\bbE}{\mathbb{E}}

\newcommand{\bbN}{\mathbb{N}}

\newcommand{\bbR}{\mathbb{R}}

\DeclareMathAlphabet{\mathbsf}{OT1}{cmss}{bx}{n}
\DeclareMathAlphabet{\mathssf}{OT1}{cmss}{m}{sl}

\DeclareSymbolFont{bsfletters}{OT1}{cmss}{bx}{n}  
\DeclareSymbolFont{ssfletters}{OT1}{cmss}{m}{n}
\DeclareMathSymbol{\bsfGamma}{0}{bsfletters}{'000}
\DeclareMathSymbol{\ssfGamma}{0}{ssfletters}{'000}
\DeclareMathSymbol{\bsfDelta}{0}{bsfletters}{'001}
\DeclareMathSymbol{\ssfDelta}{0}{ssfletters}{'001}
\DeclareMathSymbol{\bsfTheta}{0}{bsfletters}{'002}
\DeclareMathSymbol{\ssfTheta}{0}{ssfletters}{'002}
\DeclareMathSymbol{\bsfLambda}{0}{bsfletters}{'003}
\DeclareMathSymbol{\ssfLambda}{0}{ssfletters}{'003}
\DeclareMathSymbol{\bsfXi}{0}{bsfletters}{'004}
\DeclareMathSymbol{\ssfXi}{0}{ssfletters}{'004}
\DeclareMathSymbol{\bsfPi}{0}{bsfletters}{'005}
\DeclareMathSymbol{\ssfPi}{0}{ssfletters}{'005}
\DeclareMathSymbol{\bsfSigma}{0}{bsfletters}{'006}
\DeclareMathSymbol{\ssfSigma}{0}{ssfletters}{'006}
\DeclareMathSymbol{\bsfUpsilon}{0}{bsfletters}{'007}
\DeclareMathSymbol{\ssfUpsilon}{0}{ssfletters}{'007}
\DeclareMathSymbol{\bsfPhi}{0}{bsfletters}{'010}
\DeclareMathSymbol{\ssfPhi}{0}{ssfletters}{'010}
\DeclareMathSymbol{\bsfPsi}{0}{bsfletters}{'011}
\DeclareMathSymbol{\ssfPsi}{0}{ssfletters}{'011}
\DeclareMathSymbol{\bsfOmega}{0}{bsfletters}{'012}
\DeclareMathSymbol{\ssfOmega}{0}{ssfletters}{'012}

\newcommand{\lea}{\stackrel{\rm(a)}{\le}}
\newcommand{\leb}{\stackrel{\rm(b)}{\le}}

\DeclareMathOperator*{\argmin}{argmin} 
\DeclareMathOperator*{\dist}{dist}

\DeclareMathOperator{\minimize}{minimize}

\usepackage[accepted]{icml2020}

\icmltitlerunning{Inertial Block Proximal Method}

\begin{document}

\twocolumn[
\icmltitle{Inertial Block Proximal Methods for \\ 
Non-Convex Non-Smooth Optimization}



\icmlsetsymbol{equal}{*}

\begin{icmlauthorlist}
\icmlauthor{Le Thi Khanh Hien}{to}
\icmlauthor{Nicolas Gillis}{to}
\icmlauthor{Panagiotis Patrinos}{goo}
\end{icmlauthorlist}

\icmlaffiliation{to}{Department of Mathematics and Operational Research, University of Mons, Belgium}
\icmlaffiliation{goo}{Department of Electrical Engineering (ESAT-STADIUS), KU Leuven, Belgium}

\icmlcorrespondingauthor{Le Thi Khanh Hien}{thikhanhhien.le@umons.ac.be}

\icmlkeywords{non-convex optimization, extrapolation,  block coordinate descent, matrix methods, tensor methods}

\vskip 0.2in
]



\printAffiliationsAndNotice{}  

\begin{abstract}
We propose inertial versions of block coordinate descent methods for solving non-convex non-smooth composite optimization problems. 
Our methods possess three main advantages compared to current state-of-the-art accelerated first-order methods: 
(1)~they allow using two different extrapolation points to evaluate the gradients and to add the inertial force (we will empirically show that it is more efficient than using a single extrapolation point), 
(2)~they allow to randomly picking the block of variables to update, 
and 
(3)~they do not require a restarting step. 
We prove the subsequential convergence of the generated sequence under mild assumptions, prove the global convergence under some additional assumptions, and provide convergence rates. We deploy the proposed methods to solve non-negative matrix factorization (NMF) and show that they compete favourably with the state-of-the-art NMF algorithms. Additional experiments on non-negative approximate canonical polyadic decomposition, also known as non-negative tensor factorization, are also provided.  
\end{abstract}
\section{Introduction}
\label{sec:Intro}
In this paper, we consider the following non-smooth non-convex optimization problem 
\begin{equation}
\label{eq:main}
\minimize_{x\in\mathbb{E}}F\left(x\right), \quad \text{ where }F(x) :=  f(x)+ r(x),
\end{equation}
and  
$\bbE =\bbE_1\times \ldots\times \bbE_s$ with $\bbE_i$, $i=1,\ldots,s$, being finite dimensional real linear spaces equipped with  norm $\norm{\cdot}_{(i)}$ and  inner product  $\iprod{\cdot}{\cdot}_{(i)}$,  
 $f:\mathbb{E}\to\mathbb{R}$ is a continuous but possibly non-smooth non-convex function, and 
 $r(x)=\sum_{i=1}^{s}r_{i}(x_{i})$ with $r_{i}:\mathbb{E}_{i}\to\mathbb{R}\cup\left\{ +\infty\right\}$ for $i=1,\ldots,s$ being proper and lower semi-continuous functions.
 
Problem~\eqref{eq:main} covers many applications including 
compressed sensing with non-convex ``norms" \cite{Attouch2010}, 
sparse dictionary learning \cite{aharon2006k,XuYin2016}, 
non-negative matrix factorization (NMF) \cite{Gillis2014}, and 
``$l_p$-norm" regularized sparse regression problems with $0\leq p<1$ \cite{BLUMENSATH2009,Natarajan1995}. 
In this paper, we will focus on NMF which is defined as follows: given  $X\in \bbR_+^{\mathbf m\times \mathbf n}$ and the integer $\mathbf r < \min(\mathbf m, \mathbf n)$, solve 
\begin{equation}
\label{eq:NMF_F}
\min_{U,V} \frac12\norm{X-UV}_F^2 \, \text{such that} \, U\in \bbR_+^{\mathbf m\times \mathbf{r}}, V\in \bbR_+^{\mathbf r\times \mathbf n}. 
\end{equation}
NMF is a key problem in data analysis and machine learning with applications in image processing, document classification, hyperspecral unmixing and audio source separation, to cite a few \cite{cichocki2009nonnegative, Gillis2014, fu2019nonnegative}. 
NMF can be written as a problem of the form~\eqref{eq:main} with $s=2$, letting   $f(U,V)=\frac12\norm{X-UV}_F^2$, and $r_1$ and $r_2$ being indicator functions  $r_1(U)=I_{\bbR_+^{\mathbf m\times \mathbf r}}(U)$, and $r_2(V)=I_{\bbR_+^{\mathbf r\times \mathbf n}}(V)$. Note that $UV=\sum_{i=1}^{r}U_{:i}V_{i:}$; 
hence NMF can also be written as a function of $2\times \mathbf r$ variables $U_{:i}$ (the columns of $U$) and $V_{i:}$ (the rows of $V$) for $i=1,\ldots,\mathbf  r$. 
\vspace{-0.1in}
\subsection{Related works} 
The Gauss-Seidel iteration scheme, also known as the block coordinate descent (BCD) method, is a standard approach to solve both convex and non-convex problems in the form of~\eqref{eq:main}. Starting with a given initial point $x^{(0)}$, the method generates a sequence $\lrbrace{x^{(k)}}_{k\geq 0}$ by cyclically updating one block of variables at a time while fixing the values of the other blocks. Let us denote $f_i^{(k)}(x_i) := f\big(x_1^{(k)},\ldots,x_{i-1}^{(k)},x_i,x_{i+1}^{(k-1)},\ldots,x_s^{(k-1)}\big)$ the value of the objective function for the $i$th block at the $k$th iteration of a BCD method. Based on how the blocks are updated, BCD methods can typically be classified into three categories: 
\vspace{-0.1in}
\begin{enumerate}
\item Classical BCD methods update~\cite{GRIPPO20001,Hildreth} using exact updates:  
\vspace{-0.05in}
\[
x_i^{(k)}=\argmin_{x_i\in \bbE_i} f_i^{(k)}(x_i)+r_i(x_i).  
\] 
\vspace{-0.2in}
\item Proximal BCD methods update  \cite{Auslender1992,GRIPPO20001,Razaviyayn2013} using exact updates along with a proximal term: 
\begin{equation} 
\label{BCD2}
x_i^{(k)}= \argmin_{x_i\in \bbE_i}    f_i^{(k)}\lrpar{x_i}+r_i(x_i) + \frac{1}{2\beta_i^{(k)}} \|x_i-x_i^{(k-1)}\|^2,
\end{equation}
where $\beta_i^{(k)}$ is  referred to as the stepsize. 
\item Proximal gradient BCD methods update \cite{Bolte2014,Razaviyayn2013,Tseng2009} using a linearization of $f$
\begin{equation}
\label{BCD3}
\begin{split}
x_i^{(k)}=&\argmin_{x_i \in \bbE_i} \langle\nabla f_i^{(k)}(x_i^{(k-1)}),x_i-x_i^{(k-1)}\rangle  \\
&\qquad+r_i(x_i)+ \frac{1}{2\beta_i^{(k)}} \|x_i-x_i^{(k-1)}\|^2.
\end{split} 
\end{equation}
\end{enumerate} 
Incorporating inertial force is a popular and efficient method to accelerate the convergence of first-order methods.  The inertial term was first introduced by Polyak's heavy ball method \cite{POLYAK1964}, which adds to the new direction a momentum term equal to the difference of the two previous iterates; this is also known as extrapolation. 
While the gradient evaluations used in Polyak's method are not affected by the momentum, the famous accelerated gradient method of \citet{Nesterov1983,Nesterov1998,Nesterov2004,Nesterov2005} evaluates the gradients at the points which are extrapolated. In the convex setting, these methods are proved to achieve the optimal convergence rate, while the computational cost of each iteration is essentially unchanged. In the non-convex setting, the heavy ball method was first considered by \citet{Zavriev1993} to solve an unconstrained smooth minimization problem. Two inertial proximal gradient methods were proposed by~\citet{Ochs2014} and~\citet{Bot2016} to solve~\eqref{eq:main} with $s=1$. 
The method considered by~\citet{Ochs2014}, referred to as iPiano, makes use of the inertial force but does not use the extrapolated points to evaluate the gradients. 
iPiano was extended for $s>1$ and analysed by \citet{Ochs2019}. \citet{Pock2016} proposed iPALM to solve~\eqref{eq:main} with $s=2$. \citet{Xu2013,Xu2017} proposed inertial versions of proximal BCD, cf.~\eqref{BCD3}. Xu \& Yin's methods need restarting steps to guarantee the decrease of the objective function. As stated by \citet{Nesterov2004}, this relaxation property  for some problem classes is too expensive and may not allow optimal convergence. In another line of works, it is worth mentioning the randomized BCD methods for solving convex problems; see \citet{Fercoq2015,Nesterov2001}. The analysis of this type of algorithms considers the convergence of the function values in expectation. This is out of the scope of this work.  

\subsection{Contribution}
In this paper, we propose inertial versions for the proximal and proximal gradient BCD methods~\eqref{BCD2} and~\eqref{BCD3}, for solving the non-convex non-smooth problem~\eqref{eq:main} with \emph{multiple blocks}.   
For the inertial version of the proximal gradient BCD~\eqref{BCD3}, two extrapolation points can be used to evaluate gradients and add the inertial force so that the corresponding scheme is more flexible and may lead to significantly better numerical performance compared with the inertial methods using a single extrapolation point; this will be confirmed with some numerical experiments (see Section~\ref{sec:computation} and the supplementary material). The idea of using two different extrapolation points was first used for iPALM to solve~\eqref{eq:main} with two blocks; however, the parameters of the implemented version of iPALM in the experiments by  \citet{Pock2016} are chosen outside the theoretical bounds established in the paper.   
Our methods for solving~\eqref{eq:main} with multiple blocks allow picking deterministically or randomly the block of variables to update; it was empirically observed that randomization may lead to better solutions and/or faster convergence \cite{Xu2017}. 
Another key feature of our methods is that they do not require restarting steps. We extend our methods in the framework of Bregman divergence so that they are more general hence admit potentially more applications.  
To prove the convergence of the whole sequence to a critical point of $F$ and derive its convergence rate, we combine a modification of the convergence proof recipe by~\citet{Bolte2014} with the technique of using auxiliary functions \cite{Ochs2014}. By choosing appropriate parameters that guarantee the convergence, we apply the methods to NMF. We also apply it to non-negative canonical polyadic decomposition (NCPD) in the supplementary material.
\vskip -0.1in
\section{The proposed methods: IBP and IBPG} 
\label{sec:algorithms}


Algorithm~\ref{alg:Fro} describes our two proposed methods: (1) the inertial block proximal method (IBP) which is a proximal BCD method with one extrapolation point,  and (2) the inertial block proximal gradient method (IBPG) which is a proximal gradient BCD method with two  extrapolation points. 

Algorithm~\ref{alg:Fro} includes an outer loop which is indexed by $k$ and an inner loop which is indexed by $j$. At the $j$th iteration of an inner loop, only one block is updated. Table~\ref{notations} summarizes the notation used in the paper.  
\begin{table}[t]
\caption{Notation}
\label{notations}
\vskip 0.05in
\begin{center}
\begin{small}
\begin{tabularx}{0.475\textwidth}{lX}
\toprule
Notation & Definition\\
\midrule
$x^{(k,j)}$ &  $x$ at the $j$th iteration within the $k$th outer loop 
\\
$\tilde{x}^{(k)}$ & the main generated sequence (the output) \\ 
$T_k$ & number of iterations within the $k$th outer loop \\
$f_{i}^{(k,j)}(x_{i})$ & a function of the $i$th block while fixing the latest updated values of the other blocks, i.e., \\ 

 $= f(x_{1}^{(k,j-1)},\ldots, $ & $x_{i-1}^{(k,j-1)},x_{i},x_{i+1}^{(k,j-1)},\ldots,x_{s}^{(k,j-1)})$ \vspace{0.1cm} \\
$F_{i}^{(k,j)}(x_{i})$ & $F_{i}^{(k,j)}(x_{i}) = f_{i}^{(k,j)}(x_{i}) + r_i(x_{i})$\\
$\bar{x}_i^{(k,m)}$  & the value of block $i$ after it has been updated $m$ times during the $k$th outer loop \\
 $d_i^k$    & the total number of times the $i$th block is updated during the $k$th outer loop  \\
 $\bar\alpha_i^{(k,m)}$ & the values of $\alpha_i^{(k,j)} $,\\
$\bar \beta_i^{(k,m)}$ &  the values of $\beta_i^{(k,j)}$, \\
$\bar \gamma_i^{(k,m)}$ & and the values of $ \gamma_i^{(k,j)}$ that are used in 
 \eqref{eq:extrapol1}, \eqref{eq:prox1-Fro}, \eqref{eq:extrapol1-1}, \eqref{eq:prox1-Fro-1}, \eqref{eq:prox1} and \eqref{eq:prox1-1} to update block $i$ from $\bar{x}_i^{(k,m-1)} $ to $\bar{x}_i^{(k,m)}$\\
 $\{\bar{x}_i^{(k,m)}\}_{k\geq 1}$ & the sequence that contains the updates of the $i$th block, i.e., $\{\ldots,\bar{x}_i^{(k,1)},\ldots, \bar{x}_i^{(k,d_i^k)},\ldots\}$ \\
\bottomrule
\end{tabularx}
\end{small}
\end{center}
\vskip -0.2in
\end{table} 
\begin{algorithm}[tb]
\caption{IBP and IBPG} 
\label{alg:Fro} 
\begin{algorithmic} 
\STATE {\bf Initialize}: Choose $\tilde{x}^{(0)}= \tilde{x}^{(-1)}$. Choose a method: IBP or IBPG. Parameters are chosen as in Section \ref{sec:subsequential}.
\FOR{$k=1,\ldots$}
\STATE $x^{(k,0)}=\tilde{x}^{(k-1)}$.
\FOR{$j=1,\ldots,T_k$}
\STATE  Choose $i\in \{1,\ldots,s\}$ deterministically or randomly such that Assumption \ref{assum:loop} is satisfied.  
Let $y_i$ be the value of the $i$th block before it was updated to $x_i^{(k,j-1)}$. 

\STATE {For \textbf{IBP}: extrapolate}  
 \begin{equation}
\label{eq:extrapol1}
\hat{x}_{i}=x_{i}^{(k,j-1)}+\alpha_{i}^{(k,j)}\left(x_{i}^{(k,j-1)}-y_{i}\right),
\end{equation}
\quad and compute
\begin{equation}
\label{eq:prox1-Fro}
x_{i}^{(k,j)}=\argmin_{x_i} F_{i}^{(k,j)}(x_i) + \frac{1}{2\beta_{i}^{(k,j)}} \norm{x_i-\hat x_i}^2.
\end{equation} 

\STATE {For \textbf{IBPG}: extrapolate} 
 \begin{equation}
\label{eq:extrapol1-1}
\begin{array}{ll}
\hat{x}_{i}&=x_{i}^{(k,j-1
)}+\alpha_{i}^{(k,j)}\left(x_{i}^{(k,j-1)}-y_{i}\right),\\
\grave{x}_{i}&=x_{i}^{(k,j-1)}+\gamma_{i}^{(k,j)}\left(x_{i}^{(k,j-1)}-y_{i}\right),
\end{array}
\end{equation}
\quad and compute
\begin{equation}
\label{eq:prox1-Fro-1}
\begin{split}
x_i^{(k,j)}=&\argmin_{x_i} \langle\nabla f_i^{(k,j)}(\grave{x}_{i}),x_i-x_i^{(k,j-1)}\rangle  \\
&\qquad+r_i(x_i)+ \frac{1}{2\beta_i^{(k,j)}} \|x_i-\hat{x}_{i}\|^2. 
\end{split} 
\end{equation} 

\STATE  Let $x_{i'}^{(k,j)}=x_{i'}^{(k,j-1)}$ for $i'\ne i$.
\ENDFOR 
\STATE  Update $\tilde{x}^{(k)}=x^{(k,T_k)}$.
\ENDFOR
\end{algorithmic}
\end{algorithm}
The choice of the parameters $\alpha^{k,j}$, $\beta^{k,j}$ and $\gamma^{k,j}$ in Algorithm~\ref{alg:Fro} that guarantee the convergence will be discussed in Section~\ref{sec:subsequential}. We can observe from \eqref{eq:extrapol1-1} and \eqref{eq:prox1-Fro-1} that using two extrapolation points do not bring extra computation cost when compared with using a single extrapolation point (which happens when $\alpha_i^{(k,j)}=\gamma_i^{(k,j)}$).  We make the following standard assumption throughout this paper. 
\vspace{-0.1in} 
\begin{assumption}
\label{assum:loop}
For all $k$, all blocks are updated after the $T_k$ iterations  performed within the $k$th outer loop, and there exists a positive constant $\bar T$ such that 
$s \leq T_k \leq \bar T$. 
\end{assumption}

\vspace{-0.1in}

\paragraph{Illustration with NMF}  Let us illustrate the proposed methods for NMF; see the supplementary material for the application to NCPD. We will use IBPG for NMF with 2 blocks of variables, namely $U$ and $V$, and IBP with $2\times \mathbf{r}$ blocks of variables, namely $U_{:,i}$ and $V_{i,:}$ 
$(1 \leq i  l \leq \mathbf{r}$). We choose the Frobenius norm for the proximal terms in \eqref{eq:prox1-Fro} and \eqref{eq:prox1-Fro-1}.  
We have $
\nabla_{U} f=UVV^{T}-XV^{T} 
$ and $ 
\nabla_{V}f = U^{T}UV-U^{T}X,
$
hence the inertial proximal gradient step \eqref{eq:prox1-Fro-1} of IBPG is a projected gradient step. If we choose $T_k=2$ for all $k$ then each inner loop of IBPG updates $U$ and $V$ once. Our algorithm also allows to choose $T_k>2$, hence updating $U$ or $V$ several times before updating the other one. As explained by \citet{Gillis2012}, repeating the update of $U$ and $V$ accelerates the algorithm compared to the pure cyclic update rule, because the terms $VV^T$ and $XV^T$ (resp.\@ $U^TU$ and $U^TX$) in the gradient of $U$ (resp.\@ $V$) do not need to be recomputed hence the second evaluation of the gradient is much cheaper; 
namely,  $O(\mathbf m \mathbf r^2)$ (resp.\@ $O(\mathbf n \mathbf r^2)$) vs.\@ $O(\mathbf{mnr})$ operations; while  
$\mathbf r \ll \min(\mathbf m,\mathbf n)$ for most applications.  
Regarding IBP, the inertial proximal step \eqref{eq:prox1-Fro}  has a closed form: 
\vspace{-0.1in} 
\begin{align*}
&\argmin_{U_{:i}\geq0} \sum\frac{1}{2}\big\| X-\sum_{q=1}^{i-1}U_{:q}V_{q:}-\sum_{q=i+1}^{r}U_{:q}V_{q:}-U_{:i}V_{i:}\big\|^2\\
&\qquad\qquad\qquad+\frac{1}{2\beta_{i}}\big\Vert U_{:i}-\hat{U}_{:i}\big\Vert ^{2}\\
&=\max\Big(0,\frac{XV_{i:}^{T}-(UV)V_{i:}^{T}+U_{:i}V_{i:}V_{i:}^{T}+1/\beta_{i}\hat{U}_{:i}}{V_{i:}V_{i:}^{T}+1/\beta_{i}}\Big),
\end{align*}
and a similar update for the rows of $V$ can be derived by symmetry since $\|X-UV\|_F^2 = \|X^T-V^T U^T\|_F^2$. 
For the same reason as above, IBP should update the columns of $U$ and the rows of $V$ several times before doing so for the other one. 
\vspace{-0.1in}
\subsection{Extension to Bregman divergence}
The inertial proximal steps \eqref{eq:prox1-Fro} and \eqref{eq:prox1-Fro-1} can be generalized by replacing $\|.\|$ by a Bregman divergence. Let $H_{i}: \mathbb{E}_{i}\to\mathbb{R}$ be a strictly
convex function that is continuously differentiable. The Bregman distance associated with $H_i$ is defined as: 
$$
D_{i}(u,v)=H_{i}(u)-H_{i}(v)-\left\langle \nabla H_{i}(v),u-v\right\rangle ,\forall u, v \in \mathbb{E}_{i}.
$$ 
The squared Euclidean distance $D_i(u,v)=\frac{1}{2}\|u-v\|_{2}^{2}$ corresponds to $H_i(u)=\frac{1}{2}\|u\|_{2}^{2}$. 
\vspace{-0.05in}
\begin{definition} \label{def:prox} 
 For a given $v\in \mathbb{E}_i$, and a positive number $\beta$,  the Bregman proximal map of a function $\phi$ is defined by
\begin{equation}
\label{eq:prox}
{\rm prox}^{H_i}_{\beta, \phi}(v) := \argmin_{u\in\mathbb{E}_i}\big\{ \phi(u)+\frac{1}{\beta}D_i(u,v)\big\} .
\end{equation}
\end{definition}
\begin{definition}  
\label{def:Gprox}  
For given $u_1 \in {\rm int\, dom } \, g, u_2\in \mathbb{E}_i$  and $\beta>0$, the Bregman proximal gradient map of a pair of functions $(\phi,g)$ ($g$ is continuously differentiable) is defined by  
\begin{equation}
\label{eq:proxg}
\begin{split}
&{\rm Gprox}^{H_i}_{\beta, \phi,g}(u_{1},u_{2})\\
& := \argmin_{u\in\mathbb{E}_i}\big\{ \phi(u)+\left\langle \nabla g(u_{1}),u\right\rangle +\frac{1}{\beta}D_i(u,u_{2})\big\} 
\end{split}
\end{equation} 
\end{definition}\vskip -0.1in
For notation succinctness, whenever the generating function is clear from the context, we omit $H_i$ in the notation of the corresponding Bregman proximal maps. 
As $\phi$ can be non-convex,  ${\rm prox}_{\beta, \phi}(v)$  and ${\rm Gprox}_{\beta, \phi,g}(u_{1},u_{2})$ are set-valued maps in general. Various types of assumptions can be made to guarantee  their well-definedness; see \citet{Eckstein1993}, \citet{Teboulle1997,Teboulle2018} for the well-posedness of  \eqref{eq:prox}, and   \citep[Lemma 3.1]{Bolte2018}, \citep[Lemma 2]{Bauschke2017} for the well-posedness of \eqref{eq:proxg}.   Note that the proximal gradient maps in \cite{Bauschke2017,Bolte2018} use the same point for evaluating the gradient and the Bregman distance while ours allow using two different  points $u_1$ and $u_2$. This modification is important for our analysis; however, it does not affect the proofs of the lemmas in those papers. 

Algorithm \ref{alg:framework} describes IBP and IBPG in the framework of Bregman divergence.  
\begin{algorithm}[tb]
\caption{IBP and IBPG with Bregman divergence } 
\label{alg:framework} 
\begin{algorithmic} 
\STATE {\bf Initialize}: Choose $\tilde{x}^{(0)}= \tilde{x}^{(-1)}$. Choose a method: IBP or IBPG. Parameters are chosen as in Section \ref{sec:subsequential}.
\FOR{$k=1,\ldots$}
\STATE $x^{(k,0)}=\tilde{x}^{(k-1)}$.
\FOR{$j=1,\ldots,T_k$}
\STATE  Choose $i \in \{1,\ldots,s\}$ deterministically or randomly such that Assumption \ref{assum:loop} is satisfied.  
\STATE \textbf{Update of IBP:} extrapolate as in \eqref{eq:extrapol1}
and compute 
\begin{equation}
\label{eq:prox1}
x_{i}^{(k,j)}\in{\rm prox}^{H_i}_{\beta_{i}^{(k,j)},F_{i}^{(k,j)}}\left(\hat{x}_{i}\right).
\end{equation}
\STATE \textbf{Update of IBPG:} extrapolate as in \eqref{eq:extrapol1-1}
and compute
\begin{equation}
\label{eq:prox1-1}
x_{i}^{(k,j)}\in {\rm Gprox}^{H_i}_{\beta_{i}^{(k,j)},r_{i},f_{i}}\lrpar{\grave{x}_{i},\hat{x}_{i}}.
\end{equation}
\STATE  Let $x_{i'}^{(k,j)}=x_{i'}^{(k,j-1)}$ for $i'\ne i$.
\ENDFOR 
\STATE  Update $\tilde{x}^{(k)}=x^{(k,T_k)}$.
\ENDFOR
\end{algorithmic}
\end{algorithm}

Throughout this paper, we assume the following. 
\vspace{-0.1in}
\begin{assumption}
\label{assump:wellposed}
\noindent (A1) The function $H_i$, $i=1,\ldots,s$, is $\sigma_i$-strongly convex, continuously differentiable and $\nabla H_i$ is $L_{H_i}$-Lipschitz continuous.

\noindent  (A2) The proximal maps \eqref{eq:prox} and \eqref{eq:proxg} are well-defined. 
\end{assumption} 
Note that (A1) holds if $H_i$ satisfies $L_{H_i}\mathcal I \preceq \nabla^2 H_i \preceq \sigma_i \mathcal I$. 
A quadratic entropy distance  is a typical example of a Bregman divergence that satisfies (A1) \cite{Deem2019}.  More discussion about important properties and how to evaluate \eqref{eq:prox} and \eqref{eq:proxg} are given in the supplementary material. 
\vspace{-0.1in}
\section{Subsequential convergence}
\label{sec:subsequential}

Before providing the subsequential convergence guarantees, let us elaborate on the notation, in particular $\bar{x}_i^{(k,m)}$ and $d_i^k$ which will be used much in the upcoming analysis, see Table~\ref{notations} for a summary of the notation. 
The elements of the sequence ${x}_i^{(k,j)}$ remain unchanged during many iterations since only one block is updated within each inner loop of Algorithm~\ref{alg:framework}, that is, we will have ${x}_i^{(k,j+1)} = {x}_i^{(k,j)}$ for many $j$'s.  
To simplify the analysis, we introduce the 
subsequence $\bar{x}_i^{(k,m)}$ of ${x}_i^{(k,j)}$ that  will only record the value of the $i$th block when it is actually updated. More precisely, there exists a subsequence $\{ i_1, i_2, \dots, i_{d_i^k} \}$ of $\{1,2,\dots, T_k\}$ such that 
$\bar{x}_i^{(k,m)} = x_{i}^{(k,i_m)}$ for all $m=1,2,\dots,d_i^k$. 
 The previous value of block $i$ before it is updated to $\bar{x}_i^{(k,m)}$ is $\bar{x}_i^{(k,m-1)}$. 
We have $\bar{x}_i^{(k,0)}=\bar{x}_i^{(k-1,d_i^{k-1})}=\tilde{x}_i^{(k-1)}$ and $\bar{x}_i^{(k,d_i^k)}=\tilde{x}_i^{(k)}$.
As for ${x}_i^{(k,j)}$, we use the notation $\bar{x}_i^{(k,-1)}=\bar{x}_i^{(k-1,d_i^{k-1}-1)}$. 

\subsection{Choosing parameters}
\label{sec:choosepara}
We first explain how to choose the parameters for  IBP and IBPG within Algorithm \ref{alg:framework} (note that Algorithm \ref{alg:Fro} is a special case of Algorithm \ref{alg:framework})  such that their subsequential convergence is guaranteed. Let us point out that $\bar\alpha_i^{(k,m)}$, $\bar \beta_i^{(k,m)}$, and $\bar \gamma_i^{(k,m)}$ are the  values of $\alpha_i^{(k,j)} $, $\beta_i^{(k,j)}$ and $ \gamma_i^{(k,j)}$ that are used to update block $i$ from $\bar{x}_i^{(k,m-1)} $ to $\bar{x}_i^{(k,m)}$. 
 
\vspace{-0.1in}
\paragraph{Parameters for IBP}  Let $0<\nu<1,  \delta>1$.  For  $m=1,\ldots,d_i^k$ and $i=1,\ldots,s$, denote $ \theta_i^{(k,m)}=\frac{\lrpar{L_{H_i} \bar\alpha_i^{(k,m)}}^2}{2\nu\sigma_i\bar\beta_i^{(k,m)}}$. Let $\theta_i^{(k,d_i^k+1)}=\theta_i^{(k+1,1)}$. We choose $\bar\alpha_i^{(k,m)}$ and $\bar\beta_i^{(k,m)}$ such that, for $m=1,\ldots,d_i^k$, 
\begin{equation}
\label{eq:par1}
\frac{(1-\nu)\sigma_{i}}{2\bar\beta_{i}^{(k,m)}}\geq \delta\theta_i^{(k,m+1)}. 
\end{equation}
\vspace{-0.2in}

\paragraph{Parameters for IBPG}
Considering IBPG, we need to assume that $\nabla f_{i}^{(k,j)} $ is  $L_{i}^{(k,j)}$-Lipschitz continuous, with $L_{i}^{(k,j)}>0$. 
For notational clarity, we correspondingly use $\bar L_i^{(k,m)}$ for $L_{i}^{(k,j)}$ when updating block $i$ from $\bar{x}_i^{(k,m-1)} $ to $\bar{x}_i^{(k,m)}$. To simplify the upcoming analysis, we  choose $\bar \beta_i^{(k,m)}=\frac{\sigma_i}{\kappa\bar L_i^{(k,m)}}$ with $\kappa>1$. Let $0<\nu<1, \delta>1$. Denote 
\[ 
\lambda_i^{(k,m)} = \frac{ 1}{2}\Big(\bar\gamma_i^{(k,m)}+\frac{\kappa L_{H_i}\bar \alpha_i^{(k,m)}}{\sigma_i} \Big)^2\frac{\bar L_i^{(k,m)}}{\nu(\kappa-1)}, 
\] 
for $m=1, \ldots, d_i^k$  and $i=1,\ldots,s$. 
Let $\lambda_i^{(k,d_i^k+1)}=\lambda_i^{(k+1,1)}$. We choose $\bar \alpha_i^{(k,m)}$, $\bar \beta_i^{(k,m)}$ and $\bar \gamma_i^{(k,m)}$ such that, for  $m=1,\ldots,d_i^k$, 
\begin{equation}
\label{par2}
\frac{(1-\nu)(\kappa-1)\bar L_i^{(k,m)}}{2} \geq \delta\lambda_i^{(k,m+1)}. 
\end{equation}
We make the following standard assumption for the boundedness of the parameters; see \citep[Assumption 2]{Xu2013}, \citep[Assumption 2]{Bolte2014}. 
\vspace{-0.1in}
\begin{assumption}
\label{assump:parameters}
For IBP,  there exist  positive numbers  $W_1$, $\overline{\alpha}$ and  $\underline{\beta}$ such that $\theta_i^{(k,m)}\geq W_1$,  $ \bar{\alpha}_i^{(k,m)}\leq\overline{\alpha}$ and $\underline{\beta}\leq\bar\beta_i^{(k,m)}$,  $\forall \,k\in \bbN$, $m=1, \ldots, d_i^k$, $i=1,\ldots,s$. 

For IBPG, there exist  positive numbers  $W_1$,  $\overline L>0$, $\overline{\alpha}$ and  $\overline{\gamma}$ such that $ \lambda_i^{(k,m)}\geq W_1 $, $\bar L_i^{(k,m)}\leq \overline{L}$, $ \bar{\alpha}_i^{(k,m)}\leq\overline{\alpha}$ and $\bar \gamma_i^{(k,m)} \leq \overline{\gamma}$ for all $k\in \bbN$, $m=1, \ldots, d_i^k$ and $i=1,\ldots,s$. 
\end{assumption}
The algorithm iPALM of~\citet{Pock2016} is a special case of IBPG when $D$ is the Euclidean distance, $s=2$ and the two blocks are cyclically updated; however, our chosen parameters are different. 
In particular, the stepsize $\bar \beta_i^{(k,m)}$ of iPALM depends on the inertial parameters~\citep[Formula 4.9]{Pock2016}, while we choose $\bar \beta_i^{(k,m)}$  independently of $\bar \alpha_i^{(k,m)}$ and $\bar \gamma_i^{(k,m)}$. Our parameters allow using dynamic inertial parameters (see Section \ref{sec:NMF-subconvergence}). As also experimentally tested by \cite{Pock2016}, choosing the inertial parameters dynamically leads to a significant improvement of the algorithm performance. The analysis by \citet{Pock2016} does not support this choice of parameters, while ours guarantee subsequential convergence. 

  \vskip -0.1in
\subsection{Subsequential convergence theory}

The following proposition serves as a cornerstone to prove the subsequential convergence. 

\vspace{-0.1in}
\begin{proposition} 
\label{prop:SubConverge1} 
Let $\{\tilde{x}^{(k)}\}$ be a sequence generated by Algorithm~\ref{alg:framework}, and consider $\big\{\tilde{x}^{(k)}_{\rm prev}\big\}$ with  $\big(\tilde{x}^{(k)}_{\rm prev}\big)_i=\bar x_i^{(k,d_i^{k}-1)}$. Suppose Assumption \ref{assum:loop}, \ref{assump:wellposed} and \ref{assump:parameters} are satisfied.  

(i) We have
$\sum_{k=1}^\infty \big\|\tilde{x}^{(k)}-\tilde x_{\rm prev}^{(k)}\big\|^2<\infty$ and $\sum_{k=1}^\infty\sum_{i=1}^{s}\sum_{m=1}^{d_{i}^{k}}\big\Vert \bar{x}_{i}^{(k,m)}-\bar{x}_{i}^{(k,m-1)}\big\Vert ^{2} <\infty.
$

(ii) If there exists a limit point $x^{*}$ of  $\lrbrace{\tilde{x}^{(k)}}$ (that is, there exists a subsequence $\left\{ \tilde x^{(k_{n})}\right\}$  converging to $x^{*}$), then we have 
$
\lim_{n\to\infty} r_i\lrpar{\bar{x}_{i}^{(k_n,m)}} = r_i\lrpar{x^*_i}. 
$ 
\end{proposition}
\vspace{-0.1in}
\begin{remark}[Relax~\eqref{eq:par1} for block-convex $F$]
\label{remark:bothconvex-a1}
For IBP, if $F$ is block-wise convex then we can 
choose  $\bar\alpha_i^{(k,m)}$ and $\bar\beta_i^{(k,m)}$ satisfying
\begin{equation}
\label{eq:par1-1}
\frac{2(1-\nu)\sigma_{i}}{\bar\beta_{i}^{(k,m)}}\geq  \delta\theta_i^{(k,m+1)}, \quad \text{for}\,\, m=1,\ldots,d_i^k,
\end{equation}
and Proposition~\ref{prop:SubConverge1} still holds. 
Compared to~\eqref{eq:par1}, Condition \eqref{eq:par1-1} allows larger values of the extrapolation parameters $\bar\alpha_i^{(k,m)}$ when using the same stepsize  $\bar\beta_i^{(k,m)}$. 
\end{remark}
\vspace{-0.1in}
\begin{remark}[Relax~\eqref{par2}  for convex $r_i$'s]
\label{remark:convex}
If the functions $r_i$'s are convex (note that $f$ is not necessary block-wise convex) then we can use a larger stepsize. Specifically, we can use   $\bar \beta_i^{(k,m)}=\sigma_i/\bar L_i^{(k,m)}$ and $$ \lambda_i^{(k,m)}= \frac{ 1}{2}\Big(\bar\gamma_i^{(k,m)}+\frac{ L_{H_i}\bar \alpha_i^{(k,m)}}{\sigma_i}\Big)^2\frac{\bar L_i^{(k,m)}}{\nu},$$ and choose  $\bar \alpha_i^{(k,m)}$ and $\bar \gamma_i^{(k,m)}$ satisfying 
\begin{equation}
\label{par2-convex}
\frac{(1-\nu)\bar L_i^{(k,m)}}{2} \geq \delta \lambda_i^{(k,m+1)}, \text{for}\, \, m=1,\ldots,d_i^k, 
\end{equation}
and Proposition~\ref{prop:SubConverge1} still holds.
\end{remark}

\begin{remark}[Relax~\eqref{par2}  for block-convex $f$ and convex $r_i$'s] 
\label{remark:bothconvex} 
If the $r_i$'s are convex and $f(x)$ is block-wise convex, then we can use larger extrapolation parameters. Specifically, we choose $H_i(x_i)=\tfrac{1}{2}\norm{x_i}^2$ and let  $\bar \beta_i^{(k,m)}=1/\bar L_i^{(k,m)}$ and   
\[\lambda_i^{(k,m)}= \lrpar{\lrpar{\bar\gamma_i^{(k,m)}}^2+\frac{\lrpar{\bar\gamma_i^{(k,m)}-\bar \alpha_i^{(k,m)}}^2}{\nu}}\frac{\bar L_i^{(k,m)}}{2},
\]
where $0<\nu<1$, and choose  $\bar \alpha_i^{(k,m)}$ and $\bar \gamma_i^{(k,m)}$ satisfying 
\[\frac{1-\nu}{2}\bar L_i^{(k,m)}\geq \delta \lambda_i^{(k,m+1)}, \text{for}\, \, m=1,\ldots,d_i^k. 
\]
For these values, Proposition~\ref{prop:SubConverge1} still holds.  
In Section~\ref{sec:computation} we numerically show that choosing $\bar{\gamma}_i^{(k,m)}\ne \bar \alpha_i^{(k,m)}$ can significantly improve the performance of the algorithm.
\end{remark}
We now state the local convergence result. The definitions of critical points can be found in the supplementary material.
\vspace{-0.2in}
\begin{theorem} 
\label{thm:newlocal1}
Suppose Assumption \ref{assum:loop}, \ref{assump:wellposed} and \ref{assump:parameters} are satisfied.

(i) For IBP, if $F$ is regular then every limit point of the sequence $\lrbrace{\tilde{x}^{(k)}}$ generated by Algorithm~\ref{alg:framework}  is a critical point type I  of $F$. If $f$ is continuously differentiable  then every limit point is a critical point  type II of $F$.

(ii) For IBPG, every limit point of the sequence $\lrbrace{\tilde{x}^{(k)}}$ generated by Algorithm~\ref{alg:framework} is a critical point type II  of $F$.     
\end{theorem}

\subsection{Choice of parameters for NMF}
\label{sec:NMF-subconvergence}

Let us illustrate the choice of parameters for NMF. 
In the remainder of this paper, in the context of NMF, 
we will refer to IBPG as Algoritm~\ref{alg:Fro} with the choice $T_k=2$ (cyclic update of $U$ and $V$), and to IBPG-A with the choice  $T_k>2$ ($U$ and $V$ are updated several times). IBPG-A is expected to be more efficient; see the discussion in Section~\ref{sec:algorithms}. 
For IBPG and IBPG-A, we take $\bar L_1^{(k,m)}=\tilde L_1^{(k)}=\big\|(\tilde V^{(k-1)})^T \tilde V^{(k-1)}\big\|$ and $\bar L_2^{(k,m)}=\tilde L_2^{(k)}= \big\|(\tilde U^{(k)})^T \tilde U^{(k)}\big\|$ for $m\geq 1$. We take $\bar\beta_i^{(k,m)}=1/\tilde L_i^{(k)}$, $ 
\bar\gamma_i^{(k,m)}=\min\big\{\frac{\tau_k-1}{\tau_k},\tilde{\gamma}\sqrt{\frac{\tilde L_i^{(k-1)}}{\tilde L_i^{(k)}}}\big\}
$
and
$
\bar \alpha_i^{(k,m)}=\breve \alpha \bar\gamma_i^{(k,m)}, 
$
where $\tau_0=1$, $\tau_k=\frac{1}{2}(1+\sqrt{1+4 \tau^2_{k-1}})$, $\tilde{\gamma}=0.99$ and $\breve{\alpha}=1.01$. We can verify that there exists $\delta>1$ such that  $ \breve\gamma^2\lrpar{(\breve \alpha-1)^2/\nu+1} <(1-\nu)/\delta$ with $
\nu=0.0099$. Hence, our choice of parameters satisfy the conditions of Remark~\ref{remark:bothconvex}. 

Regarding IBP, we choose $1/\beta_i^{(k,m)}=0.001$ and $\alpha_i^{(k,m)}=\tilde\alpha^{(k)}=\min(\bar \beta,\gamma \tilde\alpha^{(k-1)})$, with $\bar \beta=1$, $\gamma=1.01 $ and $\tilde\alpha^{(1)}=0.6$.  This choice of parameters satisfies the conditions of Remark~\ref{remark:bothconvex-a1}.

\vspace{-0.1in}
\section{Global convergence} \label{sec:globconv} 
A key tool of the upcoming global convergence (i.e., the whole  sequence converges to a critical point) analysis is the Kurdyka-{\L}ojasiewicz (KL) function defined as follows. 
\begin{definition}
\label{def:KL}
A function $\phi(x)$ is said to have the KL property
at $\bar{x}\in{\rm dom}\,\partial\, \phi$ if there exists $\eta\in(0,+\infty]$,
a neighborhood $U$ of $\bar{x}$ and a concave function $\xi:[0,\eta)\to\mathbb{R}_{+}$
that is continuously differentiable on $(0,\eta)$, continuous at
$0$, $\xi(0)=0$, and $\xi'(s)>0$ for all $s\in(0,\eta),$ such that for all
$x\in U\cap[\phi(\bar{x})<\phi(x)<\phi(\bar{x})+\eta],$ we have
\begin{equation}
\label{ieq:KL}
\xi'\left(\phi(x)-\phi(\bar{x})\right) \, \dist\left(0,\partial\phi(x)\right)\geq1.
\end{equation}
$\dist\left(0,\partial\phi(x)\right)=\min\left\{ \|y\|:y\in\partial\phi(x)\right\}$.
If $\phi(x)$ has the KL property at each point of ${\rm dom}\, \partial\phi$ then $\phi$ is a KL function. 
\end{definition}
The class of KL functions is rich enough to cover many non-convex non-smooth functions found in practical applications. Some noticeable examples are real analytic functions, semi-algebraic functions, and locally strongly convex functions~\cite{Bochnak1998,Bolte2014}.

\vspace{-0.1in}
\subsection{Global convergence recipe} 
\label{sec:globalrecipe}

\citet{Attouch2010,Attouch2013} and \citet{Bolte2014} were the first to prove the global convergence of proximal point algorithms for solving non-convex non-smooth problems. 
We note that a direct deployment of the methodology to our proposed algorithms is not possible since the relaxation property does not hold (that is, the objective functions are not monotonically decreasing) and our methods allow for a randomized strategy. In the following theorem, we  modify the proof recipe of \citet{Bolte2014} so that it is applicable to our proposed methods.
\begin{theorem}
\label{thm:globalex}
Let $\Phi: \mathbb{R}^N\to (-\infty,+\infty]$ be a proper and lower semicontinuous function which is bounded from below. Let $\mathcal{A}$ be a generic algorithm which generates a bounded  sequence $\lrbrace{z^{(k)}}$ by $
z^{(0)}\in \bbR^N$, $z^{(k+1)}\in \calA(z^{(k)})$, $k=0,1,\ldots$
Assume that there exist positive constants  $\rho_1, \rho_2$  and $\rho_3$ and a non-negative sequence  $\lrbrace{\zeta_k}_{k\in \bbN}$ such that the following conditions are satisfied
\vspace{-0.1in}
\begin{itemize}
\item[(B1)] 
\textbf{Sufficient decrease property}: 
\[\rho_1 \|z^{(k)}-z^{(k+1)}\|^2 \leq \rho_2 \zeta_k^2 \leq \Phi(z^{(k)}) - \Phi(z^{(k+1)}),  
\] 
$\forall \,k=0,1,\ldots$
\item[(B2)] \textbf{Boundedness of subgradient}: 
\[
\|w^{(k+1)}\|\leq \rho_3 \zeta_k, w^{(k)}\in \partial \Phi (z^{(k)}), \forall k=0,1,\ldots
\]
\item[(B3)] \textbf{KL property}: $\Phi$ is a KL function.
\item[(B4)] \textbf{A continuity condition}: If a subsequence 
$\{z^{(k_n)}\}$   converges to  $\bar{z}$  then $\Phi(z^{(k_n)})\to \Phi(\bar{z})$ as $n\to \infty$.
\end{itemize}
Then we have $
\sum_{k=1}^\infty \zeta_k<\infty
$
and $\{z^{(k)}\}$ converges to a critical point of $\Phi$. 
\end{theorem}
We remark that if we take $\zeta_k=\|z^{(k)}-z^{(k+1)}\|$ then Theorem \ref{thm:globalex} recovers the proof recipe of \citet{Bolte2014}. 
The following theorem establish the convergence rate under {\L}ojasiewicz property.

\begin{theorem}
\label{thrm:rate}
Suppose $\Phi$ is a KL function and $\xi(a)$ of Definition \ref{def:KL} has the form $\xi(a)= C a^{1-\omega}$ for some $C>0$ and $\omega\in [0,1)$. Then we have
\begin{itemize}
\item[(i)] If  $\omega=0$ then $\{z^{(k)}\}$ converges after a finite number of steps.
\item[(ii)]  If $\omega \in (0, 1/2]$ then there exists $\omega_1 >0$ and $\omega_2 \in [0,1)$ such that $\norm{z^{(k)}- \bar z}\leq \omega_1 \omega_2 ^k$. 
\item[(iii)]  If $\omega \in (1/2, 1)$ then there exists $\omega_1 >0$  such that $\norm{z^{(k)}- \bar z}\leq \omega_1 k^{-(1-\omega)/(2\omega-1)}$. 
\end{itemize}
\end{theorem}
\vspace{-0.1in}
\subsection{Global convergence of IBP and IBPG}
We need the use of the following auxiliary function $$ 
\Psi(\acute{y},\breve{y}) := F(\acute y)+\rho D(\acute y,\breve{y}),
$$ 
where $\rho>0$ and $D(\acute y,\breve{y})=\sum_{i=1}^s D_{i}(\acute y_{i},\breve{y}_{i})$. Recall that $(\tilde{x}^{(k)}_{\rm prev})_i=\bar x_i^{(k,d_i^{k}-1)}$. Let us consider the sequence $\lrbrace{ Y^{(k)}}$ with 
$Y^{(k)}=\big(\acute{y}^{(k)},\breve{y}^{(k)}\big)=\big(\tilde{x}^{(k)},\tilde x_{\rm prev}^{(k)}\big).
$
We then have 
\begin{equation}
\label{eq:PhiYkdef}
\Psi(Y^{(k)})=F(\tilde{x}^{(k)})+\rho D (\tilde{x}^{(k)},\tilde x_{\rm prev}^{(k)}),
\end{equation}
\[\|Y^{(k)}-Y^{(k+1)}\|^2=\|\tilde{x}^{(k)}-\tilde{x}^{(k+1)}\|^2 + \|\tilde x_{\rm prev}^{(k)}-\tilde x_{\rm prev}^{(k+1)}\|^2.
\]
We define   
\begin{align*}
\varphi_{k}^2 &:= \sum_{i=1}^{s}\sum_{m=0}^{d^{(k+1)}_{i}}\Vert \bar{x}_{i}^{(k+1,m)}-\bar{x}_{i}^{(k+1,m-1)}\Vert^2\\
&=\sum_{i=1}^{s}\sum_{m=1}^{d^{(k+1)}_{i}}\Vert \bar{x}_{i}^{(k+1,m)}-\bar{x}_{i}^{(k+1,m-1)}\Vert^2 \\
& \qquad\qquad + \|\tilde{x}^{(k)}-\tilde x_{\rm prev}^{(k)}\|^2. 
\end{align*}
We make the following additional assumption.
\vspace{-0.1in}
\begin{assumption}
\label{assump:global2}
The sequences $\lrbrace{\tilde x^{(k)}}_{k\in \bbN}$ generated  by Algorithm \ref{alg:framework} are bounded. 
\end{assumption} 

In Proposition \ref{prop:gradbound}  we will prove that $\Psi(Y^{(k)})$ is non-increasing; thus,  $\Psi(Y^{(k)})$ is upper bounded by $\Psi(Y^{(-1)})$. Moreover, note that $D(\tilde{x}^{(k)},\tilde x_{\rm prev}^{(k)}) \geq 0$. Hence, from \eqref{eq:PhiYkdef} this implies that $F\lrpar{\tilde{x}^{(k)}}$ is also upper bounded by $\Psi\lrpar{Y^{(-1)}}$. Therefore, we can say that Assumption \ref{assump:global2} is satisfied when $F$ has bounded level sets. 
Denote $\sigma=\min\lrbrace{\sigma_1,\ldots,\sigma_s}$ and $L_H=\max\lrbrace{L_{H_1},\ldots,L_{H_s}}$.  

The following proposition  gives an upper bound for the subgradients and a sufficient decrease property for $\lrbrace{\Psi\lrpar{Y^{(k)}}}$.
\vspace{-0.2in}    
\begin{proposition} 
\label{prop:gradbound}
Suppose Assumption \ref{assum:loop}, \ref{assump:wellposed}, \ref{assump:parameters} and \ref{assump:global2} hold. 

(i) Suppose $f$ is continuously differentiable and $\nabla f$ is Lipschitz continuous on bounded subsets of $\bbE$ (this is a standard assumption, see \citep[Lemma 2.6]{Xu2013}, \citep[Assumption 2 iv]{Bolte2014}).  We have 
$\|\hat q^{(k+1)}\| = O(\varphi_k)$ for some $\hat q^{(k)}\in \partial \Psi\lrpar{Y^{(k)}}$.

(ii) Together with the condition in Proposition \ref{prop:gradbound} (i), assume that there exists a constant $W_2$ such that $\forall \,k\in\bbN$, $m=1,\ldots,d_i^k$ and $i=1,\ldots,s$, we have $\theta_i^{(k,m)} \leq W_2$ for IBP, $ \lambda_i^{(k,m)} \leq W_2$ for IBPG and $\delta> (L_H W_2)/(\sigma W_1)$.  
Let $\rho = \frac{\delta W_1}{L_H}+\frac{W_2}{\sigma}$ in~ \eqref{eq:PhiYkdef} and let $\rho_2=\frac{\delta \sigma W_1}{2L_H}-\frac{W_2}{2}$. 
Then 
\[
\Psi(Y^{(k)})-\Psi(Y^{(k+1)})\geq \rho_2 \varphi_k^2.
\]
\end{proposition}
\vspace{-0.1in}
We are now ready to state our global convergence result. 
\begin{theorem}
\label{thm:global}
Assume $F$ is a KL-function and the conditions of Proposition~\ref{prop:gradbound} are satisfied.  Then the whole sequence $\{\tilde{x}^{(k)}\}$ generated by IBP or IBPG converges to a critical point type II of $F$. 
\end{theorem}
We note that $\norm{Y^{(k)}-Y^*}\geq \norm{\tilde x^{(k)}-x^*}$, hence the convergence rate of the sequence $\{\tilde x^{(k)}\} $ is at least the same order as the rate of $\{Y^{(k)}\}$. If $\Psi$ is a KL function with $\xi(a)=C a^{1-\omega}$, then we can apply Theorem  \ref{thrm:rate} to derive the convergence rate of $\{Y^{(k)}\}$. 
\begin{remark} 
\label{rem:existdelta}
Note that we need the additional condition $\delta>\frac{L_{H}W_{2}}{\sigma W_{1}}$ in order to obtain the global convergence in Theorem \ref{thm:global}. Therefore, it makes sense to show that there exists  $\delta$ such that Condition \eqref{par2} for IBPG (or Condition \eqref{eq:par1} for IBP) is also satisfied. See the supplementary material for the proof. 
\end{remark}
The parameters of IBP  for NMF in Section \ref{sec:NMF-subconvergence} satisfy the conditions for global convergence. 
\vspace{-0.1in}
\section{Numerical results for NMF}
\label{sec:computation} 
\begin{figure*}[ht]
\begin{center}
\begin{tabular}{cc}
\includegraphics[width=0.44\textwidth]{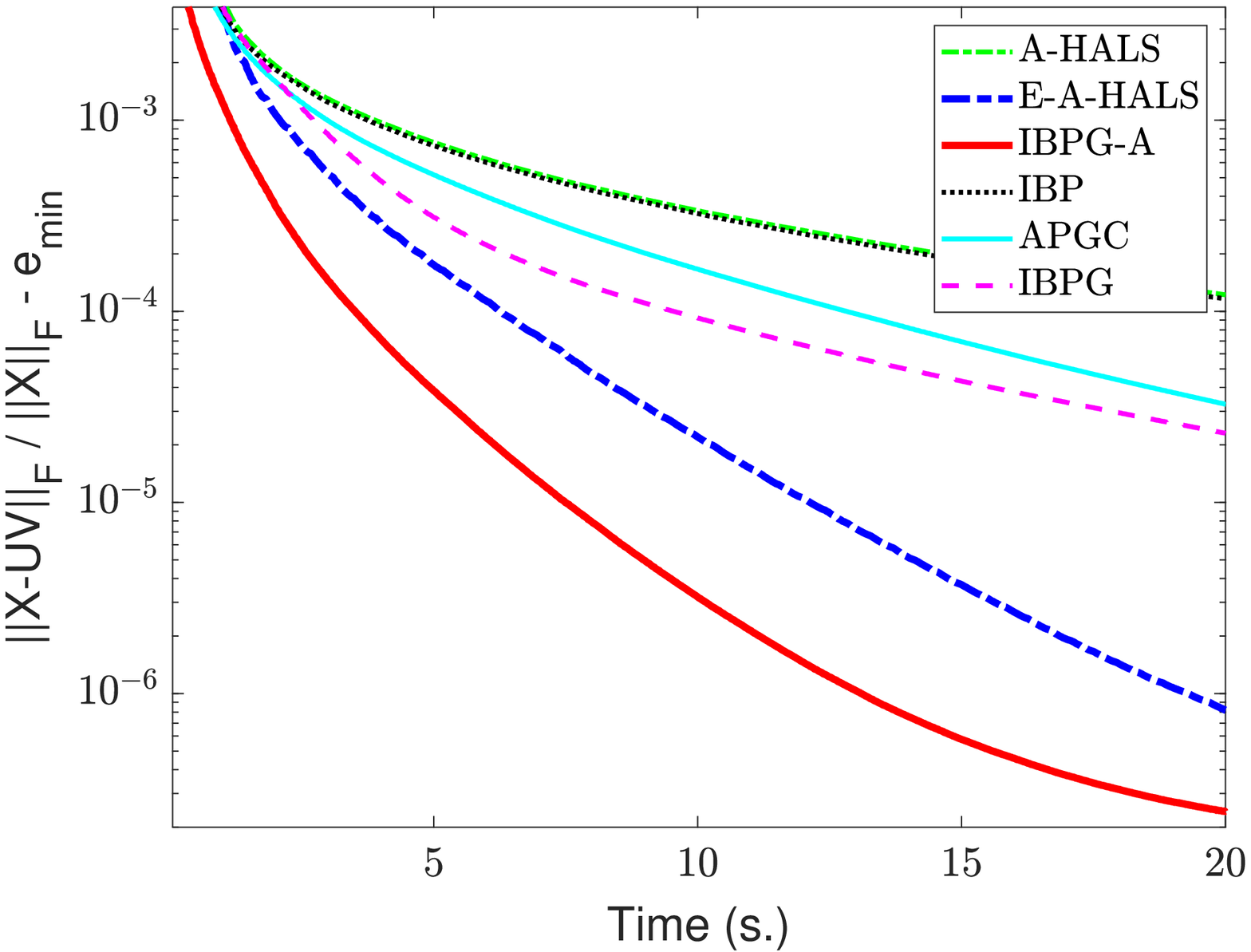}  & 
\includegraphics[width=0.44\textwidth]{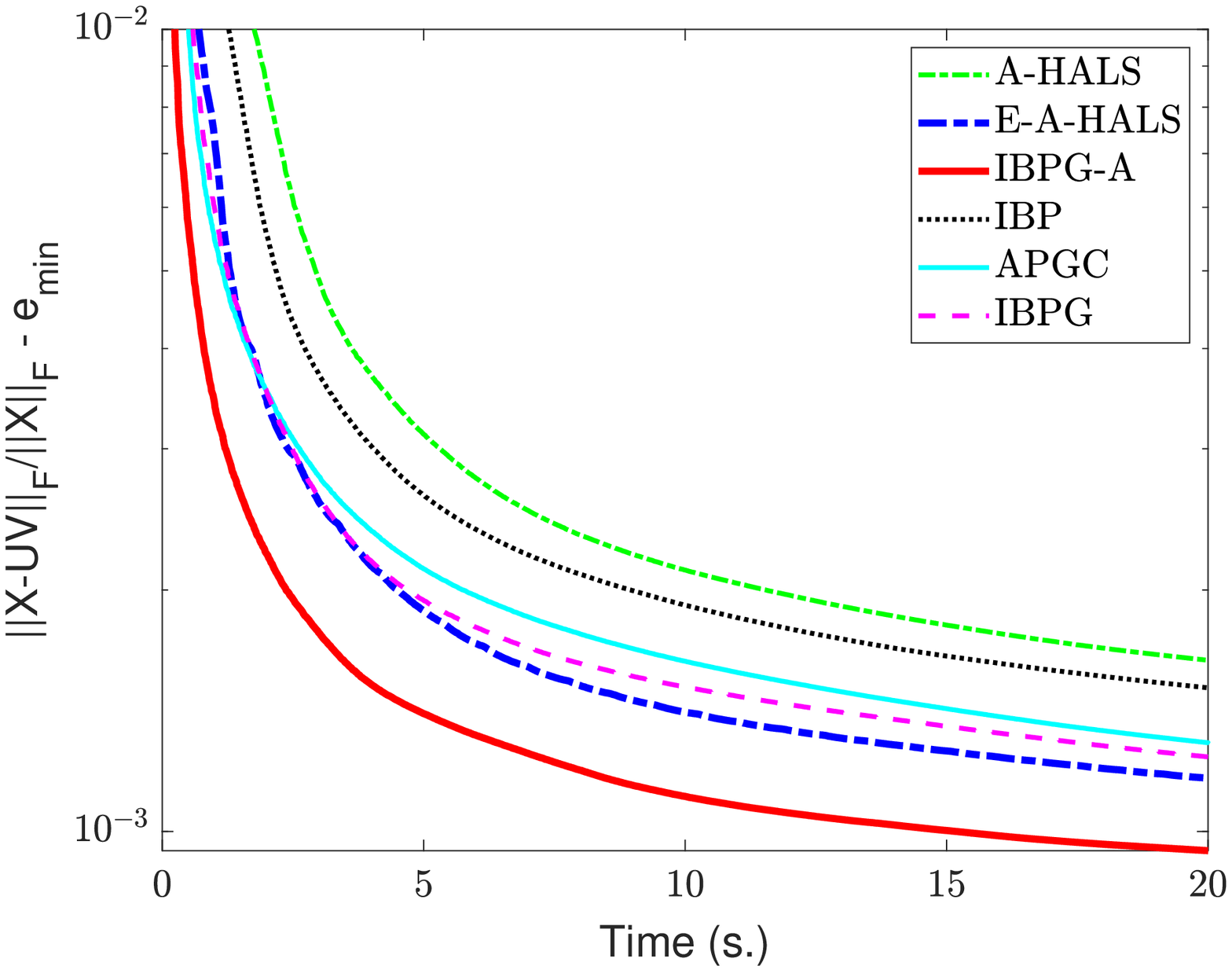} 
\end{tabular}
\caption{Average value of $E(k)$ with respect to time on 2 random low-rank matrices: $200\times 200$ (the left) and $200\times 500$ (the right).
\label{fig:synt}} 
\end{center}
  \vskip -0.2in
\end{figure*} 
\begin{figure*}[ht]
\begin{center}
\begin{tabular}{cc}
\includegraphics[width=0.44\textwidth]{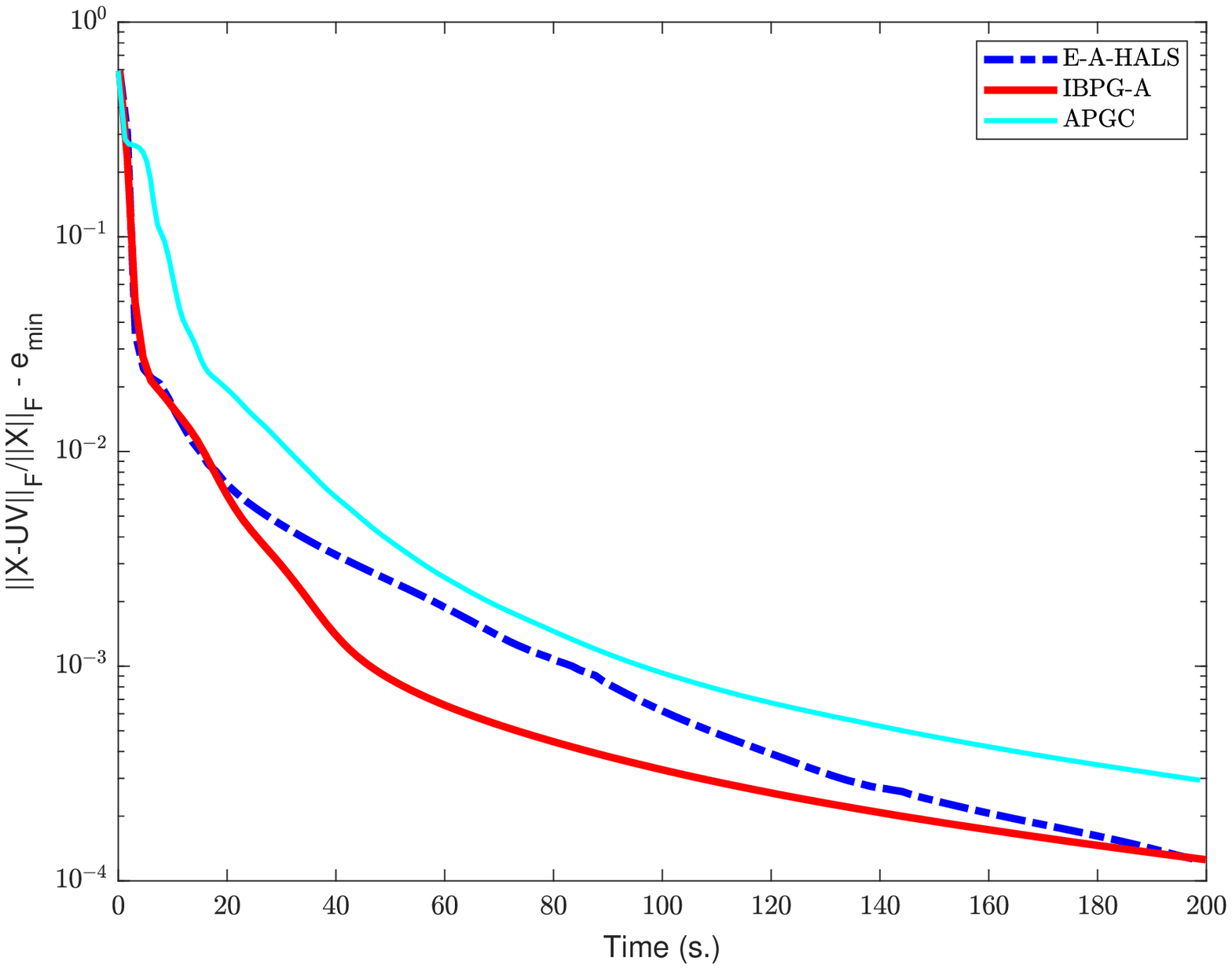}  &
\includegraphics[width=0.44\textwidth]{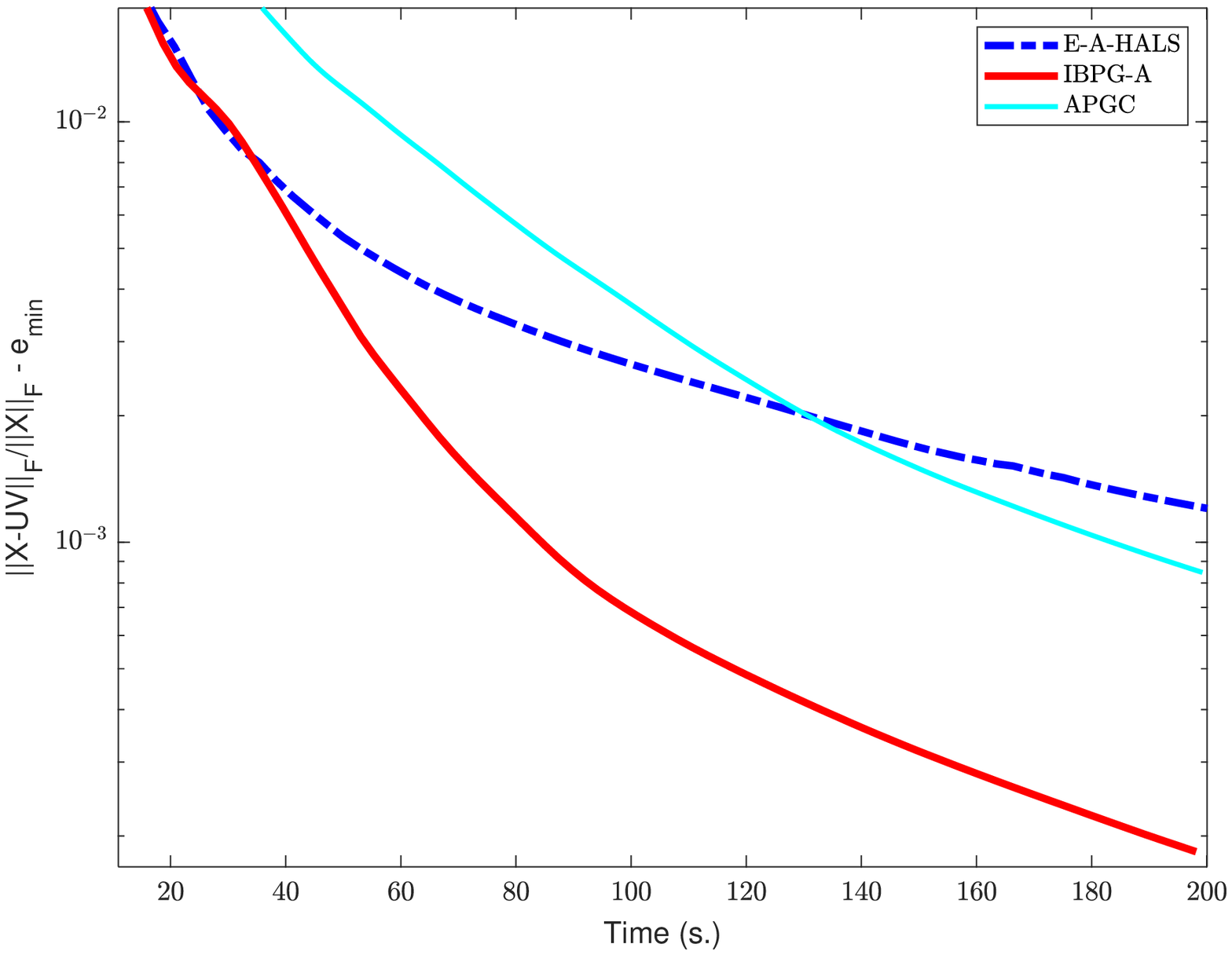} 
\end{tabular}
\caption{Average value of $E(k)$ with respect to time on 2 hyperspectral images: urban (the left) and SanDiego (the right).
\label{fig:urbanSanDiego}} 
\end{center}
  \vskip -0.15in
\end{figure*}
In this section, we compare our IBP, IBPG and IBPG-A (see Section \ref{sec:NMF-subconvergence}) with the following NMF algorithms:  

+ A-HALS: the accelerated hierarchical alternating least squares algorithm \cite{Gillis2012}. A-HALS outperforms standard projected gradient, the popular multiplicative updates and alternating non-negative least squares \cite{kim2014algorithms, Gillis2014}. 

+ E-A-HALS: the acceleration version of A-HALS proposed by~\citet{Ang2018}. This algorithm was experimentally shown to outperform A-HALS. This is, as far as we know, one of the most efficient NMF algorithms. Note that E-A-HALS is a heuristic with no convergence guarantees.

+  APGC: the accelerated proximal gradient coordinate descent method proposed by~\citet{Xu2013} which corresponds exactly to IBPG with $\tilde{\gamma}=\breve{\alpha}=0.9999$. 

We define the relative errors ${\rm relerror}_k=\frac{\norm{X-\tilde U^{(k)} \tilde V^{(k)}}_F}{\norm{X}_F}.$ 
We let $e_{\min}=0$ for the experiments with low-rank synthetic data sets, and in the other experiments $e_{\min}$ is the lowest relative error obtained by any algorithms with any  initializations. We define $E(k)= {\rm relerror}_k - e_{\min}$.
These are the same settings as in~\cite{Gillis2012}. 
All tests are preformed using Matlab
R2015a on a laptop Intel CORE i7-8550U CPU @1.8GHz 16GB RAM. 

\textbf{Experiments with synthetic data sets.
}
Two low-rank matrices of size $200\times 200$ and $200 \times 500$ are generated by letting $X=UV$, where $U$ and $V$ are generated by commands $rand(\mathbf m, \mathbf r)$ and $rand(\mathbf r,\mathbf n)$ with $\mathbf r=20$. 
For each $X$, we run all algorithms with the same 50 random initializations $U_0=rand(\mathbf m,\mathbf r)$, $V_0=rand(\mathbf r,\mathbf n)$, and for each initialization we run each algorithm for 20 seconds. Figure~\ref{fig:synt} illustrates the evolution of the average of $E(k)$ over 50 initializations with respect to time.

To compare the accuracy of the solutions,  we generate 50 random  $\mathbf m\times \mathbf n$ matrices,  $\mathbf m$ and $\mathbf n$ are random integer numbers in the interval [200,500]. For each $X$ we run the algorithms for 20 seconds  with  1 random initialization. Table~\ref{table:syntlowrank} reports the average and standard deviation (std) of the errors. It also provides a ranking between the different algorithms: the $i$th entry of the ranking vector indicates how many times the corresponding algorithm obtained the $i$th best solution.

We observe that (i) in terms of convergence speed and the final errors obtained, IBPG-A outperforms the other algorithms, and (ii) APGC converges slower than IBPG and produces worse solutions. This illustrates the fact that using two extrapolated points may lead to a faster convergence. 

\textbf{Experiments with real data sets.} In these experiments, we will only keep the best performing algorithms, namely IBPG-A and E-A-HALS, along with APGC for our observation purpose. 
For each data set, we generate 35 random initializations and for each initialization we run each algorithm for 200 seconds. 
We test the algorithms on two widely used hyperspectral images, namely the Urban and San Diego data sets; see~\cite{gillis2015hierarchical}.  
We let $\mathbf r=10$. 

Figure~\ref{fig:urbanSanDiego} reports the evolution of the average value of $E(k)$, 
and Table~\ref{table:urbanSanDiego} reports the average error, standard deviation and ranking of the final value of $E(k)$ among the 70 runs (2 data sets with 35 initializations for each data set). 

We see that IBPG-A outperforms E-A-HALS and APGC both in terms of convergence speed and accuracy.
\begin{table}[t] 
 \begin{center} 
 \caption{
Average, standard deviation and ranking of the value of $E(k)$ at the last iteration among the different runs on the low-rank synthetic data sets.  
The best performance is highlighed in bold.  \label{table:syntlowrank} }   \footnotesize
 \begin{tabular}{|p{1.48cm}|p{3.36cm}|p{2.2cm}|} 
 \hline Algorithm &  mean $\pm$ std & ranking  \\ 
 \hline 
A-HALS &  $1.990 \, 10^{-3} \pm 7.910 \, 10^{-4}$ &  (0,  0,  1,  3,  6, 40)   \\ 
 E-A-HALS &  $1.486 \, 10^{-3} \pm 7.233 \, 10^{-4}$ &  (13, 22,  6,  6,  3,  0)   \\ 
IBPG-A &  $\mathbf{1.081 \, 10^{-3}} \pm 6.012 \, 10^{-4}$ 
&  (\textbf{34}, 14,  1,  1,  0,  0)   \\ 
IBP &  $1.916 \, 10^{-3} \pm 7.762 \, 10^{-4}$ &  (0,  1,  5,  9, 34,  1) \\ 
 APGC &  $1.729 \, 10^{-3} \pm 7.452 \, 10^{-4}$ &  (0,  4, 20, 16,  1,  9)   \\ 
IBPG &  $1.672 \, 10^{-3} \pm 7.313 \, 10^{-4}$ &  (3,  9, 17, 15,  6,  0) \\ 
\hline \end{tabular} 
 \end{center} 
 \vskip -0.2in
 \end{table} 
 \begin{table}[t] 
 \begin{center} 
\caption{Average error, standard deviation and ranking among the different runs for urban and SanDiego data sets.
\label{table:urbanSanDiego} }
\vspace{0.15cm}
\small 
 \begin{tabular}{|c|c|c|c|} 
 \hline Algorithm &  mean $\pm$ std & ranking  \\ 
 \hline 
E-A-HALS &  $0.018823  \pm 6.739 \; 10^{-4}$ &  (17, 28, 25)   \\ 
IBPG-A &  $\mathbf{0.018316 } \pm 9.745 \; 10^{-4}$ &  (\textbf{53}, 15,  2)   \\ 
 APGC &  $0.018728  \pm 7.779 \; 10^{-4}$ &  (0, 27, 43)   \\ 
\hline \end{tabular} 
 \end{center} 
  \vskip -0.2in
 \end{table}  
\vspace{-0.1in}
\section{Conclusion}
We have analysed inertial versions of proximal BCD and proximal gradient BCD methods for solving non-convex non-smooth composite optimization problems. Our methods do not require restarting steps, and allow the use of randomized strategies and of two extrapolation points.   
We first proved sub-sequential convergence of the generated sequence to a critical point of $F$ (Theorem~\ref{thm:newlocal1})  and then, under some additional assumptions, convergence of the whole sequence (Theorem~\ref{thm:global}).  
We showed that the proposed methods compared favourably with state-of-the-art algorithms for NMF. Additional experiments on NMF and NCPD are given in the supplementary material. Exploring other Bregman divergences for IBP and IBPG to solve NMF and NCPD  may lead to other efficient algorithms for NMF and NCPD. This is one of our future research directions. 
\newpage 

\bibliography{iBMD}
\bibliographystyle{icml2020}

\onecolumn
\appendix
\begin{center}
\textbf{SUPPLEMENTARY MATERIAL
}\end{center}

\section{Preliminaries}
\label{sec:pre}
In this section, we give important definitions and  properties that allow us to provide our convergence results. 
\subsection{Preliminaries of non-convex non-smooth optimization}
 \label{sec:prelnnopt}

Let $g: \bbE\to \bbR\cup \{+\infty\} $ be a proper lower semicontinuous function.  
\begin{definition}
\label{def:dd}
\begin{itemize}
\item[(i)] For any $x\in{\rm dom}\,g,$  and $d\in\bbE$, we denote the directional derivative of $g$ at $x$ in the direction $d$ by 
\[g'\lrpar{x;d}=\liminf_{\tau \downarrow 0}\frac{g(x+\tau d)-g(x)}{\tau}. 
\] 
\item[(ii)] For each $x\in{\rm dom}\,g,$ we denote $\hat{\partial}g(x)$ as
the Frechet subdifferential of $g$ at $x$ which contains vectors
$v\in\mathbb{E}$ satisfying 
\[
\liminf_{y\ne x,y\to x}\frac{1}{\left\Vert y-x\right\Vert }\left(g(y)-g(x)-\left\langle v,y-x\right\rangle \right)\geq 0.
\]
If $x\not\in{\rm dom}\:g,$ then we set $\hat{\partial}g(x)=\emptyset.$  
\item[(iii)] The limiting-subdifferential $\partial g(x)$ of $g$ at $x\in{\rm dom}\:g$
is defined as follows. 
\[
\partial g(x) := \left\{ v\in\mathbb{E}:\exists x^{(k)}\to x,\,g\left(x^{(k)}\right)\to g(x),\,v^{(k)}\in\hat{\partial}g\left(x^{(k)}\right),\,v^{(k)}\to v\right\} .
\]
\end{itemize}
\end{definition}
The following definition, see \citep[Section 3]{Tseng2001}, is necessary in our convergence analysis for the inertial version of \eqref{BCD2} without the smoothness assumption on $f$. 
\begin{definition}
\begin{itemize}
\item[(i)] We say that 	$x^*\in \rm{dom}\,F$ is a critical point  type I  of $F$ if $ F'(x^*;d) \geq 0, \forall \, d$.  
\item[(ii)]  $x^*\in \rm{dom}\,F$ is said to be a coordinatewise minimum of $F$ if 
\[F\lrpar{x^*+(0,\ldots,d_i,\ldots,0)}\geq F(x^*), \forall\, d_i\in \bbE_i, \forall\, i=1,\ldots,s.
\]
\item[(iii)]	We say that $F$ is regular at $x\in \rm{dom}\,F$ if for all 
$d=\lrpar{d_1,\ldots,d_s}$  such that $$ F'\lrpar{z;(0,\ldots,d_i,\ldots,0)}\geq 0, i=1,\ldots,s,$$ then $F'(x;d)\geq 0$.
\end{itemize}
\end{definition}
It is straightforward to see from the definition that if $F$ is regular at $x^*$ and $x^*$ is a coordinate-wise minimum point of $F$  then $x^*$ is also a critical point  type I of $F$.  
We refer the readers to Lemma 3.1 in \cite{Tseng2001} for the sufficient conditions that imply the regularity of $F$. When $f$ is assumed to be smooth (for the analysis of inertial version of \eqref{BCD3}), Definition \ref{def:type2} will be used. 
\begin{definition}
\label{def:type2}
We call $x^{*}\in \rm{dom}\,F$ a critical point  type II of $F$ if $0\in\partial F\left(x^{*}\right).$ 
\end{definition}
We note that if $x^{*}$ is a minimizer of $F$ then $x^{*}$ is a critical point type I and type II of $F$. 
\subsection{Kurdyka-{\L}ojasiewicz functions} \label{sec:KLprop}
The following lemma (see Lemma 6 of \citealt{Bolte2014}) is the cornerstone to establish the global convergence of our proposed methods.   
\begin{lemma}[Uniformized KL property] 
\label{lem:KL} Let $\phi$ be a proper and lower semicontinuous function.  Assuming that $\phi$ satisfies the KL property and is constant on a compact set $\Omega$.  Then there exist $\varepsilon>0$, $\eta>0$ and a function $\xi$ satisfying the conditions in Definition \ref{def:KL} such that for all $\bar x\in \Omega$ and 
\[x\in \lrbrace{x\in \bbE: {\rm{dist}}(x,\Omega)<\varepsilon} \cap \left[ \phi(\bar x) <\phi(x) < \phi(\bar x)+\eta\right] 
\] 
we have 
\[\xi'\lrpar{\phi(x)-\phi(\bar x)} {\rm{dist}}\, (0,\partial \phi (x))\geq 1.
\]
\end{lemma}
\subsection{Bregman proximal maps} 

We now recall some useful properties of Bregman distance in the following lemmas. Their proofs can be found in \cite{Chen1993,Bauschkebook}. 
\begin{lemma} 
\label{lem:DgeqE}
(i) If $H_i$ is strongly convex with constant $\sigma_i$, that is,  
$$
H_i(u)\geq H_i(v)+\left\langle \nabla H_i(v),u-v\right\rangle +\frac{\sigma_i}{2}\|u-v\|^{2}, \forall  u, v \in \mathbb{E}_{i}
$$
then $$D_i(u,v)\geq\frac{\sigma_i}{2}\|u-v\|^{2}, \forall  u, v \in \mathbb{E}_{i}.$$ 

(ii) If $\nabla H_i$ is $L_{H_i}$-Lipschitz continuous, then 
$$
D_i(u,v)\leq\frac{L_{H_i}}{2}\|u-v\|^{2}, \forall  u, v \in \mathbb{E}_{i}.
$$
\end{lemma}
\begin{lemma}Let $D_i(u,v)$ be
the Bregman distance associated with $H_i$.
\begin{itemize}
\item[(i)]\label{lem:triangle} (The three point identity) 
We have: 
\[
D_i(u,v)+D_i(v,w)=D_i(u,w)+\left\langle \nabla H_i(w)-\nabla H_i(v),u-v\right\rangle, \forall u, v , w \in \mathbb{E}_i .
\]
\item[(ii)]\label{lem:convexcase}(Property 1 of \citealt{Tseng2008}).
Let $ z^+=	\argmin_{u} \phi(u)+D_i(u,z)$, where $\phi$ is a proper convex function. Then for all $u\in \mathbb{E}_i$ we have 
\[
\phi (u) + D_i(u,z) \geq \phi(z^+) + D_i(z^+,z) + D_i(u,z^+).
\]
\end{itemize} 
\end{lemma}

The following inequality is crucial for our convergence analysis.  
\begin{lemma}
\label{lem:nonincrease0} 
For a  given $\hat{w}\in \bbE_i$, if $w^{+}\in {\rm prox}^{H_i}_{\beta, \phi}(\hat w) $
then for all $w \in \mathbb{E}_i$ we have
\[
\phi(w^{+})+\frac{1}{\beta}D_i\left(w^{+},w\right)\leq \phi (w)+\frac{1}{\beta}\left\langle \nabla H_i(\hat{w})-\nabla H_i(w),w^{+}-w\right\rangle .
\]
\end{lemma}  
\proof
It follows from the definition of $w^{+}$ that 
\[
\phi \left(w^{+}\right)+\frac{1}{\beta}D_i\left(w^{+},\hat{w}\right)\leq \phi(w)+\frac{1}{\beta}D_i\left(w,\hat{w}\right).
\]
On the other hand, by Lemma \ref{lem:triangle} (i) we get 
\[
D_i\left(w^{+},\hat{w}\right)-D_i\left(w,\hat{w}\right)=D_i\left(w^{+},w\right)-\left\langle \nabla H_i\left(\hat{w}\right)-\nabla H_i\left(w\right),w^{+}-w\right\rangle .
\]
The result follows. 
\endproof

If $\phi$ is convex,  applying Lemma \ref{lem:convexcase} (ii), we get the following lemma. Lemma \ref{lem:better} will be used when the function $r_i$ is convex. 
\begin{lemma}
\label{lem:better}
For a  given $\hat{w}\in \bbE_i$, if $w^{+}\in {\rm prox}^{H_i}_{\beta, \phi}(\hat w) $ and $\phi$ is convex
then for all $w \in \mathbb{E}_i$ we have
\[
\phi(w^{+})+\frac{1}{\beta}D_i\left(w^{+},\hat w\right) + \frac{1}{\beta} D_i\left(w,w^{+}\right)\leq \phi (w)+\frac{1}{\beta} D_i\lrpar{w,\hat w}.
\]
\end{lemma}

It is crucial to be able to compute efficiently the Bregman proximal maps in~\eqref{eq:prox} and~\eqref{eq:proxg}. 
When $D_i$ is the Euclidean distance, the maps reduce to the classical proximal/proximal gradient maps. We refer the readers to \cite{Boyd} for a comprehensive discussion on how to evaluate the classical maps.  

In~\citep[Section 3.1]{Bauschke2017}, the authors present a splitting mechanism to evaluate~\eqref{eq:proxg} when $u_1$ and $u_2$ are identical.  
Following their methodology, we first define a Bregman gradient operator as follows:  
\[
{\rm p}_{\beta, g}(u_{1},u_{2}) := \argmin\left\{ \left\langle \nabla g(u_{1}),u\right\rangle +\frac{1}{\beta}D_i(u,u_{2}):u\in\mathbb{E}_i\right\} .
\]
Writing the optimality conditions for (\ref{eq:proxg}) together with formal computations (see \citep[Section 3.1]{Bauschke2017} for the details), we can prove that 
\[
{\rm Gprox}_{\beta, \phi,g}(u_{1},u_{2})={\rm prox}_{\beta, \phi}\left({\rm p}_{\beta, g}(u_{1},u_{2})\right),
\]
and 
\begin{equation}
\label{eq:operatorp}
{\rm p}_{\beta, g}(u_{1},u_{2})
=\nabla H_i^{*}\left(\nabla H_i(u_{2})-\beta\nabla  g(u_{1})\right),
\end{equation}
where $H_i^{*}$ is the conjugate function of $H_i$. From \eqref{eq:operatorp}, we see that the calculation of ${\rm p}_{\beta, g}(u_{1},u_{2})$ depends on the calculation of $\nabla H_i^*$. Hence, once we can evaluate $H_i^*$, it is straighforward to evaluate ${\rm p}_{\beta, g}(u_{1},u_{2})$.   A very simple example is the case $D_i\lrpar{u,u_2}=\frac12\norm{u-u_2}^2_2$ for which    ${\rm p}_{\beta, g}(u_{1},u_{2})=u_{2}-\beta\nabla g(u_{1})$; see \cite{Auslender2006,Bauschke2017,Teboulle2018}  for more examples. 
Regarding to the evaluation of~\eqref{eq:prox} in the general setting of Bregman distances, we note that the evaluation can be very difficult and refer the readers to \citep[Section 5]{Bauschke2017}, \citep[Section 5]{Bolte2018} and \citep[Section 6]{Teboulle2018} for some specific examples and discussions. 

\section{Proofs}
\subsection{Proof of Proposition \ref{prop:SubConverge1}}
\label{proof:Subconverge}
\subsection*{Proof for IBP} 
(i)  Applying Lemma \ref{lem:nonincrease0} for (\ref{eq:prox1}) with
$\beta=\beta_{i}^{(k,j)},$ $w=x_{i}^{(k,j-1)}$, $w^{+}=x_{i}^{(k,j)},$
$\hat{w}=\hat{x}_{i}$ we have
\begin{equation*}
\begin{array}{ll}
&F_{i}^{(k,j)}\left(x_{i}^{(k,j)}\right)+\frac{1}{\beta_{i}^{(k,j)}}D_{i}\left(x_{i}^{(k,j)},x_{i}^{(k,j-1)}\right)\\
&\leq F_{i}^{(k,j)}\left(x_{i}^{(k,j-1)}\right)+\frac{1}{\beta_{i}^{(k,j)}}\left\langle \nabla H_{i}\left(\hat{x}_{i}\right)-\nabla H_{i}\left(x_{i}^{(k,j-1)}\right),x_{i}^{(k,j)}-x_{i}^{(k,j-1)}\right\rangle \\
&\lea F_{i}^{(k,j)}\left(x_{i}^{(k,j-1)}\right)+\frac{L_{H_{i}}}{\beta_{i}^{(k,j)}}\norm{\hat{x}_{i}-x_{i}^{(k,j-1)}}\norm{x_{i}^{(k,j)}-x_{i}^{(k,j-1)}}\\
& \leb F_{i}^{(k,j)}\left(x_{i}^{(k,j-1)}\right)+\frac{\lrpar{L_{H_{i}}\alpha_{i}^{(k,j)}}^2}{2\nu\sigma_i\beta_{i}^{(k,j)}}\norm{x_{i}^{(k,j-1)}-y_{i}}^2+\frac{\sigma_i\nu}{2\beta_i^{(k,j)}}\norm{x_{i}^{(k,j)}-x_{i}^{(k,j-1)}}^2,
\end{array}
\end{equation*}
where  we use the Lipschitz continuity of $\nabla H_i$ in (a), use \eqref{eq:extrapol1} and the inequality $ab\leq a^2/(2s) + s b^2/2$ in (b). Together with the inequality $D_{i}\left(w_1,w_2\right)\geq\frac{\sigma_{i}}{2}\left\Vert w_1-w_2\right\Vert ^{2}$ (see Lemma \ref{lem:DgeqE}) and noting that $F\lrpar{x^{(k,j)}}=F_{i}^{(k,j)}\lrpar{x_{i}^{(k,j)}} $, we get
\begin{equation}
\label{eq:recursive1}
\begin{array}{ll}
F\big(x^{(k,j)}\big)+\frac{\sigma_{i}(1-\nu)}{2\beta_{i}^{(k,j)}}\left\Vert x_{i}^{(k,j)}-x_{i}^{(k,j-1)}\right\Vert ^{2}
\leq F\big(x^{(k,j-1)}\big)+\frac{\lrpar{L_{H_{i}}\alpha_{i}^{(k,j)}}^2}{2\nu\sigma_i\beta_{i}^{(k,j)}}\left\Vert x_{i}^{(k,j-1)}-y_{i}\right\Vert ^{2}.
\end{array}
\end{equation}
Note that $y_{i}$, $x_{i}^{(k,j-1)}$ and $x_{i}^{(k,j)}$ are 3 consecutive iterates of $\bar{x}_{i}^{(k,-1)},\ldots,\bar{x}_{i}^{(k,d_{i}^k)}$. Summing up Inequality \eqref{eq:recursive1} for $j=1$ to $T_k$, and combining with  \eqref{eq:par1}   we obtain
\begin{equation}
\label{eq:recursive0}
\begin{array}{ll}
&F\left(x^{(k,T_k)}\right)+\sum_{i=1}^{s}\sum_{m=1}^{d_{i}^{k}} \delta\theta_i^{(k,m+1)}\left\Vert \bar{x}_{i}^{(k,m)}-\bar{x}_{i}^{(k,m-1)}\right\Vert ^{2}
\\
&\leq F\left(x^{(k,0)}\right)+\sum_{i=1}^{s}\sum_{m=1}^{d_{i}^{k}}\theta_i^{(k,m)}\left\Vert \bar{x}_{i}^{(k,m-1)}-\bar{x}_{i}^{(k,m-2)}\right\Vert ^{2}, 
\end{array}
\end{equation}
which implies 
\begin{equation}
\label{recursive1}
\begin{array}{ll}
&F\lrpar{\tilde{x}^{(k)}}+\sum_{i=1}^{s}\delta\theta_i^{(k,d_i^k+1)}\norm{\bar{x}_{i}^{(k,d_i^k)}-\bar{x}_{i}^{(k,d_i^k-1)}}^2\\
&\qquad+\sum_{i=1}^{s}\sum_{m=1}^{d_{i}^{k}-1} (\delta-1)\theta_i^{(k,m+1)}\left\Vert \bar{x}_{i}^{(k,m)}-\bar{x}_{i}^{(k,m-1)}\right\Vert ^{2}
\\
&\leq F\lrpar{\tilde{x}^{(k-1)}}+\sum_{i=1}^s \theta_i^{(k,1)} \norm{\tilde{x}_i^{(k-1)}-\lrpar{\tilde x_{\rm prev}^{(k-1)}}_i}^2, 
\end{array}
\end{equation}
where $\sum_{i=a}^b (.)_i = 0$ if $a>b$. Note that $
\norm{\bar{x}_{i}^{(k,d_i^k)}-\bar{x}_{i}^{(k,d_i^k-1)}}^2=\norm{\tilde x_i^{(k)}-\lrpar{\tilde x_{\rm prev}^{(k)}}_i}^2.
$
Hence from \eqref{recursive1} we get
\begin{equation}
\label{recursive2}
\begin{array}{ll}
&F\lrpar{\tilde{x}^{(k)}}+ \delta\sum_{i=1}^{s} \theta_i^{(k+1,1)} \norm{\tilde x_i^{(k)}-\lrpar{\tilde x_{\rm prev}^{(k)}}_i}^2  \\
&\qquad+\sum_{i=1}^{s}\sum_{m=1}^{d_{i}^{k}-1} (\delta-1)\theta_i^{(k,m+1)}\left\Vert \bar{x}_{i}^{(k,m)}-\bar{x}_{i}^{(k,m-1)}\right\Vert ^{2}
\\\\
&\leq F\lrpar{\tilde{x}^{(k-1)}}+\sum_{i=1}^s \theta_i^{(k,1)} \norm{\tilde{x}_i^{(k-1)}-\lrpar{\tilde x_{\rm prev}^{(k-1)}}_i}^2.
\end{array}
\end{equation}
Summing up Inequality \eqref{recursive2} from $k=1$ to $k=K$ we obtain 
\begin{equation}
\label{recursive3}
\begin{split}
&F\lrpar{\tilde{x}^{(K)}} +(\delta-1)\sum_{k=1}^{K} \sum_{i=1}^s \theta_i^{(k+1,1)} \norm{\tilde{x}_i^{(k)}-\lrpar{\tilde x_{\rm prev}^{(k)}}_i}^2 \\
&\qquad+\sum_{k=1}^{K}\sum_{i=1}^{s}\sum_{m=1}^{d_{i}^{k}-1} (\delta-1)\theta_i^{(k,m+1)}\left\Vert \bar{x}_{i}^{(k,m)}-\bar{x}_{i}^{(k,m-1)}\right\Vert ^{2}
\\ \\
&\leq F\lrpar{\tilde{x}^{(0)}} + \sum_{i=1}^s\theta_i^{(1,1)} \norm{\tilde{x}_i^{(0)}-\tilde{x}_i^{(-1)}}^2.
\end{split}
\end{equation}
Note that $F$ is lower bounded and  $\theta_i^{(k,m)} \geq W_1>0$.  We deduce the result from \eqref{recursive3}.

(ii)
We derive from Proposition \ref{prop:SubConverge1} (i)  that  
\begin{equation}
\label{temp4}
 \lrbrace{\norm{\tilde{x}^{(k)}-\tilde{x}_{\rm {prev}}^{(k)}}}_{k\geq 1} \,\text{and}\, \lrbrace{\sum_{m=1}^{d_i^k}\norm{\bar x_i^{(k,m)}-\bar x_i^{(k,m-1)}}}_{k\geq 1} \, \text{converge to } \, 0.
\end{equation}
By Assumption \ref{assum:loop}, $ d_i^k\geq 1 $ and $d_i^k$ is finite. 
We also note that $\bar{x}_i^{(k_n,d_i^{k_n})}=\tilde{x}_i^{(k_n)}$. Therefore, as $\sum_{m=1}^{d_i^{k_n}}\norm{\bar x_i^{(k_n,m)}-\bar x_i^{(k_n,m-1)}}\to 0$, we  deduce that $\lrbrace{\bar{x}_i^{(k_n,m)}}_{m=0,\ldots,d_i^{k_n}}$ also converges to $x_i^*$. Then,  let $k$ in \eqref{temp4} be $k_n-1$ and note that  $\bar{x}_i^{(k_n,0)}=\tilde{x}_i^{(k_n-1)}$,  $\bar{x}_i^{(k_n,-1)}=\lrpar{\tilde{x}_{\rm {prev}}^{(k_n-1)}}_i$. We thus have  $\bar{x}_i^{(k_n,-1)} \to x_i^*$.
At the $k_n$-th inner loop, we recall that $ \hat x_i=\bar x_i^{(k_n,m-1)}+\bar \alpha_i^{(k_n,m)}\lrpar{\bar{x}_i^{(k_n,m-1)}-\bar{x}_i^{(k_n,m-2)}}.$
Hence  $\hat x_i$
also converges to $x^*_i$. 
From \eqref{eq:prox1}, for all $x_{i}\in\mathbb{E}_{i}$  we have 
\begin{equation}
\label{temp1}
\begin{array}{ll}
&f\left({x}^{(k_n,j)}\right)+r_i\lrpar{\bar{x}_{i}^{(k_n,m)}}+\frac{1}{\bar\beta_{i}^{(k_n,m)}}D_i\left(\bar{x}_{i}^{(k_n,m)},\hat{x}_{i}\right)\\
&\quad
\leq f_{i}^{(k_n,j)}\left(x_{i}\right)+r_i\lrpar{x_i}+\frac{1}{\bar\beta_{i}^{(k_n,m)}}D_i\left(x_{i},\hat{x}_{i}\right).
\end{array} 
\end{equation}
In \eqref{temp1}, let $x_i=x_i^*$ and  let $n\to \infty$ to get 
$
\limsup_{n\to \infty} r_i\lrpar{\bar{x}_{i}^{(k_n,m)}} \leq r_i\lrpar{x^*_i}.
$
Furthermore, as $r_i$ is lower semicontinuous, we have $\liminf_{n\to \infty} r_i\lrpar{\bar{x}_{i}^{(k_n,m)}} \geq r_i\lrpar{x^*_i}$. This completes the proof. 

\subsection*{Proof for IBPG}

(i)  From the assumption that $\nabla f_{i}^{(k,j)}$ is $\bar L_i^{(k,m)}$-Lipschitz continuous, we have 
\begin{equation}
\label{eq:ie3}
\begin{split}
f_{i}^{(k,j)}\left(\bar{x}_i^{(k,m)}\right)&\leq f_{i}^{(k,j)}\left(\bar{x}_i^{(k,m-1)}\right)+\left\langle \nabla f_{i}^{(k,j)}\left(\bar{x}_i^{(k,m-1)}\right),\bar{x}_i^{(k,m)}-\bar{x}_i^{(k,m-1)}\right\rangle \\
&+\frac{\bar L_i^{(k,m)}}{2}\left\Vert \bar{x}_i^{(k,m)}-\bar{x}_i^{(k,m-1)}\right\Vert ^{2}.
\end{split}
\end{equation}
Applying Lemma \ref{lem:nonincrease0} with $\phi(w)=\left\langle \nabla f_{i}^{(k,j)}\left(\grave{x}_{i}\right),w-\bar{x}_i^{(k,m-1)}\right\rangle +r_{i}\left(w\right)$,
$w^{+}=\bar{x}_i^{(k,m)},$ $\hat{w}=\hat{x}_{i},$ and $w=\bar{x}_i^{(k,m-1)}$
we get 
\begin{equation}
\label{eq:ie4}
\begin{split}
\left\langle \nabla f_{i}^{(k,j)}\left(\grave{x}_{i}\right),\bar{x}_i^{(k,m)}-\bar{x}_i^{(k,m-1)}\right\rangle +r_{i}\left(\bar{x}_i^{(k,m)}\right)+\frac{1}{\bar\beta_{i}^{(m)}}D_{i}\left(\bar{x}_i^{(k,m)},\bar{x}_i^{(k,m-1)}\right)\\
\leq r_{i}\left(\bar{x}_i^{(k,m-1)}\right)+\frac{1}{\bar\beta_{i}^{(m)}}\left\langle \nabla H_{i}\left(\hat{x}_{i}\right)-\nabla H_{i}\left(\bar{x}_i^{(k,m-1)}\right),\bar{x}_i^{(k,m)}-\bar{x}_i^{(k,m-1)}\right\rangle .
\end{split}
\end{equation}
 Note that $\frac{\sigma_i}{2}\left\Vert \bar{x}_i^{(k,m)}-\bar{x}_i^{(k,m-1)}\right\Vert ^{2}\leq D_{i}\left(\bar{x}_i^{(k,m)},\bar{x}_i^{(k,m-1)}\right).$
From (\ref{eq:ie3}) and (\ref{eq:ie4}), we get
\begin{align*}
&f_{i}^{(k,j)}\left(\bar{x}_i^{(k,m)}\right)+r_{i}\left(\bar{x}_i^{(k,m)}\right) \\
&\leq f_{i}^{(k,j)}\lrpar{\bar{x}_i^{(k,m-1)}}+\iprod{\nabla f_{i}^{(k,j)} \lrpar{\bar{x}_i^{(k,m-1)}}}{\bar{x}_i^{(k,m)}-\bar{x}_i^{(k,m-1)}}+r_i\lrpar{\bar{x}_i^{(k,m)}}\\
&\qquad+\frac{\bar L_i^{(k,m)}}{\sigma_i}D_i\lrpar{\bar{x}_i^{(k,m)},\bar{x}_i^{(k,m-1)}}\\
&\leq f_{i}^{(k,j)}\lrpar{\bar{x}_i^{(k,m-1)}} + \iprod{\nabla f_{i}^{(k,j)} \lrpar{\bar{x}_i^{(k,m-1)}}-\nabla f_{i}^{(k,j)} \lrpar{\grave{x}_i}}{\bar{x}_i^{(k,m)}-\bar{x}_i^{(k,m-1)}}\\
&\qquad+\lrpar{\frac{\bar L_i^{(k,m)}}{\sigma_i}-\frac{1}{\bar \beta_i^{(k,m)}}}D_i\lrpar{\bar{x}_i^{(k,m)},\bar{x}_i^{(k,m-1)}} + r_i\lrpar{\bar{x}_i^{(k,m-1)}}\\
&\qquad +\frac{1}{\bar\beta_{i}^{(k,m)}}\left\langle \nabla H_{i}\left(\hat{x}_{i}\right)-\nabla H_{i}\left(\bar{x}_i^{(k,m-1)}\right),\bar{x}_i^{(k,m)}-\bar{x}_i^{(k,m-1)}\right\rangle. 
\end{align*}
This implies 
\begin{align*}
&f_{i}^{(k,j)}\left(\bar{x}_i^{(k,m)}\right)+r_{i}\left(\bar{x}_i^{(k,m)}\right) + \frac{(\kappa-1)\bar L_i^{(k,m)}}{2}\norm{\bar{x}_i^{(k,m)}-\bar{x}_i^{(k,m-1)}}^2 \\
&\leq f_{i}^{(k,j)}\lrpar{\bar{x}_i^{(k,m-1)}} + r_i\lrpar{\bar{x}_i^{(k,m-1)}} + \bar L_i^{(k,m)}\norm{\grave x_i-\bar{x}_i^{(k,m-1)}}\norm{\bar{x}_i^{(k,m-1)}-\bar{x}_i^{(k,m)}} \\
&\qquad+ \frac{\kappa L_{H_i}\bar L_i^{(k,m)}}{\sigma_i}\norm{\hat x_i-\bar{x}_i^{(k,m-1)}}\norm{\bar{x}_i^{(k,m-1)}-\bar{x}_i^{(k,m)}}\\
&=f_{i}^{(k,j)}\lrpar{\bar{x}_i^{(k,m-1)}}+  r_i\lrpar{\bar{x}_i^{(k,m-1)}} \\
&\qquad+  \lrpar{\bar\gamma_i^{(k,m)}+\frac{\kappa L_{H_i}\bar \alpha_i^{(k,m)}}{\sigma_i} }\bar L_i^{(k,m)}\norm{\bar{x}_i^{(k,m-1)}-\bar{x}_i^{(k,m-2)}}\norm{\bar{x}_i^{(k,m-1)}-\bar{x}_i^{(k,m)}}
\end{align*}
Note that $x_i^{(k,j)}=\bar x_i^{(k,m)}$, $x_i^{(k,j-1)}=\bar x_i^{(k,m-1)}$ . We apply the Young inequality to get 
\begin{align*}
&F\lrpar{x^{(k,j)}}  + \frac{(\kappa-1)\bar L_i^{(k,m)}}{2}
\norm{\bar{x}_i^{(k,m-1)}-\bar{x}_i^{(k,m)}}^2 \\
&\leq F\lrpar{x^{(k,j-1)}} + \frac{\nu(\kappa-1)\bar L_i^{(k,m)}}{2}\norm{\bar{x}_i^{(k,m-1)}-\bar{x}_i^{(k,m)}}^2\\
&+ \frac{ 1}{2}\lrpar{\bar\gamma_i^{(k,m)}+\frac{\kappa L_{H_i}\bar \alpha_i^{(k,m)}}{\sigma_i} }^2\frac{\bar L_i^{(k,m)}}{\nu(\kappa-1)} \norm{\bar{x}_i^{(k,m-2)}-\bar{x}_i^{(k,m-1)}}^2, 
\end{align*}
where $0<\nu<1$.
We then have  
\begin{equation}
\label{eq:recur2}
\begin{split}
&F\lrpar{x^{(k,j)}}  + \frac{(1-\nu)(\kappa-1)\bar L_i^{(k,m)}}{2} \norm{\bar{x}_i^{(k,m-1)}-\bar{x}_i^{(k,m)}}^2 \\
&\leq F\lrpar{x^{(k,j-1)}}  + \frac{ 1}{2}\lrpar{\bar\gamma_i^{(k,m)}+\frac{\kappa L_{H_i}\bar \alpha_i^{(k,m)}}{\sigma_i} }^2\frac{\bar L_i^{(k,m)}}{\nu(\kappa-1)} \norm{\bar{x}_i^{(k,m-2)}-\bar{x}_i^{(k,m-1)}}^2.
\end{split}
\end{equation}
Summing up Inequality \eqref{eq:recur2} from $j=1$ to $T_k$ we obtain
\begin{equation*}
\begin{split}
&F\lrpar{x^{(k,T_k)}}+ \sum_{i=1}^s\sum_{m=1}^{d_i^{k}} \frac{(1-\nu)(\kappa-1)\bar L_i^{(k,m)}}{2} \norm{\bar x_{i}^{(k,m-1)}-\bar x_{i}^{(k,m)}}^2\\
&\leq F\lrpar{x^{(k,0)}} + \sum_{i=1}^s\sum_{m=1}^{d_i^{k}} \lambda_i^{(k,m)}\norm{\bar x_{i}^{(k,m-2)}-\bar x_{i}^{(k,m-1)}}^2.
\end{split}
\end{equation*}
Together with Condition \eqref{par2}, we see that this inequality is similar to \eqref{eq:recursive0}. Hence, we can use the same technique as in the proof for IBP to obtain the result. 

(ii)  For all $x_i\in \bbE_i$, from \eqref{eq:prox1-1} have 
\begin{equation}
\label{temp7}
\begin{split}
&\iprod{\nabla f_i^{(k,j)}\lrpar{\grave x_i}}{\bar{x}_i^{(k,m)}} + r_i\lrpar{\bar{x}_i^{(k,m)}}+ \frac{1}{\bar{\beta}_i^{k,m}}D_i\lrpar{\bar{x}_i^{(k,m)},\hat{x}_i}\\
 &\leq \iprod{\nabla f_i^{(k,j)}\lrpar{\grave x_i}}{x_i} + r_i\lrpar{x_i} + \frac{1}{\bar{\beta}_i^{k,m}}D_i\lrpar{x_i,\hat{x}_i}.
\end{split}
\end{equation}
Similarly to the proof for IBP, we can prove $\grave{x}_i\to x_i^*$, $\hat x_i\to x_i^*$; and consequently,  by choosing $x_i=x_i^*$ in \eqref{temp7} we have $r_i\lrpar{\bar x_i^{(k_n,m)}} \to r_i(x_i^*)$ as $n\to \infty$.
\subsection{Proof of Remark \ref{remark:bothconvex-a1}}
\label{bothconvex}
Applying Lemma \ref{lem:better} for \eqref{eq:prox1} we have  
\[
\begin{array}{ll}
&F_i^{(i,j)}\lrpar{\bar x_i^{(k,m)}} +\frac{1}{\bar \beta_i^{(k,m)}} D_i\lrpar{\bar x_i^{(k,m)},
\hat x}   +\frac{1}{\bar \beta_i^{(k,m)}} D_i\lrpar{\bar x_i^{(k,m-1)},\bar x_i^{(k,m)}}\\
&\leq F_i^{(i,j)}\lrpar{\bar x_i^{(k,m)}}  +\frac{1}{\bar \beta_i^{(k,m)}} D_i\lrpar{\bar x_i^{(k,m-1)},\hat x}.
\end{array}
\]
Applying Lemma \ref{lem:triangle}, we get 
\begin{equation}
\label{eq:3point}
\begin{array}{ll}
 &D_i\lrpar{\bar x_i^{(k,m)},\hat x_i}-D_i\lrpar{\bar x_i^{(k,m-1)},\hat x_i}\\
 &=D_i\lrpar{\bar x_i^{(k,m)},\bar x_i^{(k,m-1)}}-\iprod{\nabla H_i\lrpar{\hat x_i}-\nabla H_i\lrpar{\bar x_i^{(k,m-1)}}}{\bar x_i^{(k,m)}-\bar x_i^{(k,m-1)}}.
\end{array}
\end{equation}
Therefore, we have 
\[
\begin{array}{ll}
&F_i^{(k,j)}\lrpar{\bar x_i^{(k,m)}}   +\frac{2}{\bar \beta_i^{(k,m)}} D_i\lrpar{\bar x_i^{(k,m-1)},\bar x_i^{(k,m)}}\\
&\leq F_i^{(k,j)}\lrpar{\bar x_i^{(k,m)}} + \frac{1}{\bar \beta_i^{(k,m)}} \iprod{\nabla H_i\lrpar{\hat x_i}-\nabla H_i\lrpar{\bar x^{(k,m-1)}_i}}{{\bar x^{(k,m)}_i}-{\bar x^{(k,m-1)}_i}}  \\
&\leq F_i^{(k,j)}\lrpar{\bar x_i^{(k,m)}} + L_{H_i} \frac{\bar\alpha_i^{(k,m)}}{\bar \beta_i^{(k,m)}}\norm{\bar x_i^{(k,m-1)}-\bar x_i^{(k,m-2)}}\norm{\bar x_i^{(k,m)}-\bar x_i^{(k,m-1)}}\\
&\leq F_i^{(k,j)}\lrpar{\bar x_i^{(k,m)}}+ \frac{\nu\sigma_i}{\bar \beta_i^{(k,m)}}\norm{\bar x_i^{(k,m)}-\bar x_i^{(k,m-1)}}^2 +\frac{\lrpar{L_{H_i}\bar \alpha_i^{(k,m)}}^2}{4\nu\sigma_i \bar \beta_i^{(k,m)}}\norm{\bar x_i^{(k,m-2)}-\bar x_i^{(k,m-1)}}^2.
\end{array}
\]
We then obtain 
\begin{equation}
\begin{array}{ll}
&F\lrpar{x^{(k,j)}}+\frac{(1-\nu)\sigma_i}{\bar\beta_i^{(k,m)}}\norm{\bar x_i^{(k,m)}-\bar x_i^{(k,m-1)}}^2\\
&\leq F\lrpar{x^{(k,j-1)}}+\frac{\lrpar{L_{H_i}\bar \alpha_i^{(k,m)}}^2}{4\nu \sigma_i\bar \beta_i^{(k,m)}}\norm{\bar x_i^{(k,m-2)}-\bar x_i^{(k,m-1)}}^2.
\end{array}
\end{equation}
We have obtained an inequality which is similar to \eqref{eq:recursive1}. 
We therefore continue with the same technique as in the proof of Proposition \ref{prop:SubConverge1} to get the result. 

\subsection{Proof of Remark \ref{remark:convex}}
\label{convexcase}
If $r_i$ is convex then $\iprod{\nabla f_i^{(k,j)}\lrpar{\grave x_i}}{w-\bar x_i^{(k,m-1)}} + r_i(w)$ is also convex. Applying Lemma \ref{lem:better} for \eqref{eq:prox1-1} we have 
\begin{equation}
\label{bound1}
\begin{array}{ll}
&\iprod{\nabla f_i^{(k,j)}\lrpar{\grave x_i}}{\bar x_i^{(k,m)}-\bar x_i^{(k,m-1)}} + r_i(\bar x_i^{(k,m)}) + \frac{\bar L_i^{(k,m)}}{\sigma_i} D_i\lrpar{\bar x_i^{(k,m)},\hat x_i} \\
&\qquad+\frac{\bar L_i^{(k,m)}}{\sigma_i} D_i\lrpar{\bar x_i^{(k,m-1)},\bar x_i^{(k,m)}}\\
& \leq  r_i(\bar x_i^{(k,m-1)}) + \frac{\bar L_i^{(k,m)}}{\sigma_i} D_i\lrpar{\bar x_i^{(k,m-1)},\hat x_i}. 
\end{array}
\end{equation}
Together with \eqref{eq:ie3} we have  
\begin{equation}
\label{eq:temptemp}
\begin{array}{ll}
&f_i^{(k,j)}\lrpar{\bar x_i^{(k,m)}} + r_i\lrpar{\bar x_i^{(k,m)}} + \frac{\bar L_i^{(k,m)}}{\sigma_i} D_i\lrpar{\bar x_i^{(k,m)},\hat x_i} \\
 &\leq f_i^{(k,j)}\lrpar{\bar x_i^{(k,m-1)}} + r_i \lrpar{\bar x_i^{(k,m-1)}}  
+ \frac{\bar L_i^{(k,m)}}{\sigma_i} D_i\lrpar{\bar x_i^{(k,m-1)},\hat x_i}\\
&\qquad+ \iprod{\nabla f_i^{(k,j)}\lrpar{\bar x_i^{(k,m-1)}}-\nabla f_i^{(k,j)}\lrpar{\grave x_i}}{\bar x_i^{(k,m)}-\bar x_i^{(k,m-1)}}.
\end{array}
\end{equation}
 Together with \eqref{eq:3point}  we obtain 
 
 $
\begin{array}{ll}
&F\lrpar{x^{(k,j)}} +\frac{\bar L_i^{(k,m)}}{\sigma_i} D\lrpar{\bar x_i^{(k,m)},\bar x_i^{(k,m-1)}}\\
&\leq F \lrpar{x^{(k,j-1)}}  
+ \frac{\bar L_i^{(k,m)}}{\sigma_i}\iprod{\nabla H_i\lrpar{\hat x_i}-\nabla H_i\lrpar{\bar x_i^{(k,m-1)}}}{\bar x_i^{(k,m)}-\bar x_i^{(k,m-1)}}\\
&\qquad+ \iprod{\nabla f_i^{(k,j)}\lrpar{\bar x_i^{(k,m-1)}}-\nabla f_i^{(k,j)}\lrpar{\grave x_i}}{\bar x_i^{(k,m)}-\bar x_i^{(k,m-1)}},
\end{array}
$

from which we have 

 $
\begin{array}{ll}
&F\lrpar{x^{(k,j)}} +\frac{\bar L_i^{(k,m)}}{2} \norm{\bar x_i^{(k,m)}-\bar x_i^{(k,m-1)}}^2\\
&\leq F \lrpar{x^{(k,j-1)}}  
+ \frac{\bar L_i^{(k,m)}L_{H_i}}{\sigma_i}\bar \alpha_i^{(k,m)}\norm{\bar x_i^{(k,m-1)}-\bar x_i^{(k,m-2)}}\norm{\bar x_i^{(k,m)}-\bar x_i^{(k,m-1)}}\\
&\qquad+ \bar L_i^{(k,m)}\bar \gamma_i^{(k,m)}\norm{\bar x_i^{(k,m-1)}-\bar x_i^{(k,m-2)}}\norm{\bar x_i^{(k,m)}-\bar x_i^{(k,m-1)}}\\
&\leq F \lrpar{x^{(k,j-1)}}  
  +\frac{\nu\bar L_i^{(k,m)}}{2} \norm{\bar x_i^{(k,m)}-\bar x_i^{(k,m-1)}}^2\\
  &\qquad+\frac12\lrpar{\frac{L_{H_i}\bar \alpha_i^{(k,m)}}{\sigma_i}+\bar \gamma_i^{(k,m)}}^2 \frac{\bar L_i^{(k,m)}}{\nu}\norm{\bar x_i^{(k,m-1)}-\bar x_i^{(k,m-2)}}^2,
\end{array}
$

where $0<\nu<1$. Therefore, we have 
\begin{equation}
\label{eq:newrecursive}
\begin{array}{ll}
&F\lrpar{x^{(k,j)}} + \frac{(1-\nu)\bar L_i^{(k,m)}}{2} \norm{\bar x_i^{(k,m)}-\bar x_i^{(k,m-1)}}^2\\
&\leq F \lrpar{x^{(k,j-1)}}  
  +\frac12\lrpar{\frac{L_{H_i}\bar \alpha_i^{(k,m)}}{\sigma_i}+\bar \gamma_i^{(k,m)}}^2 \frac{\bar L_i^{(k,m)}}{\nu}\norm{\bar x_i^{(k,m-1)}-\bar x_i^{(k,m-2)}}^2.
\end{array}
\end{equation}
We get a similar inequality with \eqref{eq:recur2}. Summing up Inequality~\eqref{eq:newrecursive} from $j=1$ to $T_k$ and continuing with the same techniques as in the proof of Proposition~\ref{prop:SubConverge1}, we get the result. 

\subsection{Proof of Remark \ref{remark:bothconvex}}
\label{proofbothconvex}
Using the technique in \citep[Lemma 2.1]{Xu2013}, we first prove that 
\begin{equation}
\label{eq:XuYin}
\begin{array}{ll}
&F\lrpar{x^{(k,j)}}+\frac{(1-\nu)\bar L_i^{(k,m)}}{2}\norm{\bar x_i^{(k,m-1)}-\bar x_i^{(k,m)}}^2\\
&\leq  F\lrpar{x^{(k,j-1)}} + \lrpar{\lrpar{\bar\gamma_i^{(k,m)}}^2+\frac{\lrpar{\bar\gamma_i^{(k,m)}-\bar \alpha_i^{(k,m)}}^2}{\nu}}\frac{\bar L_i^{(k,m)}}{2}\norm{\bar x_i^{(k,m-1)}-\bar x_i^{(k,m-2)}}^2.
\end{array}
\end{equation}
Indeed, we derive from \eqref{eq:prox1-1} that 
\begin{equation}
\label{temp9}
\iprod{\nabla f_i^{(k,j)}\lrpar{\grave x_i}+ \bar g_i^{(k,m)}+\frac{\nabla H_i\lrpar{\bar x_i^{(k,m)}}-\nabla H_i \lrpar{\hat x_i} }{\bar \beta_i^{(k,m)}}}{\bar x_i^{(k,m-1)}-\bar x_i^{(k,m)}}\geq 0,
\end{equation}
where $\bar g_i^{(k,m)}\in \partial r_i \lrpar{\bar x_i^{(k,m)}}$. On the other hand, since $r_i$ is convex and $ f_i^{(k,j)}$ is $\bar L_i^{(k,m)}$-smooth , we have
$
r_i\lrpar{\bar x_i^{(k,m-1)}} - r_i\lrpar{\bar x_i^{(k,m)}}\geq \iprod{\bar g_i^{(k,m)}}{\bar x_i^{(k,m-1)}-\bar x_i^{(k,m)}},
$
and 
$$
f_i^{(k,j)}\lrpar{\grave x_i}-f_i^{(k,j)}\lrpar{\bar x_i^{(k,m)}}+\frac{\bar L_i^{(k,m)}}{2}\norm{\grave x_i-\bar x_i^{(k,m)}}^2\geq \iprod{\nabla f_i^{(k,j)}\lrpar{\grave x_i}}{\grave x_i-\bar x_i^{(k,m)}}.
$$
Together with \eqref{temp9} we have 
\begin{align*}
&r_i\lrpar{\bar x_i^{(k,m-1)}} - r_i\lrpar{\bar x_i^{(k,m)}}-f_i^{(k,j)}\lrpar{\bar x_i^{(k,m)}}+f_i^{(k,j)}\lrpar{\grave x_i}+\frac{\bar L_i^{(k,m)}}{2}\norm{\grave x_i-\bar x_i^{(k,m)}}^2\\
&\geq \frac{1}{\bar \beta_i^{(k,m)}}\iprod{{\nabla H_i \lrpar{\hat x_i} -\nabla H_i\lrpar{\bar x_i^{(k,m)}}}}{\bar x_i^{(k,m-1)}-\bar x_i^{(k,m)}}-\iprod{\nabla f_i^{(k,j)}\lrpar{\grave x_i}}{\bar x_i^{(k,m-1)}-\grave x_i}.
\end{align*}
We then apply the convexity property of $f_i^{(k,j)}$ to obtain 
$$
\begin{array}{lll}
&r_i\lrpar{\bar x_i^{(k,m)}}+f_i^{(k,j)}\lrpar{\bar x_i^{(k,m)}} +\frac{1}{\bar \beta_i^{(k,m)}}\iprod{{\nabla H_i \lrpar{\hat x_i} -\nabla H_i\lrpar{\bar x_i^{(k,m)}}}}{\bar x_i^{(k,m-1)}-\bar x_i^{(k,m)}}\\
&\leq r_i\lrpar{\bar x_i^{(k,m-1)}}+f_i^{(k,j)}\lrpar{\bar x_i^{(k,m-1)}}+\frac{\bar L_i^{(k,m)}}{2}\norm{\grave x_i-\bar x_i^{(k,m)}}^2.
\end{array}
$$
Note that $H_i(x_i)=\norm{x_i}^2/2$. We  then have
$$
\begin{array}{lll}
&F\lrpar{x^{(k,j)}}+\frac{1}{\bar \beta_i^{(k,m)}}\norm{\bar x_i^{(k,m-1)}-\bar x_i^{(k,m)}}^2
\\
&\leq F\lrpar{x^{(k,j-1)}}+\frac{\bar L_i^{(k,m)}}{2}\norm{\bar x_i^{(k,m-1)}-\bar x_i^{(k,m)}}^2+\frac{\bar L_i^{(k,m)}\lrpar{\bar\gamma_i^{(k,m)}}^2}{2}\norm{\bar x_i^{(k,m-1)}-\bar x_i^{(k,m-2)}}^2\\
&\qquad+\bar L_i^{(k,m)} \lrpar{\bar\gamma_i^{(k,m)}-\bar \alpha_i^{(k,m)}}\iprod{\bar x_i^{(k,m-1)}-\bar x_i^{(k,m-2)}}{\bar x_i^{(k,m-1)}-\bar x_i^{(k,m)}},
\end{array}
$$
which implies that
$$
\begin{array}{lll}
&F\lrpar{x^{(k,j)}}+\frac{\bar L_i^{(k,m)}}{2}\norm{\bar x_i^{(k,m-1)}-\bar x_i^{(k,m)}}^2\\
&\leq F\lrpar{x^{(k,j-1)}}+\frac{\bar L_i^{(k,m)}}{2}\lrpar{\bar\gamma_i^{(k,m)}}^2\norm{\bar x_i^{(k,m-1)}-\bar x_i^{(k,m-2)}}^2\\
&+\nu\frac{\bar L_i^{(k,m)}}{2} \norm{\bar x_i^{(k,m-1)}-\bar x_i^{(k,m)}}^2+ \frac{\bar L_i^{(k,m)}\lrpar{\bar\gamma_i^{(k,m)}-\bar \alpha_i^{(k,m)}}^2}{2\nu}\norm{\bar x_i^{(k,m-1)}-\bar x_i^{(k,m-2)}}^2.
\end{array}
$$
Hence we get Inequality \eqref{eq:XuYin}. In other words,  we have a similar inequality with \eqref{eq:recur2}. We then continue with the same techniques as in the proof of Proposition~\ref{prop:SubConverge1} to get the result. 

\subsection{Proof of Theorem \ref{thm:globalex}} 
\label{proof:global}
To prove Theorem \ref{thm:globalex}, we use the same methodology established in \citep{Bolte2014} (see the proof of ~\citep[Theorem 1~(i)]{Bolte2014}). It is worth noting that the same techniques were used in the recent paper~\cite{Ochs2019} to prove an abstract inexact convergence theorem, see Section 3 of \cite{Ochs2019}.

We first prove that $\Phi$ is constant on the set $w\lrpar{z^{(0)}}$ of all limit points of $\lrbrace{z^{(k)}}$. Indeed, from Condition (B1), we derive that $\Phi \lrpar{z^{(k)}}$ is non-increasing. Together with the fact that it is bounded from below,   we deduce that  $\Phi \lrpar{z^{(k)}}$    converges to some value $\bar \Phi$. Therefore, Condition (B4) shows that if $\bar z\in w\lrpar{z^{(0)}}$ then $ \Phi\lrpar{\bar z}=\bar \Phi$. 

 Condition (B1) and the fact that $\Phi$ is bounded from below imply $\norm{z^{(k)}-z^{(k+1)}}\to 0$. As proved in \citep[Lemma 5]{Bolte2014}, we then have $w\lrpar{z^{(0)}}$ is connected and compact. 
  
 If there exists an integer $\bar{k}$ such that $\Phi\lrpar{z^{(\bar{k})}}=\bar\Phi$ is trivial due to Condition (B1).  Otherwise $\Phi\lrpar{\bar{z}} <\Phi \lrpar{z^{(k)}}  $ for all $k>0$. As   $\Phi \lrpar{z^{(k)}} \to \bar \Phi$, we derive that for any $\eta>0$, there exists a positive integer $k_0$ such that $\Phi\lrpar{z^{(k)}}< \Phi\lrpar{\bar{z}}+\eta$ for all $k>k_0$. On the other hand, there exists a positive integer $k_1$ such that $\dist \lrpar{z^{(k)},w\lrpar{z^{(0)}}}<\varepsilon$ for all $k>k_1$.  Applying Lemma \ref{lem:KL} we have 
\begin{equation}
\label{eq:KL}
\xi'\lrpar{\Phi \lrpar{z^{(k)}}-\Phi \lrpar{\bar{z}}} \dist \lrpar{0, \partial \Phi \lrpar{z^{(k)}}}\geq 1, \text{for any}\, k>l:=\max\{k_0,k_1\}.
\end{equation}
From Condition (B2) we get 
\begin{equation}
\label{eq:e1}
\xi'\lrpar{\Phi \lrpar{z^{(k)}}-\Phi \lrpar{\bar{z}}} \geq \frac{1}{\rho_3 \zeta_{k-1}}.
\end{equation}
Denote $ A_{i,j}=\xi\lrpar{\Phi(z^{(i)})-\Phi(\bar z)}-\xi\lrpar{\Phi(z^{(j)})-\Phi(\bar z)}$. From the concavity of $\xi$, Condition (B1) and Inequality \eqref{eq:e1} we obtain
\begin{equation*}
\begin{split}
A_{k,k+1} \geq  \xi'\lrpar{\Phi \lrpar{z^{(k)}}-\Phi \lrpar{\bar{z}}} \left[ \Phi \lrpar{z^{(k)}}-\Phi \lrpar{z^{(k+1)}} \right]\geq \frac{\rho_2\zeta_{k}^2}{\rho_3 \zeta_{k-1}}.
\end{split}
\end{equation*}
Hence we get 
$2\zeta_k\leq 2\sqrt{\frac{\rho_3}{\rho_2} A_{k,k+1}\zeta_{k-1}}\leq \zeta_{k-1}+\frac{\rho_3}{\rho_2} A_{k,k+1}$.
Summing these inequalities from $k=l+1,\ldots,K$ we obtain 
\begin{equation*}
\begin{split}
2\sum_{k=l+1}^K\zeta_k \leq \sum_{k=l}^{K-1}\zeta_k + \frac{\rho_3}{\rho_2}\sum_{k=l+1}^K A_{k,k+1}\leq \sum_{k=l+1}^K\zeta_k +\zeta_l + \frac{\rho_3}{\rho_2} A_{l+1,K+1}.
\end{split}
\end{equation*}
This implies that for all $K>l$ we have
\begin{equation}
\label{ieq:AAk}
\sum_{k=l+1}^K\zeta_k \leq \zeta_l + \frac{\rho_3}{\rho_2}\xi\lrpar{\Phi \lrpar{z^{(l+1)}}-\bar\Phi }.
\end{equation}
Hence, $\sum_{k=1}^\infty \zeta_k <+\infty $. Condition (B1) then gives us $\sum_{k=1}^\infty \norm{z^{(k+1)}-z^{(k)}}<\infty.$
The whole sequence $\lrbrace{z^{(k)}}$ thus converges to some $\bar z$. Together with Condition (B2) and the closedness property of $\partial \Phi$, we have $0\in \partial \Phi(\bar z)$, that is, $\bar z$ is a critical point of $\Phi$. 

\subsection{Proof of Theorem \ref{thrm:rate}}
Inequality \eqref{ieq:KL} becomes 
\begin{equation}
\label{ieq:Loj}
C (1-\omega) {\rm dist} (0, \partial \Phi(z)) \geq (\Phi(z)-\Phi(\bar z))^\omega.
\end{equation} 

If $\omega=0$ then let $I:=\{k\in N: \Phi(z^{(k)})> \Phi(z^*)\}$. Suppose $I$ is infinite, then the sequence $\Big \{ \Phi (z^{(k)})-\Phi(z^*)\Big\}_{k\in I}$ of the right side of \eqref{ieq:Loj} is just a constant 1.  However, the left of \eqref{eq:gradbound1} goes to 0. Hence $I$ is finite; as such, the sequence $\{z^{(k)}\}$ converges in a finite number of steps. 

Denote $\Delta_k=\sum_{p=k}^\infty \zeta_p $, we have 
\[\Delta_k \geq \sum_{p=k}^\infty  \frac{\rho_1}{\rho_2} \norm{z^{(p)}-z^{(p+1)}} \geq \frac{\rho_1}{\rho_2} \norm{z^{(k)}- \bar z}.
\] 

Let us assume $\Phi (z^{(k)})>\Phi(z^*)$ (the case $\Phi (z^{(k_0)})=\Phi(z^*)$ for some $k_0$ is trivial, see proof of Theorem \ref{thm:globalex}), and use the same notations as in the proof of Theorem \ref{thm:globalex}. From Inequality \eqref{ieq:AAk}, which yields $\sum_{l=k}^K\zeta_l \leq \zeta_{k-1} + \frac{\rho_3}{\rho_2}\xi\lrpar{\Phi \lrpar{z^{(k)}}-\bar\Phi }$, and Inequality \eqref{ieq:Loj} we have
\begin{equation*}
\begin{split}
\Delta_k &\leq \Delta_{k-1 } - \Delta_k + \frac{\rho_3}{\rho_2} C \lrpar{\Phi \lrpar{z^{(k)}}-\bar\Phi }^{1-\omega}\\
&\leq  \Delta_{k-1 } - \Delta_k + \frac{\rho_3}{\rho_2} C \lrpar{C(1-\omega) w^{(k)}}^{(1-\omega)/\omega}.
\end{split}
\end{equation*}
Together with Condition (B2) we obtain 
\[\Delta_k \leq \Delta_{k-1 }-\Delta_k + \frac{\rho_3}{\rho_2} C^{1/\omega}((1-\omega)\rho_3)^{(1-\omega)/\omega} \lrpar{\Delta_{k-1}-\Delta_k}^{(1-\omega)/\omega}
\] 
We then can follow the same technique of the proof of \citep[Theorem 2]{Attouch2009} to get the result. 
\subsection{Proof of Proposition \ref{prop:gradbound}}
\label{proof:gradbound}
We first prove the following additional proposition. 
We remind that $\sigma=\min\lrbrace{\sigma_1,\ldots,\sigma_s}$ and $L_H=\max\lrbrace{L_{H_1},\ldots,L_{H_s}}$.  
\begin{proposition} 
\label{prop:basic}
We have
\begin{itemize}
\item[(i)] $\varphi_k^2\geq \frac{1}{2(\bar T-s+1)}\norm{Y^{(k)}-Y^{(k+1)}}^2$. 
\item[(ii)]  Denote 
$$\nabla H=\lrpar{\nabla H_1,\ldots,\nabla H_s}, \textrm {and}\ \nabla^2 H = \lrpar{\nabla^2 H_1,\ldots,\nabla^2 H_s}.$$
Let $q^{(k)}\in \partial
 F\lrpar{\tilde{x}^{(k)}}$. Denote 
 $$\hat{q}^{(k)}=\lrpar{q^{(k)}+ \rho\nabla H\lrpar{\tilde{x}^{(k)}}  -\nabla H \lrpar{\tilde x_{\rm prev}^{(k)}},\rho\nabla^2 H \lrpar{\tilde x_{\rm prev}^{(k)}} \left[\tilde x_{\rm prev}^{(k)}-\tilde{x}^{(k)}\right] }.
 $$
If $f$ is smooth, then we have $\hat q^{(k)} \in \partial \Psi\lrpar{Y^{(k)}},$
 and 
\begin{equation}
\label{eq:uppernorm}
 \norm{\hat{q}^{(k)}}^2\leq 2\norm{q^{(k)}}^2 + O\lrpar{\norm{\tilde{x}^{(k)}-\tilde x_{\rm prev}^{(k)}}^2}.
 \end{equation}
\end{itemize}
\end{proposition}

\proof
(i) We use the inequality $(a_1+\ldots+a_n)^2\leq  n(a_1^2+\ldots+a_n^2)$ and  $d_i^{(k)}\leq \bar T-s+1$ to obtain
\begin{equation}
\label{tempie}
\begin{array}{ll}
\varphi_k^2 &\geq\sum_{i=1}^s \sum_{m=1}^{d_i^{(k+1)}} \norm{\bar{x}_i^{\lrpar{k+1,m}}-\bar{x}_i^{(k+1,m-1)}}^2 \\
& \geq \sum_{i=1}^s \lrpar{1/d_i^{(k+1)}}\norm{\bar{x}_i^{\lrpar{k+1,d_i^{(k+1)}}}-\bar{x}_i^{(k+1,0)}}^2\geq \frac{1}{\bar T-s+1}\norm{\tilde{x}^{(k+1)}-\tilde{x}^{(k)}}^2. 
\end{array}
\end{equation} 
Note that $ \lrpar{\tilde{x}_{\rm {prev}}^{(k+1)}}_i=\bar{x}_i^{\lrpar{k+1,d_i^{(k+1)}-1}}$. Similarly to \eqref{tempie}, we have
\begin{equation}
\label{tempie_2}
\begin{split}
\varphi_k^2=\sum_{i=1}^s \sum_{m=0}^{d_i^{(k+1)}-1} \norm{\bar{x}_i^{\lrpar{k+1,m}}-\bar{x}_i^{(k+1,m-1)}}^2 
\geq \frac{1}{(\bar T-s+1)}\norm{\tilde{x}_{\rm {prev}}^{(k+1)}-\tilde{x}_{\rm {prev}}^{(k)}}^2.
\end{split}
\end{equation}
Summing up \eqref{tempie} and \eqref{tempie_2} we get the result.

(ii) Similarly to   \citep[Proposition 2.1]{Attouch2010}, we can prove that
\begin{equation}
\label{derivative}
\partial \Psi\lrpar{\acute{y},\breve{y}}=\lrbrace{\partial F(\acute{y})+\rho\nabla H(\acute{y})-\nabla H(\breve{y})}\times \lrbrace{\rho\nabla^2 H\lrpar{\breve{y}}[\breve{y}-\acute{y}] }.
\end{equation}
Therefore, $ \hat{q}^{(k)}\in \partial \Psi \lrpar{Y^{(k)}}$. We have 
\begin{align*}
\norm{ \hat{q}^{(k)}}^2
&=\norm{q^{(k)}+\rho\nabla H\lrpar{\tilde x^{(k)}} - \rho\nabla H \lrpar{\tilde x_{\rm prev}^{(k)}}}^2+\rho^2\norm{\nabla^2 H \lrpar{\tilde x_{\rm prev}^{(k)}} \left[\tilde x_{\rm prev}^{(k)}-\tilde x^{(k)}\right]}^2\\
&\leq 2\norm{q^{(k)}}^2+2\rho^2L_H^2 \norm{\tilde x^{(k)} - \tilde x_{\rm prev}^{(k)}}^2+\rho^2\norm{\nabla^2 H \lrpar{\tilde x_{\rm prev}^{(k)}} \left[\tilde x_{\rm prev}^{(k)}-\tilde x^{(k)}\right] }^2,
\end{align*}
where we use the inequality $(a+b)^2\leq 2 a^2+2b^2$ and the Lipschitz continuity of $\nabla H_i$. Finally, note that $\nabla^2 H_i  \preceq L_H \mathbf{I}$, where $\mathbf I$ is the identity operator. Then \eqref{eq:uppernorm} follows.
\endproof

We now prove Proposition \ref{prop:gradbound}.
\subsubsection*{Proof for IBP} 

For all $k\geq 0$, we have $\norm{\tilde{x}^{(k)}}\leq C_1$ (see Assumption \ref{assump:global2}) and $ \norm{\bar x_i^{(k,m)}-\bar x_i^{(k,m-1)}}\leq C_2$ (see Proposition \ref{prop:SubConverge1}). Furthermore, 
$\norm{\bar x_i^{(k,m)}}\leq \norm{\bar x_i^{(k,0)}} + \sum_{j=1}^m\norm{\bar x_i^{(k,j)}-\bar x_i^{(k,j-1)}}
$. Hence, $ \norm{\bar x_i^{(k,m)}} \leq C_1+mC_2\leq C_1+(\bar T-s+1)C_2$. In other words, the sequence $ \lrbrace{\bar x_i^{(k,m)}}_{k\geq 0, m=1,\ldots,d_i^k}$ is bounded. Consequently, the sequence $\lrbrace{x^{(k,j)}}_{k\geq 0, j=1,\ldots,T_k} $  is bounded. 

(i) We denote $\nabla_i f(x) =\nabla_{x_i} f(x)$ and let 
\[\bar q_{i}^{(k,m)}=\frac{1}{\bar\beta_{i}^{(k,m)}}\left(\nabla H_{i}\left(\hat{x}_{i}\right)-\nabla H_{i}\left(\bar x_{i}^{(k,m)}\right)\right)+\nabla_{i} f\lrpar{x^{(k,T_k)}}-\nabla f_{i}^{(k,j)}\lrpar{x_{i}^{(k,j)}}.
\]
Let $L_G$ is the Lipschitz constant of $\nabla f$ on the bounded set containing the sequence $\lrbrace{x^{(k,j)}}_{k\geq 0, j=1,\ldots,T_k} $. From \eqref{eq:prox1} we get $\bar q_{i}^{(k,m)} \in \nabla_{i} f\lrpar{\tilde{x}^{(k)}}+\partial r_{i}\lrpar{\bar x_{i}^{(k,m)}}.$
Also note that $ \nabla f_{i}^{(k,j)}\lrpar{x_{i}^{(k,j)}} = \nabla_i f\lrpar{x^{(k,j)}}$. Hence,
\begin{equation}
\label{temp2}
\begin{array}{lll}
\norm{\bar q_{i}^{(k,m)}}^2&\leq\frac{2L_{H}\left\Vert \hat{x}_{i}-\bar x_{i}^{(k,m)}\right\Vert^2}{\underline\beta}+ 2L_G\left\Vert  x^{(k,j)}-x^{(k,T)}\right\Vert^2 \\
&\leq\frac{4L_H}{\underline\beta} \norm{\bar x_{i}^{(k,m-1)}-\bar x_{i}^{(k,m)}}^2+\frac{4L_H\overline\alpha^2}{\underline\beta}\norm{\bar x_{i}^{(k,m-1)}-\bar x_{i}^{(k,m-2)}}^2\\
& \qquad\qquad +   2L_G\left\Vert x^{(k,j)}-x^{(k,T)}\right\Vert^2.
\end{array}
\end{equation}
We also note that 
$$\begin{array}{lll}
\left\Vert x^{(k,j)}-x^{(k,T)}\right\Vert^2&=O\lrpar{\sum_{i=j}^{T-1} \norm{x^{(k,i)}-x^{(k,i+1)}}^2} \\
&=O\lrpar{\sum_{i=1}^{s}\sum_{m=1}^{d^{(k)}_{i}}\left\Vert \bar{x}_{i}^{(k,m)}-\bar{x}_{i}^{(k,m-1)}\right\Vert^2}.
\end{array}
$$
Therefore, from \eqref{temp2} we deduce that $ \norm{\bar q_{i}^{(k,m)}}^2=O\lrpar{\varphi_{k-1}^2}$.

We now let $\bar q^{(k)}=\lrpar{\bar q_1^{(k,d_1^k)},\ldots,\bar q_s^{(k,d_s^k)}}$. Since $\bar q_i^{(k,d_i^k)}\in \nabla_i f\lrpar{\tilde{x}^{(k)}}+\partial r_i\lrpar{\tilde x_i^{(k)}}$, we have $\bar q^{(k)}\in \partial F \lrpar{\tilde{x}^{(k)}}$ by Proposition \ref{prop:basic}. From $ \norm{\bar q_i^{(k,d_i^k)}}^2=O\lrpar{\varphi_{k-1}^2}$, we can easily obtain 
$
\norm{\bar q^{(k)}}=O\lrpar{\varphi_{k-1}}$. Hence, there exists a positive number $\rho_3$ such that 
\begin{equation}
\label{eq:gradbound1}
\norm{\bar{q}^{(k+1)}}\leq \rho_3 \varphi_k.
\end{equation}
Combined with Proposition \ref{prop:basic}(ii), we get the result. 

(ii) From Inequality \eqref{eq:recursive0}, we have 
\begin{equation}
\label{ie:t1}
\begin{array}{lll}
&F\lrpar{\tilde{x}^{(k-1)}} + \frac{2W_2}{\sigma} \sum_{i=1}^{s}\sum_{m=1}^{d_{i}^{k}}  D_i\lrpar{\bar{x}_{i}^{(k,m-1)},\bar{x}_{i}^{(k,m-2)} } \\
&\geq  F\lrpar{\tilde{x}^{(k)}} + W_2\sum_{i=1}^{s}\sum_{m=1}^{d_{i}^{k}}  \norm{\bar{x}_{i}^{(k,m-1)}-\bar{x}_{i}^{(k,m-2)} }^2 \\
&\geq  F\lrpar{\tilde{x}^{(k)}} + \delta W_1 \sum_{i=1}^{s}\sum_{m=1}^{d_{i}^{k}} \left\Vert \bar{x}_{i}^{(k,m)}-\bar{x}_{i}^{(k,m-1)}\right\Vert ^{2}\\
&\geq F\lrpar{\tilde{x}^{(k)}} + \frac{2\delta W_1}{L_H}\sum_{i=1}^{s}\sum_{m=1}^{d_{i}^{k}}D_i \lrpar{\bar{x}_{i}^{(k,m)},\bar{x}_{i}^{(k,m-1)}} .
\end{array}
\end{equation}
Denote $$ a_k=\sum_{i=1}^{s}\sum_{m=1}^{d_{i}^{k}}  D_i\lrpar{\bar{x}_{i}^{(k,m-1)},\bar{x}_{i}^{(k,m-2)} } \, \text{and} \, 
 b_k=\sum_{i=1}^{s}\sum_{m=1}^{d_{i}^{k}}D_i \lrpar{\bar{x}_{i}^{(k,m)},\bar{x}_{i}^{(k,m-1)}} .$$
From \eqref{ie:t1} we get
$ F\lrpar{\tilde{x}^{(k-1)}} + \frac{2W_2}{\sigma} a_k \geq   F\lrpar{\tilde{x}^{(k)}}  +\frac{2\delta W_1}{L_H} b_k$. 
We thus obtain
\begin{equation}
\label{ie:t2}
\begin{array}{ll}
F\lrpar{\tilde{x}^{(k-1)}} + \rho  a_k - F\lrpar{\tilde{x}^{(k)}} -\rho b_k & \geq \lrpar{\frac{\delta W_1}{L_H}-\frac{W_2}{\sigma}}\lrpar{a_k+b_k}\\
& \geq \lrpar{\frac{\delta W_1}{L_H}-\frac{W_2}{\sigma}}\frac{\sigma}{2}\varphi_k^2=\rho_2\varphi_k^2 .
\end{array}
\end{equation}
Note that $a_k-b_k=D\lrpar{\tilde{x}^{(k-1)},\tilde{x}_{\rm {prev}}^{(k-1)}}-D\lrpar{\tilde{x}^{(k)},\tilde{x}_{\rm {prev}}^{(k)}}$ and 
$a_k+b_k \geq \frac{\sigma}{2}\varphi_{k-1}^2$.  Hence, from \eqref{ie:t2} we deduce that 
$$F\lrpar{\tilde{x}^{(k)}} + \rho  D\lrpar{\tilde{x}^{(k)},\tilde{x}_{\rm {prev}}^{(k)} } - F\lrpar{\tilde{x}^{(k+1)}} -\rho  D\lrpar{\tilde{x}^{(k+1)},\tilde{x}_{\rm {prev}}^{(k+1)} } \geq \rho_2\varphi_k^2.
$$
Together with \eqref{eq:PhiYkdef} we obtain the result.

\subsubsection*{Proof for IBPG}
(i)  Let
$$
\bar{q}_i^{(k,m)}= \frac{1}{\bar\beta_i^{(k,m)}}\lrpar{\nabla H_i\lrpar{\hat x_i}-\nabla H_i \lrpar{\bar x_i^{(k,m)}}}+\nabla_i f\lrpar{x^{(k,T_k)}}-\nabla f_i^{(k,j)}\lrpar{\grave{x}_i}. 
$$
From \eqref{eq:prox1-1}, we get $ \bar{q}_i^{(k,m)} \in \nabla_i f\lrpar{\tilde x^{(k)}} + \partial r_i \lrpar{\bar x_i^{(k,m)}}$.  
Let us recall that the sequences $ \lrbrace{\bar x_i^{(k,m)}}_{k\geq 0, m=1,\ldots,d_i^k}$ and $\lrbrace{x^{(k,j)}}_{k\geq 0, j=1,\ldots,T_k} $  are bounded. Furthermore, we have
\begin{align*}
\norm{\grave{x}_i}&=\norm{\bar x_i^{(k,m-1)}+\bar\gamma_i^{(k,m)}\lrpar{\bar x_i^{(k,m-1)}-\bar x_i^{(k,m-2)}}}\\
&\leq \norm{\bar x_i^{(k,m-1)}}+ \overline{\gamma} \norm{\bar x_i^{(k,m-1)}-\bar x_i^{(k,m-2)}}.
\end{align*}
Hence, $ \grave{x}_i$ is also bounded. As a consequence, the value of $\grave{\mathbf x}$, which is formed by replacing the $i$-th block of $x^{(k,j-1)}$ by $\grave{x}_i = \bar x_i^{(k,m-1)}+\bar\gamma_i^{(k,m)}\lrpar{\bar x_i^{(k,m-1)}-\bar x_i^{(k,m-2)}}$, is also bounded. Let $L_G$ be the Lipschitz constant of $\nabla f$ on the bounded set containing $x^{(k,j)}$ and $\grave{x}$.  Note that $\nabla_i f \lrpar{\grave{\mathbf x}}= \nabla f_i^{(k,j)} \lrpar{\grave x_i}$. We have 
\begin{align*}
&\norm{\nabla_i f\lrpar{x^{(k,T_k)}}-\nabla f_i^{(k,j)}\lrpar{\grave{x}_i}}^2\\
&=\norm{\nabla_i f\lrpar{x^{(k,T_k)}}-\nabla_i f \lrpar{x^{(k,j)}} + \nabla_i f \lrpar{x^{(k,j)}} -\nabla_i f\lrpar{\grave{\mathbf x}}}^2\\
&\leq 2L_G \norm{x^{(k,T_k)}-x^{(k,j)}}^2 + 2L_G \norm{x^{(k,j)}-\grave{\mathbf x}}^2\\
&=2L_G \norm{x^{(k,T_k)}-x^{(k,j)}}^2 + 2L_G \norm{\bar x_i^{(k,m)}-\grave{x}_i}^2.
\end{align*}
We then continue with the same technique as in the proof for IBP to get the bound in \eqref{eq:gradbound1}.

(ii) The proof is follows exactly the same steps as for IBP. 

\subsection{Proof of Theorem \ref{thm:global}}
We now use Theorem~\ref{thm:globalex} to prove the global convergence for both IBP and IBPG. 
We verify the Conditions (B1)-(B4) in Theorem~\ref{thm:globalex} for the auxiliary function $\Psi$ and the sequence $\lrbrace{Y^{(k)}}_{k\in \bbN}$. Proposition~\ref{prop:basic}(i) and  Proposition~\ref{prop:gradbound} show that the Conditions (B1) and (B2) are satisfied. Since $F$ is a KL-function, $\Psi$ is also a KL function. Hence Condition (B3) is satisfied. 

Suppose $Y^*\in w\lrpar{Y^{(0)}}$ is a limit point of  $\lrbrace{Y^{(k)}}$, then there exists a subsequence $\lrbrace{k_n}$ such that $\lrbrace{Y^{(k_n)}}=\lrbrace{\lrpar{\tilde{x}^{(k_n)},\tilde{x}_{\rm {prev}}^{(k_n)}}}$ converges to $Y^*$.  We remind that if $\lrbrace{\tilde{x}^{(k_n)}}$ converges to $x^*$ then  $\lrbrace{\tilde{x}_{\rm {prev}}^{(k_n)}}$ also converges to $x^*$. Hence $Y^*=(x^*,x^*)$. Moreover, from Theorem \ref{thm:newlocal1}, we have $x^*$ is a critical point of $F$, that is, $0\in \partial F(x^*)$. Hence, we derive from \eqref{derivative} that  $0\in \partial \Psi(Y^*)$, that is, $Y^*$ is a critical point of $\Psi$. 
On the other hand,  from Proposition \ref{prop:SubConverge1}(ii)   we have  $F\lrpar{\tilde x^{(k_n)}}\to F(x^*)$ (choose $m=d_i^{k_n}$). Therefore, \eqref{eq:PhiYkdef} implies that  $\Psi\lrpar{Y^{(k_n)}}$ converges to $\Psi(Y^*)$. And consequently, Condition (B4) is satisfied. Applying Theorem \ref{thm:globalex}, we have that the sequence $\lrbrace{Y^{(k)}}$ converges to $(x^*,x^*)$.  Hence the sequence $\lrbrace{\tilde x^{(k)}}$ converges to $x^*$.

\subsection{Proof of Remark \ref{rem:existdelta}}
Let us prove it for IBPG, it would be similar for IBP. Indeed, such $\delta$ would exist if we have  $\frac{(1-\nu)(\kappa-1)\bar{L}_{i}^{(k,m)}}{2\lambda_{i}^{(k,m+1)}}  >\frac{L_{H}W_{2}}{\sigma W_{1}}$, which would be satisfied if 
$
\frac{\nu(1-\nu)(\kappa-1)^{2}\bar{L}_{i}^{(k,m)}}{\bar{L}_{i}^{(k,m+1)}\xi_i^{(k,m+1)}} >\frac{L_{H}W_{2}}{\sigma W_{1}}$,  where $\xi_i^{(k,m)}=\big(\bar\gamma_i^{(k,m)}+\frac{\kappa L_{H_i}\bar \alpha_i^{(k,m)}}{\sigma_i}\big)^2$. In other words, $\delta$ would exist if we have 
\begin{equation}
\label{provedelta}
\frac{\sigma\nu(1-\nu)(\kappa-1)^{2}\bar{L}_{i}^{(k,m)}}{L_{H}\bar{L}_{i}^{(k,m+1)}} \frac{W_{1}}{W_{2}}>\xi_i^{(k,m+1)}. 
\end{equation}
Suppose $\xi_{1}\leq\xi_i^{k,m}\leq\xi_{2}$ and $0<L_{1}\leq\bar{L}^{(k,m)}\leq L_{2}$. We then have $\frac{W_{1}}{W_{2}}=\frac{\xi_{1}L_{1}}{\xi_{2}L_{2}}$, and \eqref{provedelta} holds if $
\frac{\sigma\nu(1-\nu)(\kappa-1)^{2}\bar{L}_{i}^{(k,m)}L_{1}}{L_{H}\bar{L}_{i}^{(k,m+1)}L_{2}}>\xi_{2}$. Therefore, if we choose in advance two constants $\xi_1$ and $\xi_2$ such that $\xi_{2}<\frac{\sigma\nu(1-\nu)(\kappa-1)^{2}L_{1}^{2}}{L_{H}L_{2}^{2}}$
and $0<\xi_{1}<\xi_{2}$, then there always exists $\xi_i^{(k,m)}$ accordingly such that Condition \eqref{provedelta} is satisfied.


\section{Additional experiments}
\subsection{Experiments on NMF}
\subsubsection*{Full-rank synthetic data sets}
\begin{figure}[ht!]
\begin{center}
\begin{tabular}{cc}
\includegraphics[width=0.5\textwidth]{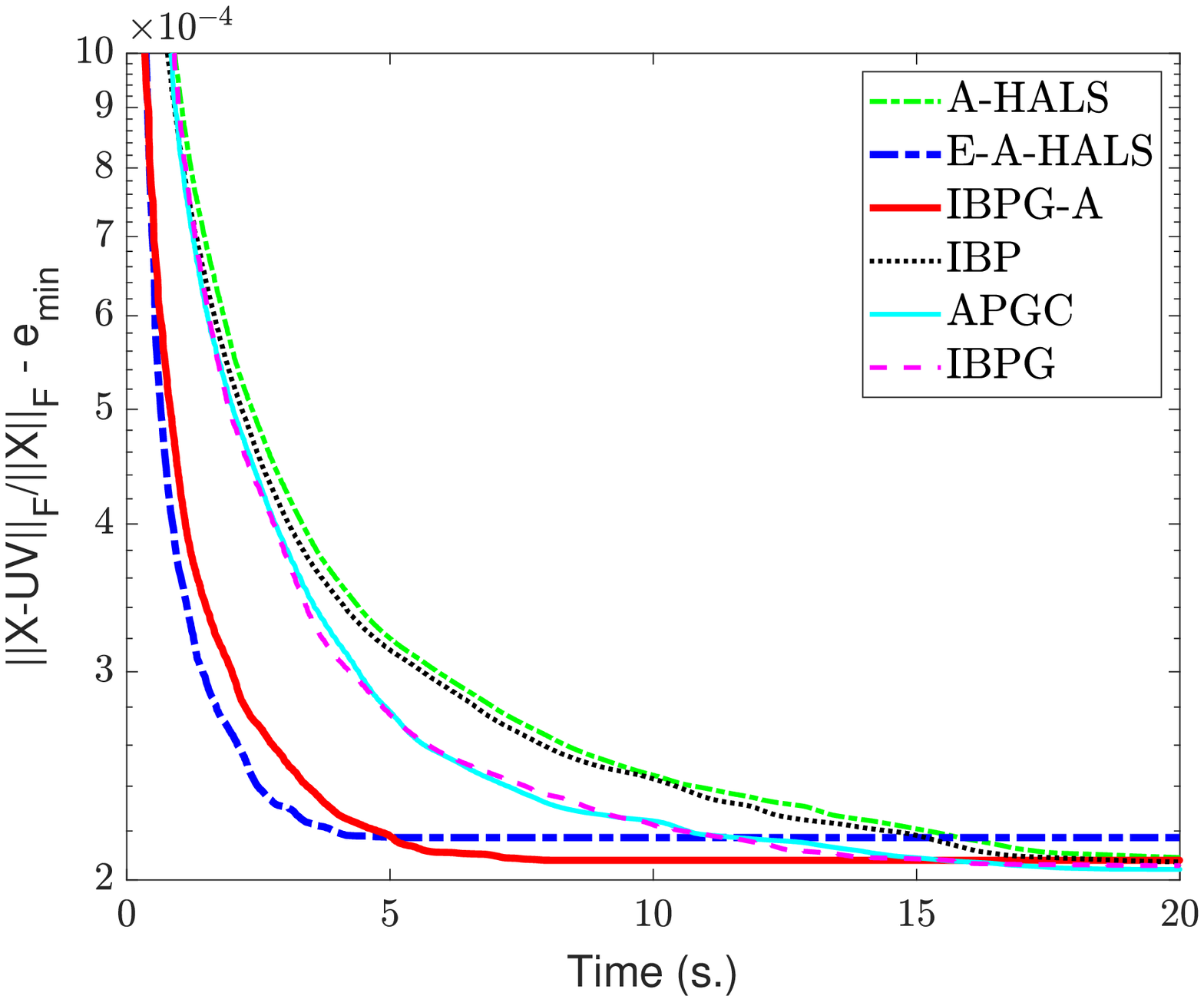}  & 
\includegraphics[width=0.5\textwidth]{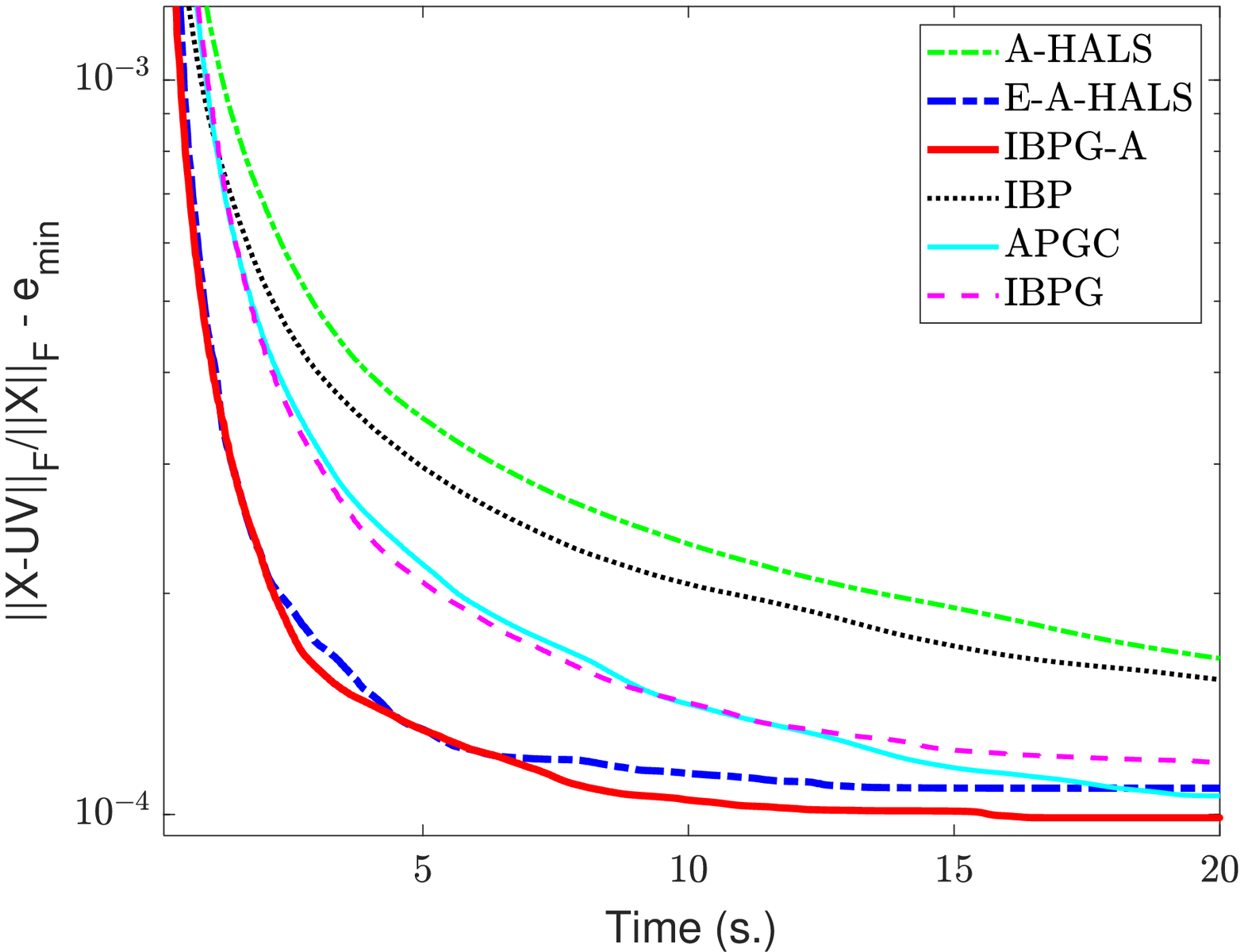} 
\end{tabular}
\caption{ Average value of $E(k)$ with respect to time on 2 random full-rank matrices: $200\times 200$ (left) and $200\times 500$ (right).
\label{fig:synt_fullrank}} 
\end{center}
\end{figure}

Two full-rank matrices of size $200\times 200$ and $200 \times 500$ are generated by MATLAB command $X=rand(m,n)$. We take $\mathbf r=20$. For each matrix $X$, we run all algorithms with the same 50 random initializations $W_0=rand(\mathbf m,\mathbf r)$ and $V_0=rand(\mathbf r,\mathbf n)$, and for each initialization we run each algorithm for 20 seconds. Figure~\ref{fig:synt_fullrank} illustrates the evolution of the average of $E(k)$ over 50 initializations with respect to time. 

We then generate 50 full-rank matrices $X=rand(m,n)$, with $\mathbf m$ and $\mathbf n$ being random integer numbers in the interval [200,500]. For each matrix $X$, we run the algorithms for 20 seconds  with a single random initialization. Table~\ref{table:fullrank_accuracy} reports the average, standard deviation (std) and ranking of the relative errors. 

\begin{center}  
 \begin{table}[h!] 
 \begin{center} 
\caption{Average, standard deviation and ranking of the value of $E(k)$ at the last iteration among the different runs on full-rank synthetic data sets.  The best performance is highlighted in bold. 
\label{table:fullrank_accuracy} } 
 \vspace{0.15cm}
 \begin{tabular}{|c|c|c|c|} 
 \hline Algorithm &  mean $\pm$ std & ranking  \\ 
 \hline 
A-HALS &  $0.450174  \pm 9.048 \, 10^{-3}$ &  (0,  2,  7,  6, 12, 23)   \\ 
E-A-HALS &  $0.450127  \pm 9.028 \, 10^{-3}$ &  (18,  8,  8,  9,  1,  6)   \\ 
IBPG-A  &  $\textbf{0.450115}  \pm 9.050 \, 10^{-3}$ &  ({\textbf{19}}, 12,  7,  9,  2,  1)   \\ 
IBP&  $0.450168  \pm 9.048 \, 10^{-3}$ &  (1,  9,  7,  9, 23,  1)   \\ 
 APGC &  $0.450146  \pm 9.056 \, 10^{-3}$ &  (4, 10,  9,  9,  3, 15)   \\ 
IBPG &  $0.450139  \pm 9.055 \, 10^{-3}$ &  (8,  9, 12,  8,  9,  4)   \\ 
\hline \end{tabular} 
 \end{center} 
 \end{table} 
 \end{center} 
We observe the following: 
\begin{itemize}

\item In both cases, IBPG-A and E-A-HALS have similar convergence rate, but IBPG-A converges to better solution than E-A-HALS more often.  These algorithms outperform the others. 

\item IBPG performs better than APGC in terms of final error obtained, while the convergence speeds are similar. 

\end{itemize}

\subsubsection*{Sparse document data sets} 

\begin{figure}[ht!]
\begin{center}
\begin{tabular}{ccc}
&\includegraphics[width=0.49\textwidth]{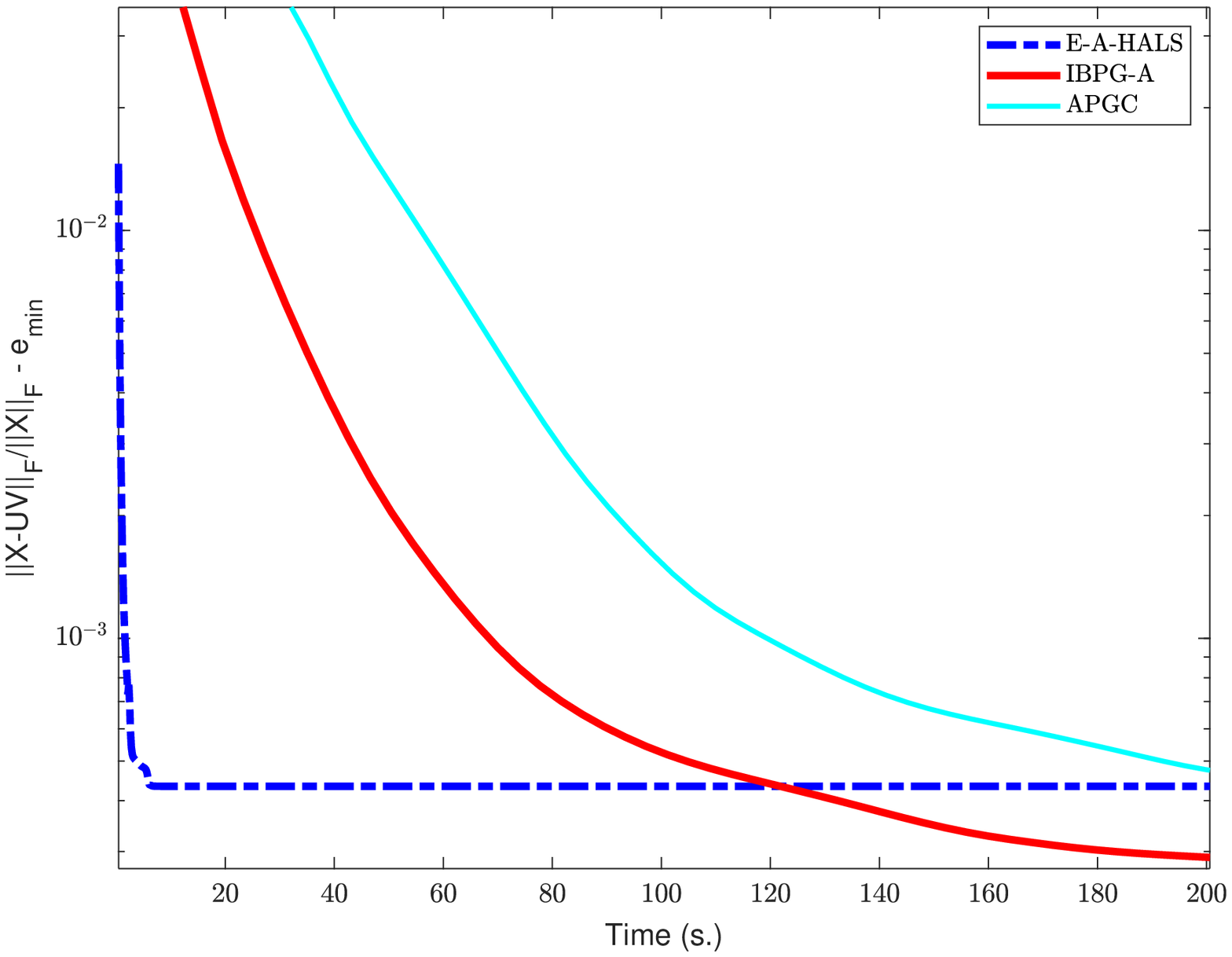}  & 
\includegraphics[width=0.49\textwidth]{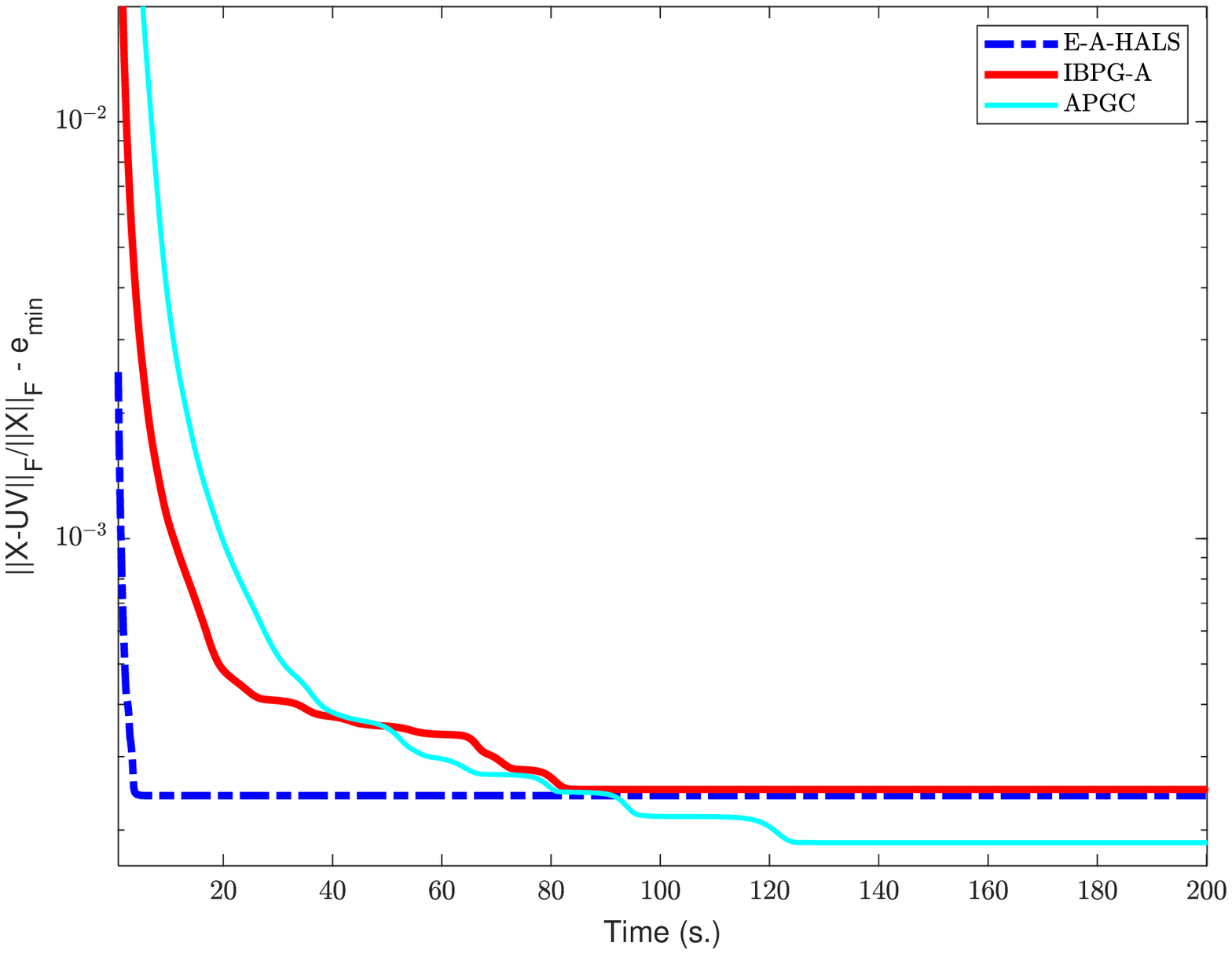}  \\
&\includegraphics[width=0.49\textwidth]{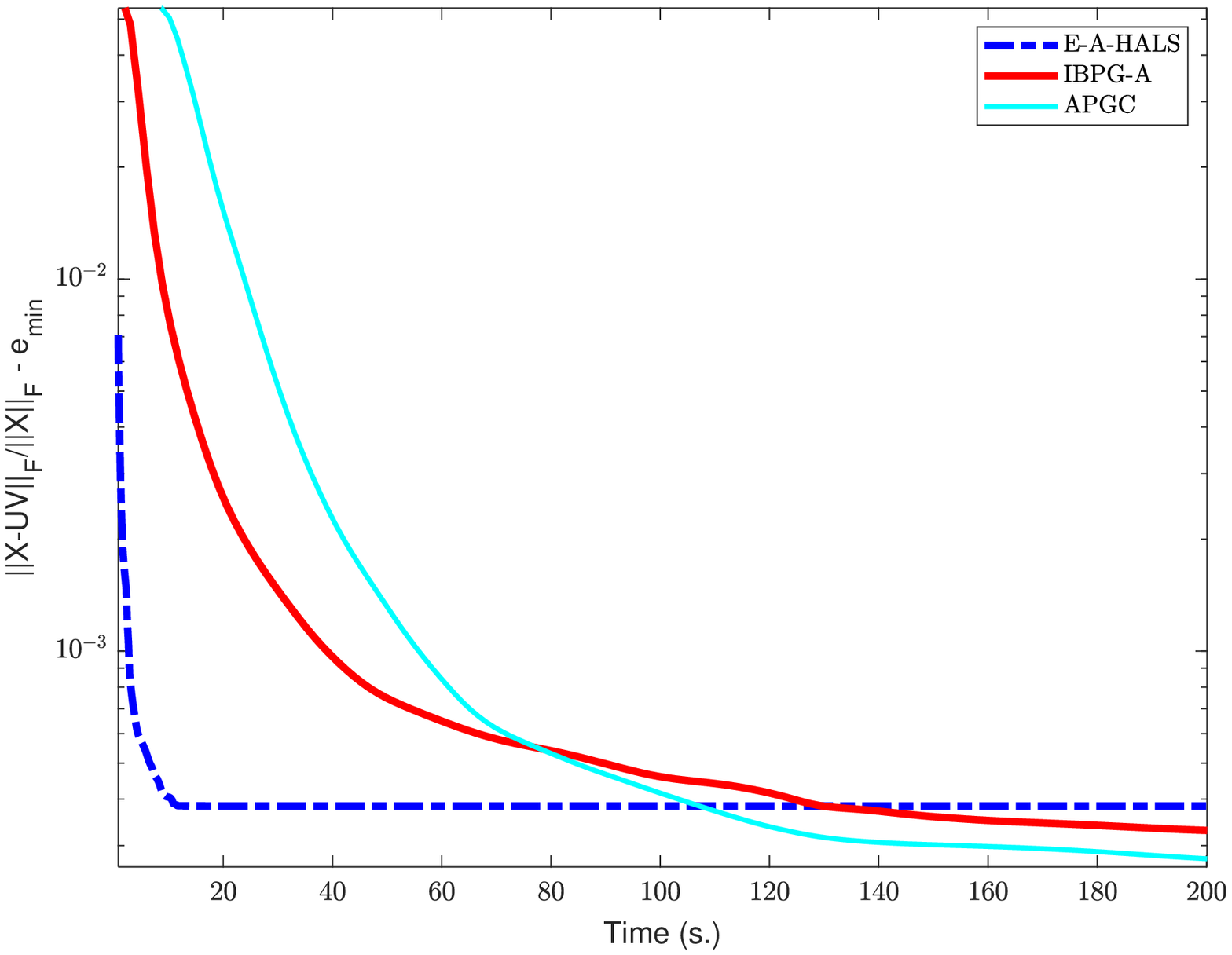}  &
\includegraphics[width=0.49\textwidth]{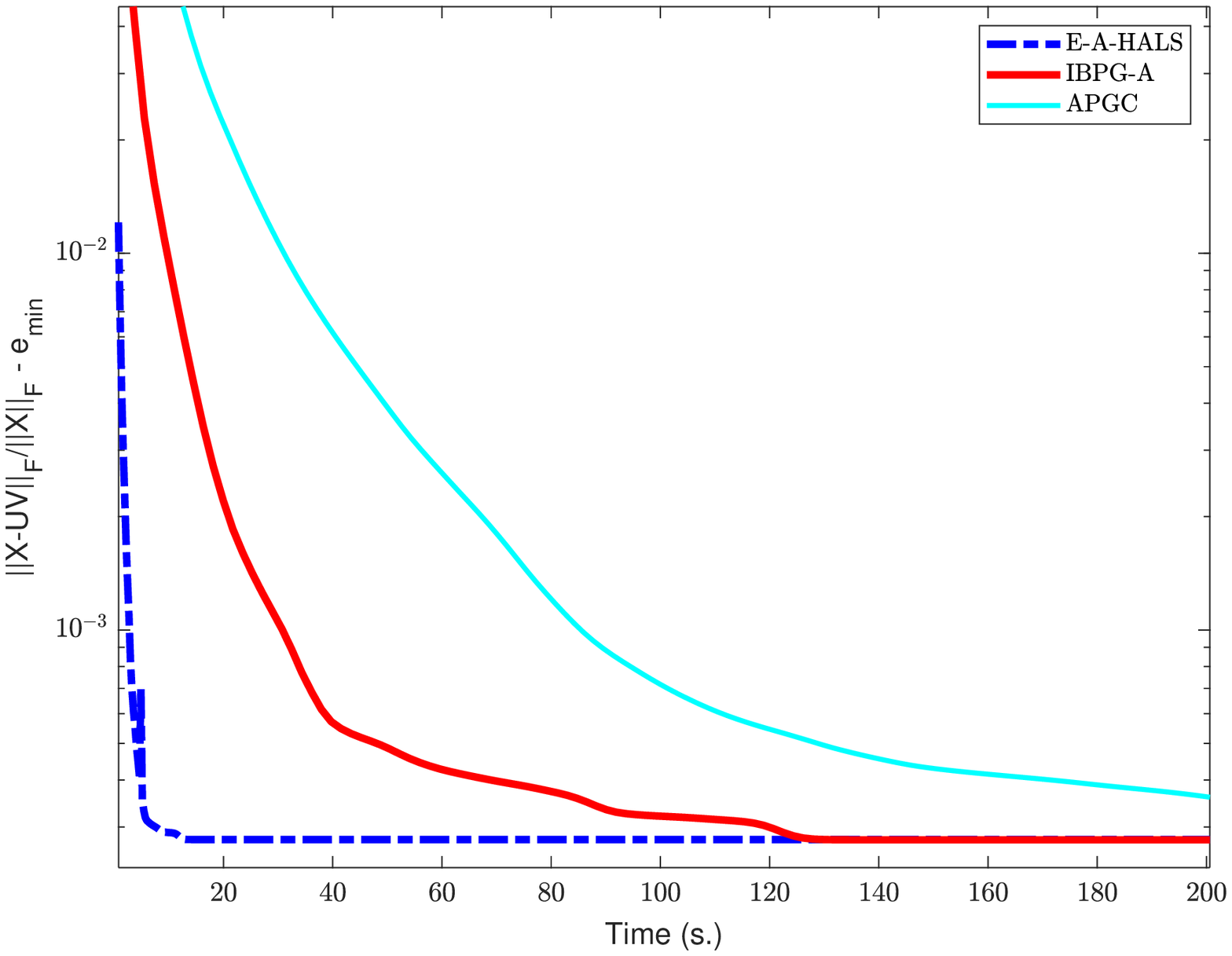}   \\
&\includegraphics[width=0.49\textwidth]{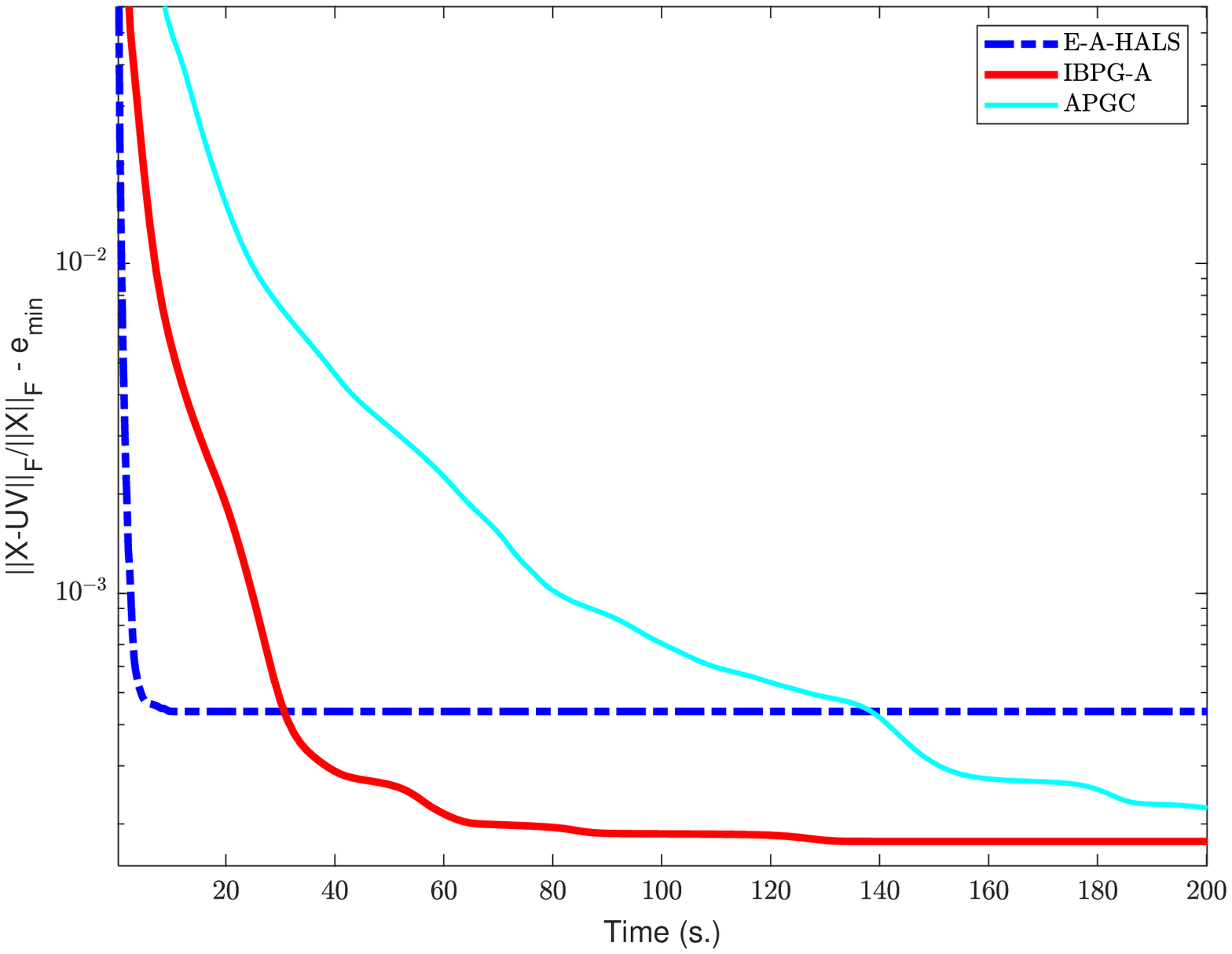} & 
\includegraphics[width=0.49\textwidth]{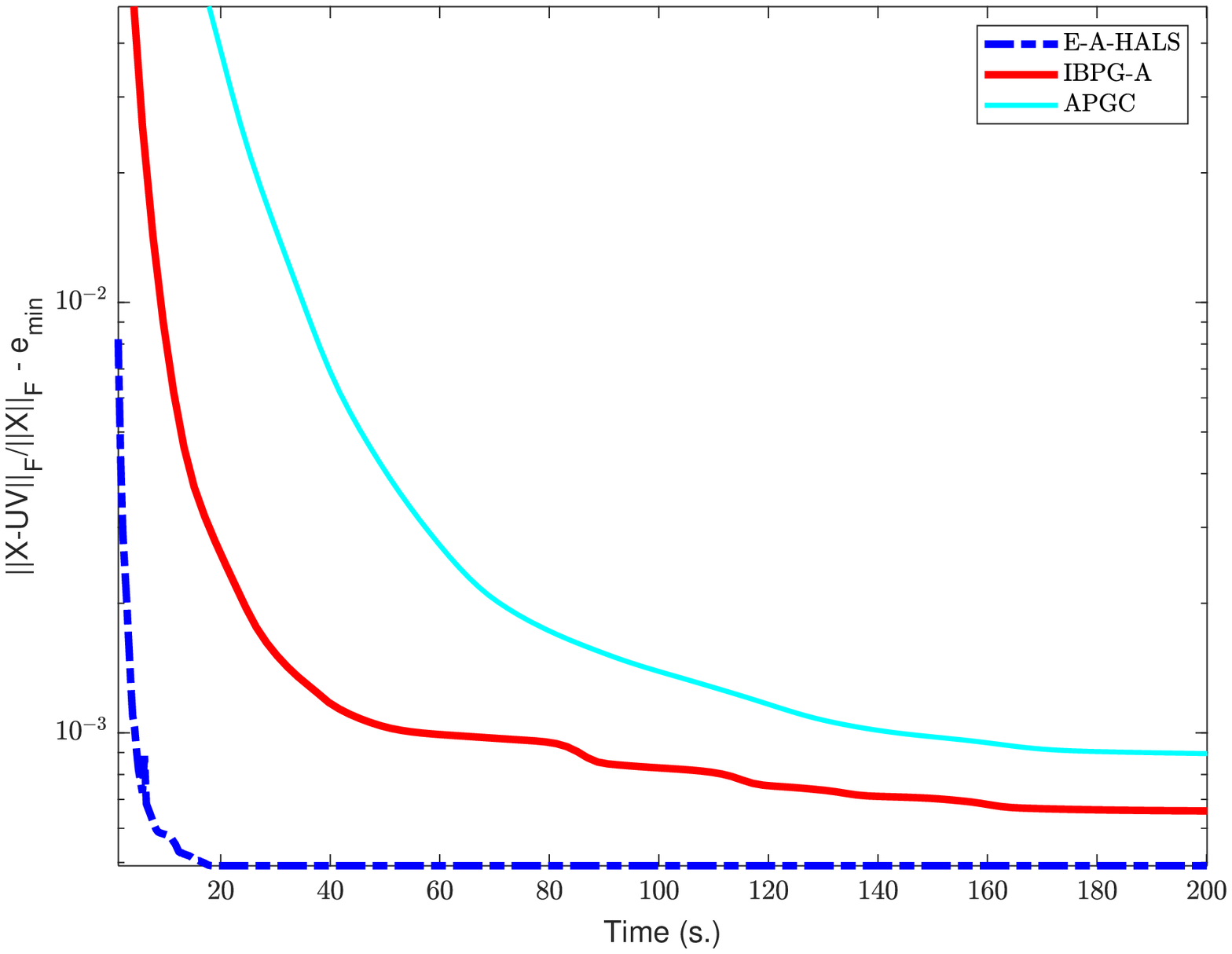} 
\end{tabular}
\caption{Average value of $E(k)$ with respect to time on 6 document data sets: 
Classic (top left), 
Hitech (top right), 
La1 (middle left), 
Ohscal (middle right), 
Reviews (bottom left) and 
Sports (bottom right).
\label{fig:textmining}} 
\end{center}
\end{figure}

We test the algorithms on the same six sparse document data sets with $r=10$ as in~\cite{Ang2018}.  
Figure~\ref{fig:textmining} reports the evolution of the average of $E(k)$ over 35 initializations, 
and Table~\ref{table:textmining} reports the average error, standard deviation and ranking of the final value of $E(k)$ among the 210 runs (6 data sets with 35 initializations for each data set). 

\begin{center}  
 \begin{table}[h!] 
 \begin{center} 
\caption{Average, standard deviation and ranking of the value of $E(k)$ at the last iteration among the different runs on the document data sets. 
The best performance is highlighted in bold. 
\label{table:textmining}}
\vspace{0.15cm}
 \begin{tabular}{|c|c|c|c|} 
 \hline Algorithm &  mean $\pm$ std & ranking  \\ 
 \hline 
E-A-HALS &  $0.881969  \pm 3.021\; 10^{-2}$ &  (73, 55, 82)   \\ 
IBPG-A  &  $\textbf{0.881921}  \pm 3.021\; 10^{-2}$ &  ({\textbf{87}}, 68, 55)   \\ 
 APGC &  $0.881992  \pm 3.019\; 10^{-2}$ &  (51, 86, 73)   \\ 
\hline
\end{tabular} 
 \end{center} 
 \end{table} 
 \end{center} 
 
 For these sparse datasets, E-A-HALS converges with the fastest rate, followed by IBPG-A.  However, IBPG-A generates in average the best final solutions. 
  
\subsection{Non-negative approximate canonical polyadic decomposition (NCPD)}
We consider the following NCPD problem: given a non-negative tensor ${T}\in \mathbb R^{I_1 \times I_2 \times\ldots\times I_N}$ and a specified order $\mathbf r$, solve 
\begin{equation}
\label{eq:NCPD}
\min_{{X}^{(1)},\ldots,{X}^{(N)}} f:=\frac12\norm{{T}-{X}^{(1)} \circ \ldots \circ {X}^{(N)}}_F^2 \quad\text{such that} \quad {X}^{(n)}\in \mathbb{R}_+^{I_n\times \mathbf{r}}, n=1,\ldots,N,
\end{equation}
where the Frobenius norm of a tensor $ {T}\in \mathbb R^{I_1 \times I_2 \times\ldots\times I_N}$ is defined as $\norm{{T}}_F=\sqrt{\sum_{i_1,\ldots,i_N}T_{i_1i_2\ldots i_N}^2}$, and the tensor product ${X}={X}^{(1)} \circ \ldots \circ {X}^{(N)}$ is defined as 
$$
X_{i_1i_2\ldots i_N}=\sum_{j=1}^{\mathbf{r}} X_{i_1j}^{(1)}  X_{i_2j}^{(2)} \ldots  X_{i_Nj}^{(N)}, \text{for} \, i_n\in \{1,\ldots,I_n\}, n=1,\ldots,N.
$$
Here $X_{ij}^{(n)}$ is the $(i,j)$-th element of ${X}^{(n)}$. 
Let us denote
\begin{eqnarray}
B^{(i)}=X^{(N)} \odot  \cdots  \odot X^{(i+1)} \odot X^{(i-1)}\odot \cdots \odot X^{(1)},
\label{def:B}
\end{eqnarray}
where  $\odot$ is the Khatri-rao product.
Then the gradient of $f$ with respect to $X^{(i)}$ is
\begin{equation}
\label{def:grad}
\nabla_{X^{(i)}} f=\Big( X^{(i)} \big( B^{(i)}\big)^T  - T_{[i]}\Big) B^{(i)}, 
\end{equation}
where $T_{[i]} $ is the mode-$i$ matricization of $T$. We see that the gradient $\nabla_{X^{(i)}} f$ is Lipschitz continuous with the constant $
L^{(i)}=\|\big(B^{(i)}\big)^T B^{(i)}\|
$,
where $B^{(i)}$ is defined in (\ref{def:B}) and $\|\cdot\|$ is the operator norm. 

As for NMF, we can write \eqref{eq:NCPD} as a problem of the form \eqref{eq:main} with $s=N$ variables. We then can apply IBPG for NCPD and use the same extrapolation parameters as NMF, see Section \ref{sec:NMF-subconvergence}. Denote 
$B_{k-1}^{(p)}=X^{(N)}_{k-1} \odot  \cdots  \odot X_{k-1}^{(p+1)} \odot
X_{k}^{(p-1)}\odot \cdots \odot X_{k}^{(1)}$ and
$L_{k}^{(p)}=\norm{\big(B_{k}^{(p)}\big)^T B_{k}^{(p)}}$. Algorithm \ref{alg:IBPG_NCPD} describes the pseudo code of IBPG when applied for solving the NCPD problem \eqref{eq:NCPD}. Step 5 of Algorithm \ref{alg:IBPG_NCPD} indicates that we cyclically update the factors $X^{(i)}$. Note that IBPG described in Algorithm \ref{alg:Fro} allows to randomly select one factor among the $N$ factors to update as long as all factors are updated  after $T_k$ iterations, leading to other variants of IBPG when applied to solve the NCPD problem. In our experiments, we implement  Algorithm \ref{alg:IBPG_NCPD}, which is the cyclic update version.  

\begin{algorithm}[tb]
	\caption{IBPG for NCPD \label{alg:IBPG_NCPD} }
	\begin{algorithmic}[1]
		\STATE Initialization: Choose $\delta_w=0.99$, $\beta=1.01 $, $t_0=1$, and initial factor matrices $\big(X^{(1)}_{-1}, \ldots, X^{(N)}_{-1}\big )=\big(X^{(1)}_{0}, \ldots, X^{(N)}_{0}\big )$. Set $k=1$.
		\STATE Set $X_{\rm{pr}}^{(i)}= X^{(i)}_{-1}$, $i=1,\ldots,N$. \% \textit{$X_{\rm{pr}}^{(i)}$ is to save the previous value of block $i$. }
		\STATE Set $X_{\rm{cr}}^{(i)}= X^{(i)}_{0}$, $i=1,\ldots,N$. \% \textit{$A_{\rm{cr}}^{(i)}$ is to save the current value of block $i$. }
		\REPEAT
		\FOR{ $i=1,\ldots,N$ }
		\STATE Compute  $t_k=\frac12\big( 1 + \sqrt{1 + 4 t_{k-1}^2}\big)$, $\hat w_{k-1}=\frac{t_{k-1}-1}{t_k}$ and  \vspace*{-2mm}
		$$w_{k-1}^{(i)}=\min \left( \hat w_{k-1}, \delta_w \sqrt{\frac{L_{k-2}^{(i)}}{L_{k-1}^{(i)}}}\right). $$  \vspace*{-3mm}
		\REPEAT
		\STATE Compute two extrapolation points  \vspace*{-2mm}
		$$\hat  X^{(i,1)}= X_{\rm{cr}}^{(i)} + w_{k-1}^{(i)} \Big( X_{\rm{cr}}^{(i)}-X_{\rm{pr}}^{(i)} \Big),$$ \vspace*{-2mm}
		and
		$$\hat X^{(i,2)}= X_{\rm{cr}}^{(i)} +\beta w_{k-1}^{(i)} \Big( X_{\rm{cr}}^{(i)}-X_{\rm{pr}}^{(i)} \Big)$$ \vspace*{-2mm}
		\STATE Set $X_{\rm{pr}}^{(i)}=X_{\rm{cr}}^{(i)}$.
		\STATE Update $X_{\rm{cr}}^{(i)}$ by projected gradient step:
		\begin{equation*}
		X_{\rm{cr}}^{(i)}=\max \left(0,\hat  X^{(i,2)}- \frac{1}{L_{k-1}^{(i)}}  \Big (\hat X^{(i,1)} \big(B_{k-1}^{(i)}\big)^T- \mathcal T_{[i]}\Big) B_{k-1}^{(i)}\right). \vspace*{-2mm}
		\end{equation*}
		\UNTIL{some criteria is satisfied}
		\STATE Set $X_{k}^{(i)}=X_{\rm{cr}}^{(i)}$.
		\ENDFOR
		\STATE Set $k=k+1$.
		\UNTIL{some criteria is satisfied}
	\end{algorithmic}
\end{algorithm}

In the following, we consider the three-way NCPD problem, i.e., 
\begin{equation*}
\min_{U,V,W} f:=\norm{{T}-U \circ V \circ W}_F^2 \quad \text{such that} \quad U\in \mathbb{R}_+^{I\times \mathbf{r}}, V\in \mathbb{R}_+^{J\times \mathbf{r}}, W\in \mathbb{R}_+^{K\times \mathbf{r}}.
\end{equation*}
and compare our algorithm IBPG-A (i.e., Step 8--10 of Algorithm \ref{alg:IBPG_NCPD} are repeated several times for updating each factor $U$, $V$, $W$ before doing so for the next factor) with  APGC  \cite{Xu2013}, A-HALS \cite{cichocki2009nonnegative,Gillis2012}, and  ADMM \cite{Huang2016flexible}. We note that it was empirically showed in \citep[Section 4.2]{Xu2013} that APGC outperforms the two ANLS based methods that are ANLS active set methods and ANLS block pivot methods.

Similarly to the NMF experiments, we define the relative errors ${\rm relerror}_k=\frac{\norm{X- U^{(k)} \circ V^{(k)} \circ W^{(k)}}_F}{\norm{X}_F}$ 
and $E(k)= {\rm relerror}_k - e_{\min}$. 

\subsubsection{Experiments with synthetic data sets} 
\begin{figure}[ht!]
\begin{center}
\begin{tabular}{ccc}
&\includegraphics[width=0.46\textwidth]{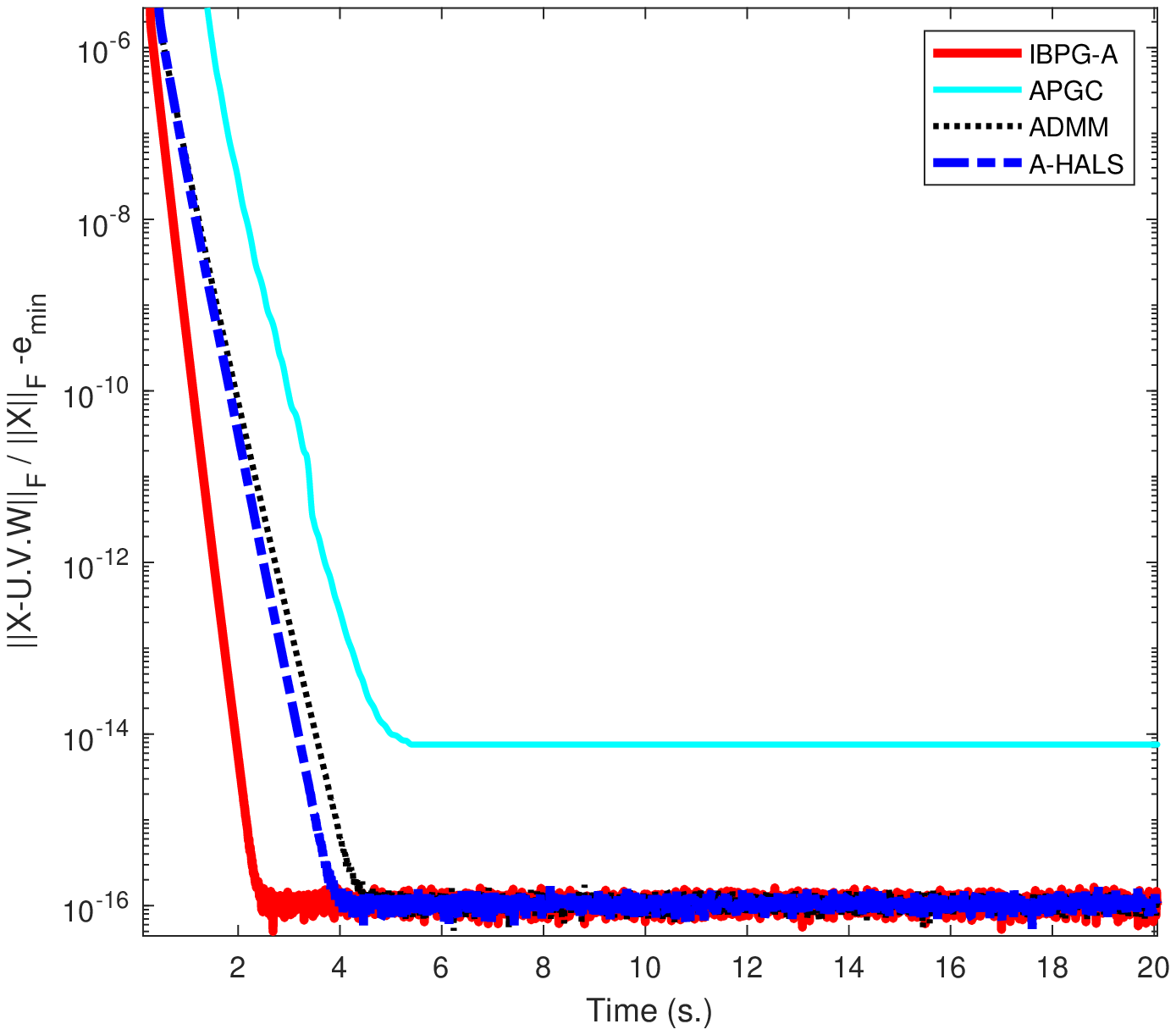}  & 
\includegraphics[width=0.46\textwidth]{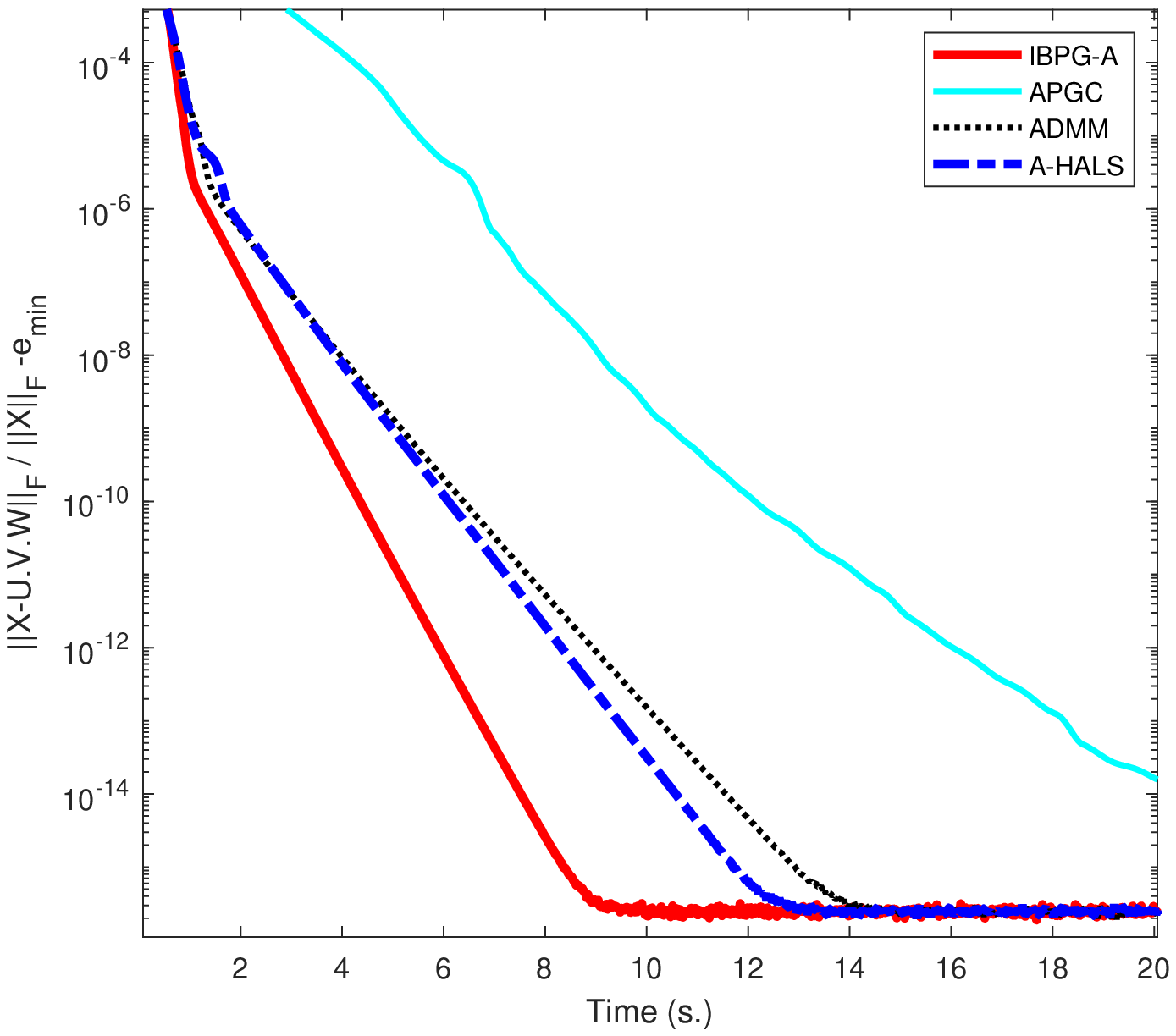}  \\
&\includegraphics[width=0.46\textwidth]{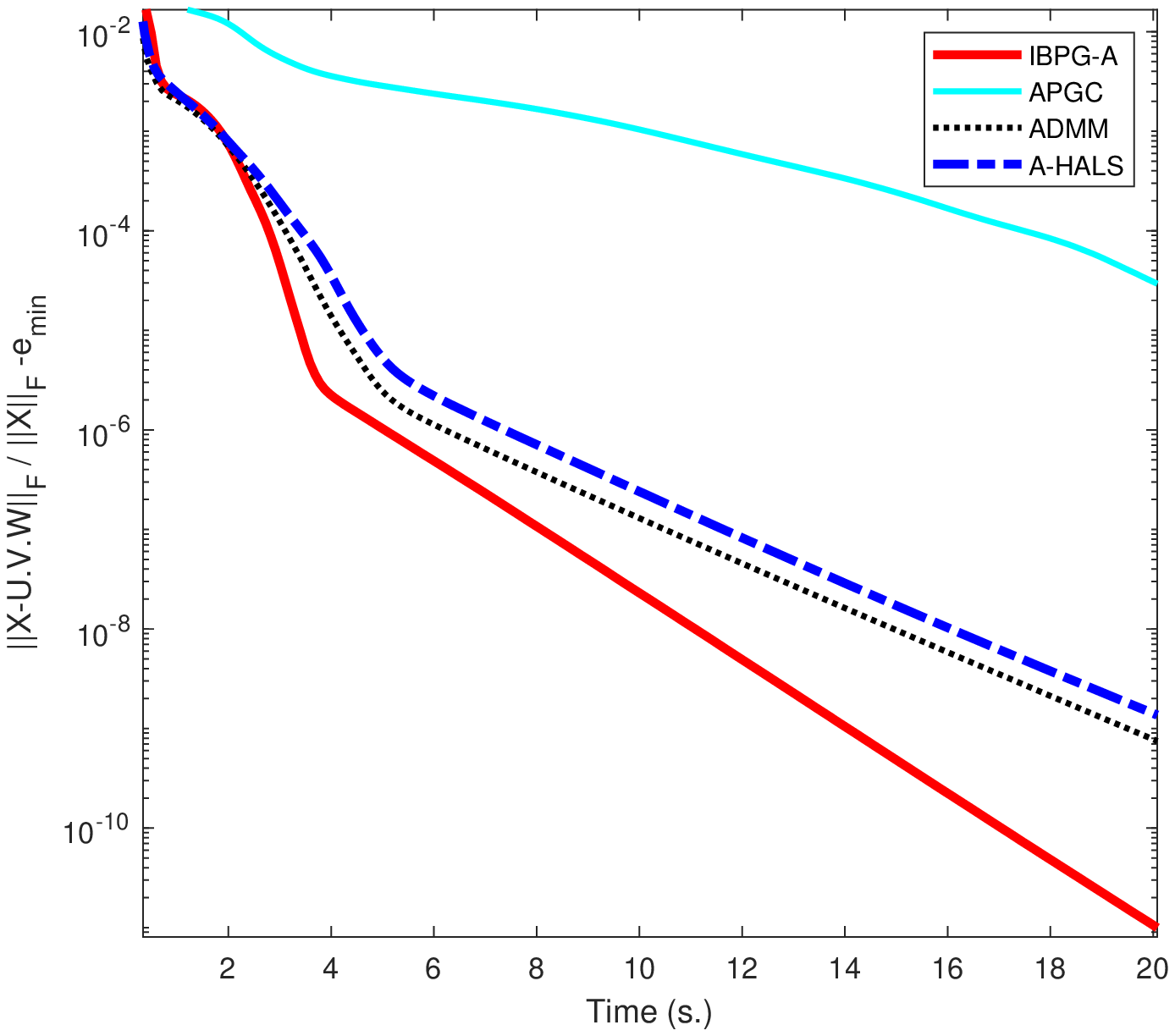}  &
\includegraphics[width=0.46\textwidth]{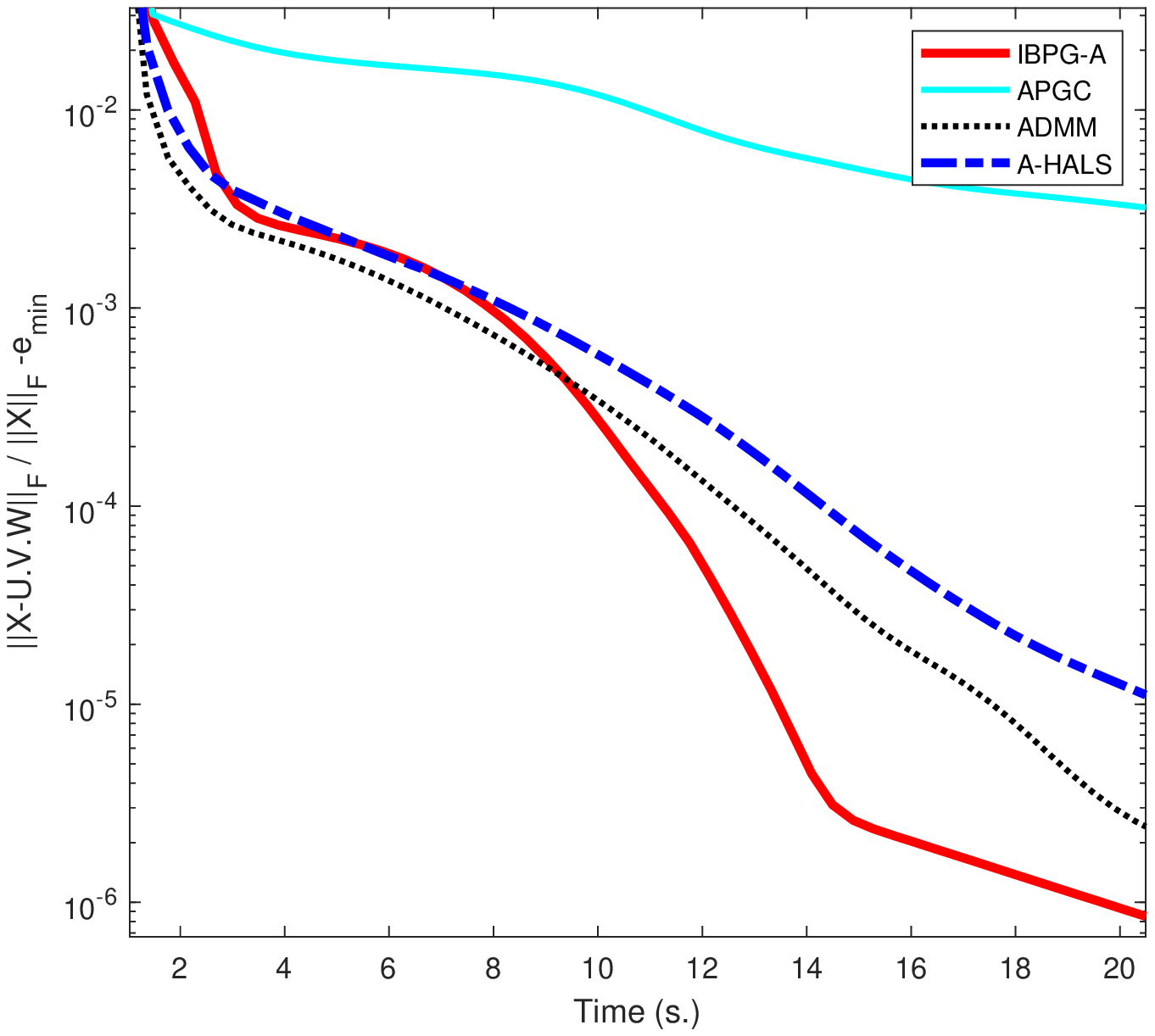}  
\end{tabular}
\caption{Average value of $E(k)$ with respect to time on 4 synthetic data sets: 
$100 \times 100 \times 100$ (top left), 
$100 \times 100 \times 500$ (top right), 
$100 \times 500 \times 500$ (bottom left) and 
$500 \times 500 \times 500$  (bottom right).
\label{fig:NCPD synthetic}} 
\end{center}
\end{figure}

Four three-way tensor of size $100 \times 100 \times 100$, $ 100 \times 100 \times 500$, $100 \times 500 \times 500$ and $500 \times 500 \times 500$ are generated by letting $X=U \circ V \circ V$, where $U$, $V$ and $W$ are generated by commands $rand(I,\mathbf r)$,  $rand(J, \mathbf r)$ and $rand(K,\mathbf r)$ with $\mathbf r=20$. For each $X$ we run all algorithms with the same 50 random initializations $U_0=rand(I,\mathbf r)$,  $V_0=rand(J, \mathbf r)$ and $W_0=rand(K,\mathbf r)$. For each initialization we run each algorithm for 20 seconds. Figure \ref{fig:NCPD synthetic} illustrates the evolution of the average of $E(k)$ over 50 initializations with respect to time. 

We observe that IBPG-A converges with the fastest rate followed by ADMM and A-HALS.

To compare the accuracy of the solutions, we generate 50 random $I \times J \times K$ tensors with $I$, $J$, $K$  being random integer numbers in the interval [100, 500]. For each tensor we run the algorithms for 20 seconds with 1 random initialization. Table \ref{table:NCPD_acc} reports the average, the standard deviation of the errors and the ranking vector of the relative errors. 

\begin{center}  
 \begin{table}[h!] 
 \begin{center} 
\caption{Average, standard deviation and ranking of the value of $E(k)$ at the last iteration among the different runs on the synthetic data sets. 
The best performance is highlighted in bold \label{table:NCPD_acc}} 
\vspace{0.1in}
 \begin{tabular}{|c|c|c|c|} 
 \hline Algorithm &  mean $\pm$ std & ranking  \\ 
 \hline 
IBPG-A &  $\mathbf{3.398 \, 10^{-9}} \pm 2.010 \, 10^{-8}$ &  (\textbf{49},  1,  0,  0)   \\ 
APGC &  $1.505 \,10^{-4} \pm 3.548 \, 10^{-4}$ &  ( 1,  0,  0, 49)   \\ 
ADMM  &  $1.108 \, 10^{-8} \pm 3.513 \, 10^{-8}$ &  (18, 22, 10,  0)   \\ 
 A-HALS &  $1.734 \, 10^{-8} \pm 5.038 \, 10^{-8}$ &  (19, 10, 21,  0)   \\ 
\hline 
\end{tabular} 
 \end{center} 
  \vskip -0.1in
 \end{table} 
 \end{center} 
 We can see that IBPG-A outperforms the other algorithms in term of the accuracy of the solutions. 
\subsubsection{Experiments with real data sets} 
We test the algorithms on two real data sets that are CBCL -- a face image data set \cite{Shashua2005} and Urban -- a hyperspectral image data set \cite{gillis2015hierarchical}. 

\begin{figure}[ht!]
\begin{center}
\begin{tabular}{cc}
&\includegraphics[width=0.46\textwidth]{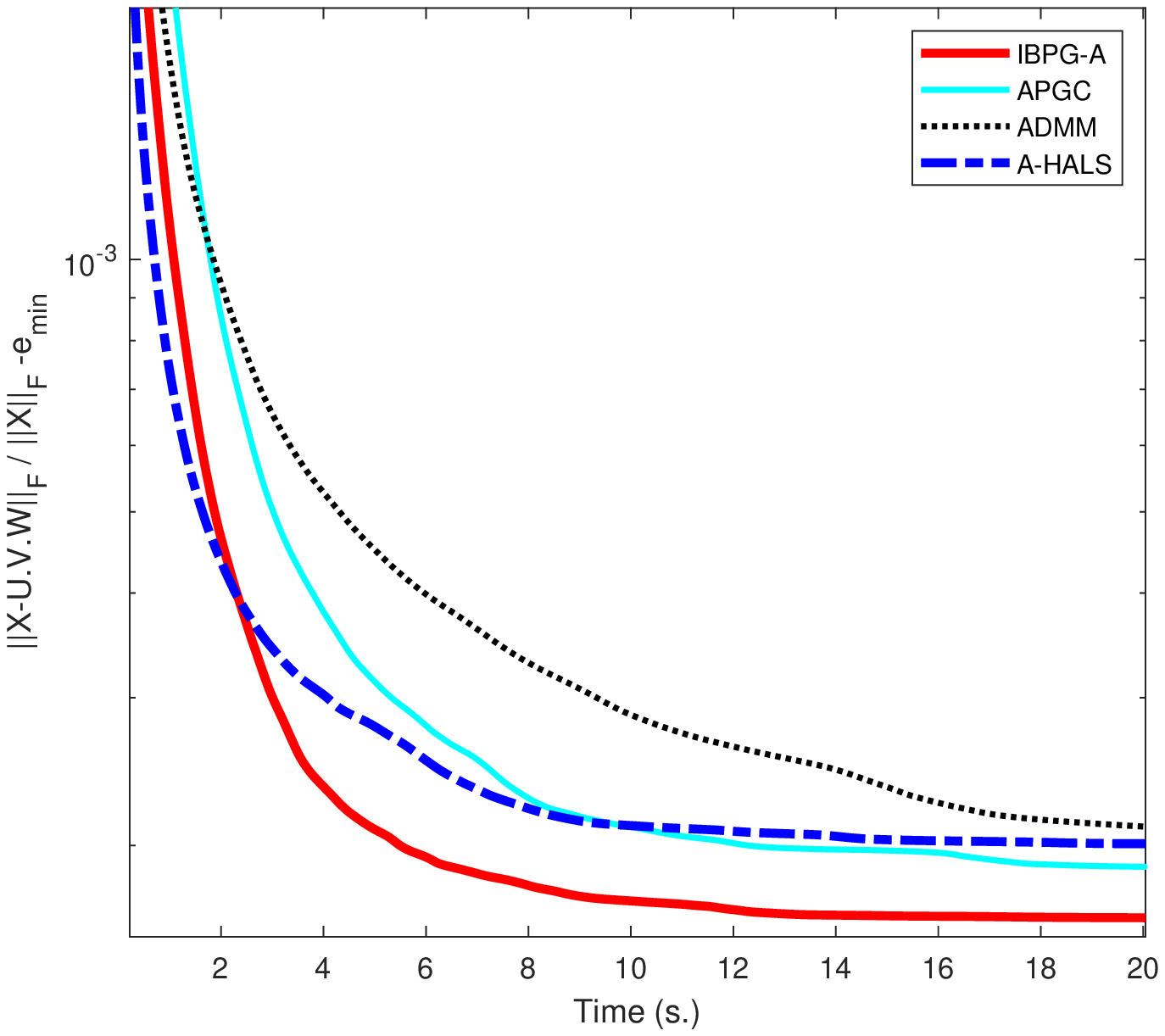}   
\includegraphics[width=0.46\textwidth]{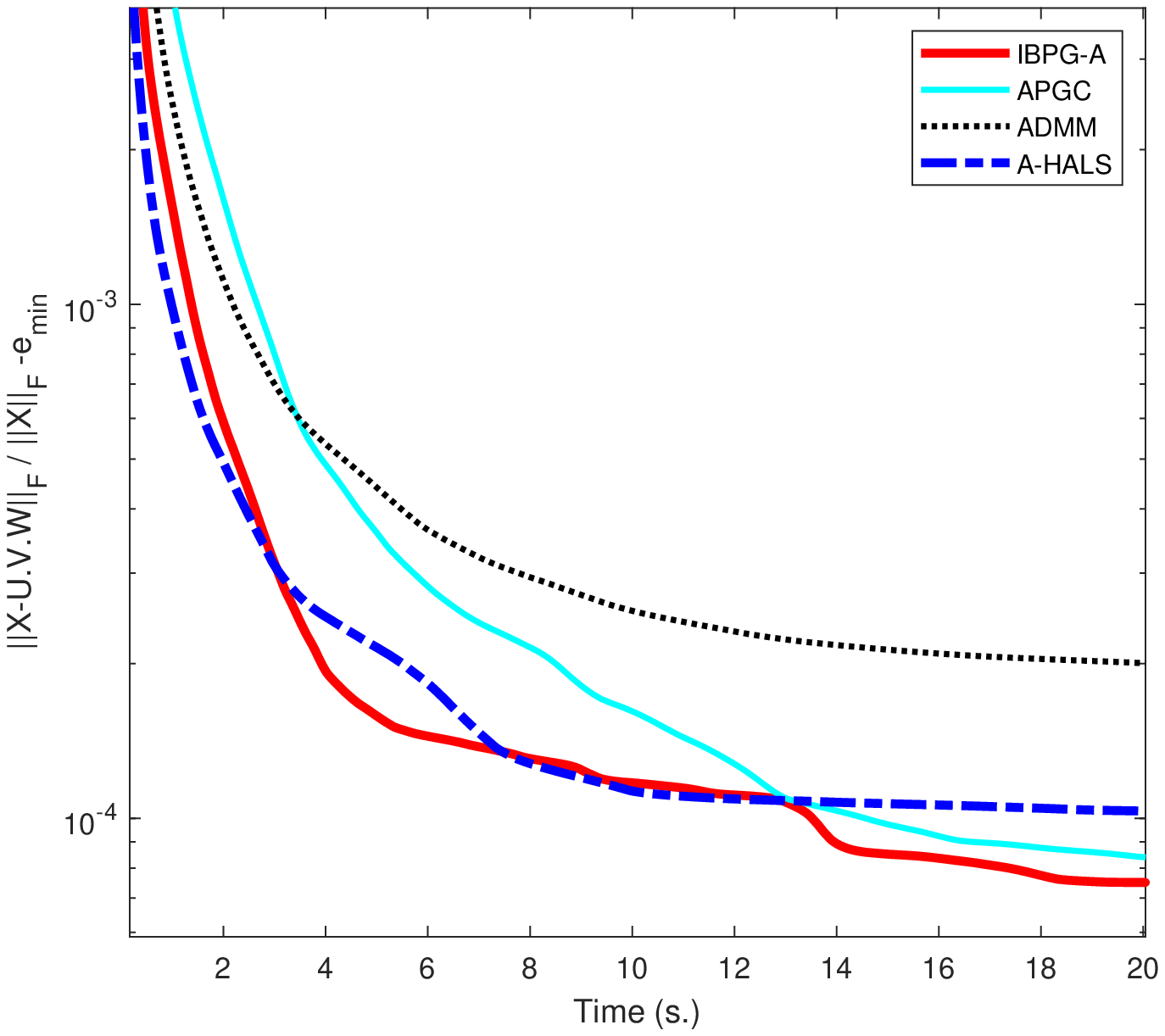}  \\
&\includegraphics[width=0.46\textwidth]{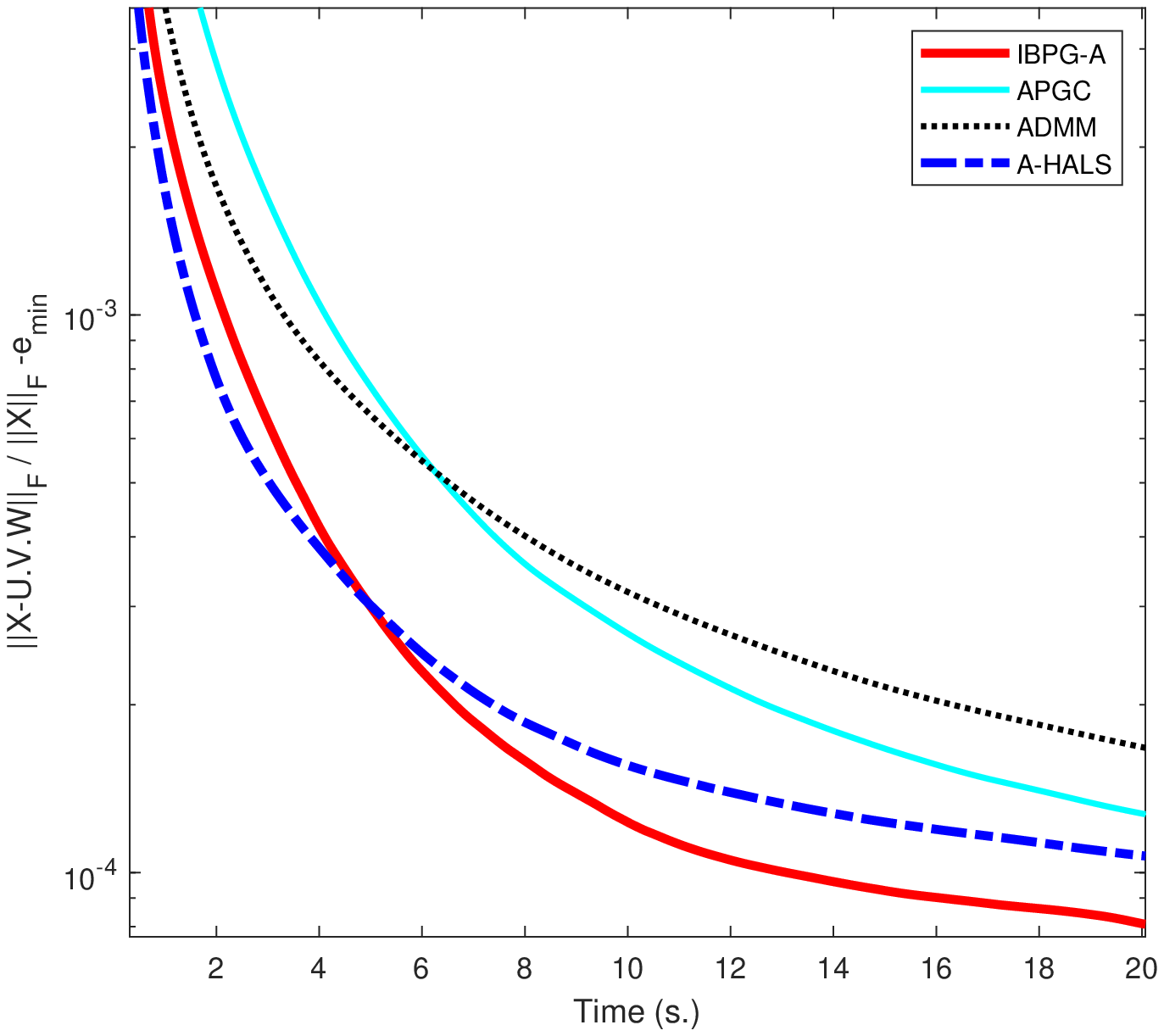}  
\end{tabular}
\caption{Average value of $E(k)$ with respect to time on the CBCL data set with different rank:  
$r=20$ (above left), 
$r=30$ (above right) and 
$r=40$ (below).
\label{fig:NCPD_CBCL}} 
\end{center}
\end{figure}

The CBCL data set can be considered as a three-way tensor of the size $19 \times 19 \times 2429$ that is formed from 2429 of $19 \times 19$ images from the MIT CBCL database. We test different ranks $r=20$, $r=30$ and $r=40$. For each $r$, we run all algorithms with the same 50 random initializations and for each initialization we run each algorithm for 20 seconds. Figure \ref{fig:NCPD_CBCL} shows the evolution of $E(k)$ over 50 initializations, and Table \ref{table:NCPD_CBCL} reports the average error, standard deviation and ranking of the final value of $E(k)$ among 150 runs (three values of $r$ with 50 initializations for each $r$). 

\begin{center}  
 \begin{table}[h!] 
 \begin{center} 
\caption{Average, standard deviation and ranking of the value of $E(k)$ at the last iteration among the different runs on the CBCL data set. 
The best performance is highlighted in bold
 \label{table:NCPD_CBCL}} 
 \vspace{0.1in}
 \begin{tabular}{|c|c|c|c|} 
 \hline Algorithm &  mean $\pm$ std & ranking  \\ 
 \hline 
IBPG-A &  $\mathbf{0.013751} \pm 3.4019 \,10^{-3}$ &  (\textbf{58}, 46, 29, 17)   \\ 
APGC &  $0.013772 \pm 3.4016 \,10^{-3}$ &  (37, 36, 49, 28)   \\ 
ADMM  &  $0.013799 \pm 3.3932 \,10^{-3}$ &  (21, 27, 39, 63)   \\ 
 A-HALS &  $0.013778 \pm 3.4078 \,10^{-3}$ &  (34, 41, 33, 42)   \\  
\hline 
\end{tabular} 
 \end{center} 
  \vskip -0.1in
 \end{table} 
 \end{center} 

\begin{figure}[ht!]
\begin{center}
\begin{tabular}{cc}
&\includegraphics[width=0.46\textwidth]{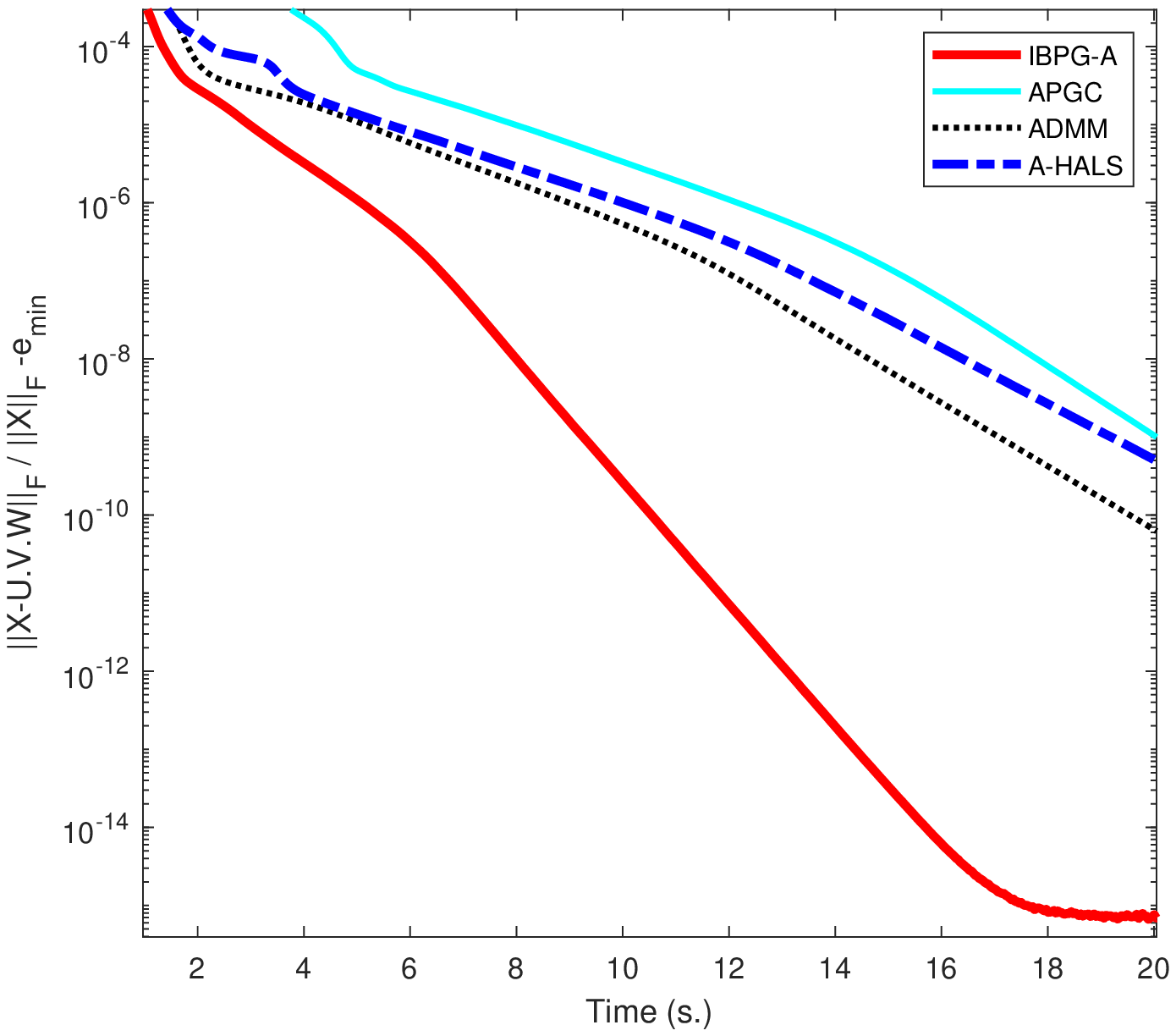}   
\includegraphics[width=0.46\textwidth]{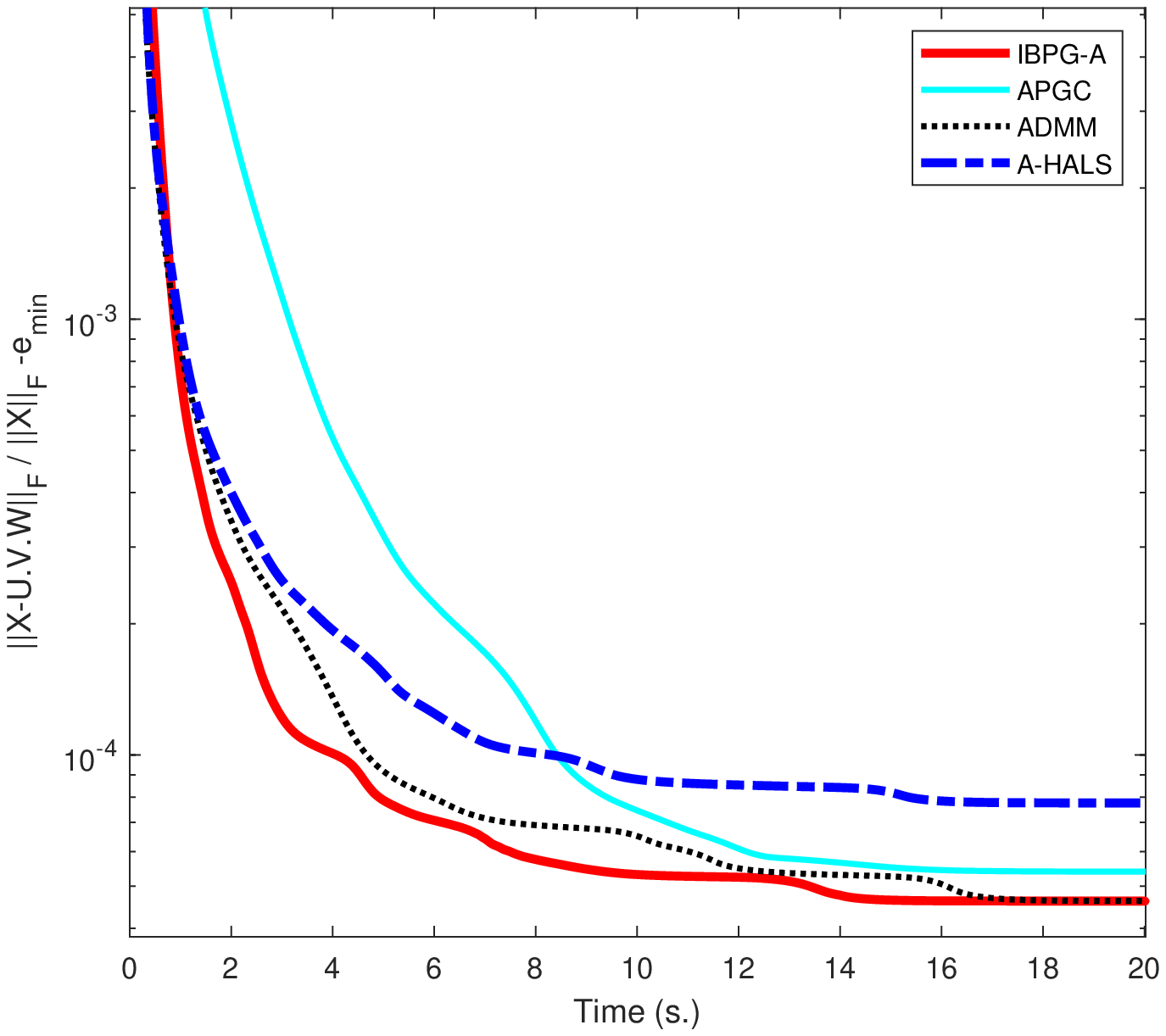}  \\
&\includegraphics[width=0.46\textwidth]{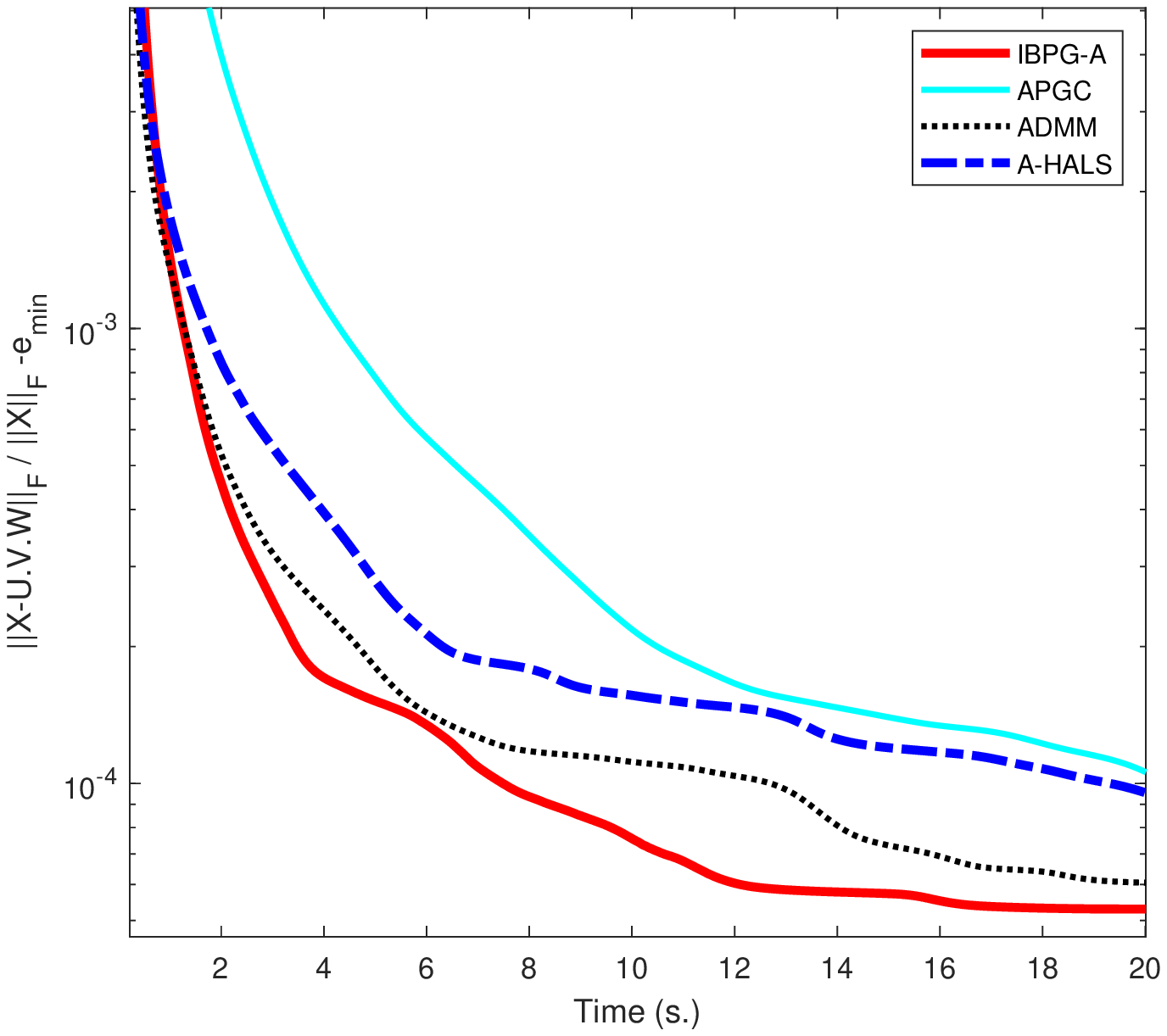}  
\end{tabular}
\caption{Average value of $E(k)$ with respect to time on the Urban data set with different rank:  
$r=6$ (above left), 
$r=8$ (above right) and 
$r=10$ (below).
\label{fig:NCPD_Urban}} 
\end{center}
\end{figure}

The Urban data set can be considered as a three-way tensor of the size $307 \times 307 \times 162$ that is formed from 162 channels of $307 \times 307$ pixels. We test different ranks $r=6$, $r=8$ and $r=10$. For each $r$, we run all algorithms with the same 50 random initializations and for each initialization we run each algorithm for 20 seconds. Figure \ref{fig:NCPD_Urban} shows the evolution of $E(k)$ over 50 initializations, and Table \ref{table:NCPD_Urban} reports the average error, standard deviation and ranking of the final value of $E(k)$ among 150 runs (three values of $r$ with 50 initializations for each $r$).
\begin{center}  
 \begin{table}[h!] 
 \begin{center} 
\caption{Average, standard deviation and ranking of the value of $E(k)$ at the last iteration among the different runs on the Urban data set. 
The best performance is highlighted in bold
\label{table:NCPD_Urban}} 
\vspace{0.1in} 
 \begin{tabular}{|c|c|c|c|} 
 \hline Algorithm &  mean $\pm$ std & ranking  \\ 
 \hline 
IBPG-A &  $\mathbf{0.130624 } \pm 7.5499\,10^{-3}$ &  (\textbf{107}, 20, 13, 10)   \\ 
APGC &  $0.130630 \pm 7.5461\,10^{-3}$ &  ( 6, 34, 56, 54)   \\ 
ADMM  &  $0.130627 \pm 7.5470\,10^{-3}$ &  (31, 50, 45, 24)   \\ 
 A-HALS &  $0.130642 \pm 7.5410\,10^{-3}$ &  (12, 40, 36, 62)   \\ 
\hline 
\end{tabular} 
 \end{center} 
  \vskip -0.1in
 \end{table} 
 \end{center} 
We can observe that IBPG-A outperforms the other algorithms both in terms of convergence speed and accuracy.
\end{document}